\documentclass[journal]{new-aiaa}
\usepackage[utf8]{inputenc}
\usepackage{textcomp}

\usepackage{graphicx}
\usepackage{longtable,tabularx}
\usepackage{amsmath}
\usepackage{amssymb}
\usepackage{mathtools}
\usepackage{bm}
\usepackage{dsfont}
\usepackage{amsthm}

\usepackage{multirow}
\usepackage{subcaption}
\usepackage{algorithm}
\usepackage{algpseudocode}

\usepackage{xcolor}
\usepackage{tikz}
\usetikzlibrary{automata,arrows,positioning}

\newtheorem{lemma}{Lemma}
\newtheorem{problem}{Problem}

\newtheorem{remark}{Remark}
\newtheorem{assumption}{Assumption}

\newlength{\clabelwidth}
\newcommand{\clabel}[1]{\mathllap{\makebox[\clabelwidth][l]{$\scriptstyle\blacktriangleleft$\;{\itshape\small #1}}}}

\usepackage{cleveref}
\crefname{lemma}{Lemma}{Lemmas}
\crefname{problem}{Problem}{Problems}
\crefname{definition}{Definition}{def}
\crefname{remark}{Remark}{rem}
\crefname{assumption}{Assumption}{Assumptions}
\Crefname{assumption}{Assumption}{Assumptions}
\crefname{section}{section}{sections}
\Crefname{section}{Section}{Sections}
\crefname{subsection}{subsection}{subsections}
\Crefname{subsection}{Subsection}{Subsections}
\crefname{subsubsection}{subsubsection}{subsubsections}
\Crefname{subsubsection}{Subsubsection}{Subsubsections}

\newcommand{\blfootnote}[1]{%
  \begingroup
  \renewcommand{\thefootnote}{}\footnote{#1}%
  \addtocounter{footnote}{-1}%
  \endgroup
}

\title{MASTO: Missed-thrust Aware Stochastic Trajectory Optimization}

\author{Shaurya Shrivastava\footnote{Ph.D. Student, School of Aeronautics and Astronautics, 701 W. Stadium Ave., West Lafayette, IN 47907; corresponding author, shriva15@purdue.edu.} and Kenshiro Oguri\footnote{Assistant Professor, School of Aeronautics and Astronautics, 701 W. Stadium Ave., West Lafayette, IN 47907.}}
\affil{Purdue University, West Lafayette, IN, 47907}

\begin{document}

\blfootnote{An earlier version of this paper was presented as paper 26-907 at the 2026 AAS/AIAA Astrodynamics Specialist Conference, Whistler, BC.}
\maketitle

\begin{abstract}

Missed Thrust Events (MTEs) are a major roadblock towards widespread adoption of low-thrust propulsion for deep space missions. Unplanned thruster outages have previously been difficult to directly incorporate in the preliminary mission design process, and current approaches mostly introduce robustness by conservative heuristics such as mandatory coasting arcs, which do not possess any safety guarantees, and usually only consider one missed-thrust event. Therefore, in this work, we utilize a stochastic hybrid system framework to model a spacecraft as a Markov jump system that can switch between nominal and failure mode. We formulate the Missed-thrust resilient trajectory optimization problem as a Markov jump covariance steering problem, and propose a sequential convex programming-based approach, which can also respect chance constraints. We demonstrate our framework on two complex three-dimensional transfers in the Earth-Moon circular restricted three body problem, showcasing robustness against a range of MTEs, while respecting chance constraints. We further analyze controller performance across a range of MTE-probability distributions, offering valuable insights into the MTE-resilient trajectory design trade space.
\end{abstract} 

\section{Introduction}

\lettrine{A}{s} the demand for deep space missions continues to increase rapidly, there is a growing need for more autonomous operations\cite{starekSpacecraftAutonomyChallenges2016}, which significantly reduce the reliance on traditional human-in-the-loop control strategies. This necessitates the need for spacecraft control policies that are robust to uncertainties. Missed Thrust events (MTEs) represent one of the most critical uncertainties\cite{mccartyMissedThrustAnalysis2020}, which can seriously undermine mission success, or even risk complete failure. This is exacerbated by the use of low-thrust propulsion systems, which are increasingly popular for their efficiency\cite{levTechnologicalCommercialExpansion2019}, but their longer duty cycles make them more susceptible to missed thrust events. MTEs are quite common in deep space missions, with recent examples such as the Dawn spacecraft experiencing more than 15 MTEs\cite{mccartyMissedThrustAnalysis2020}. Current mission design practices to account for missed thrust events are mostly ad hoc approaches, such as mandatory coasting arcs\cite{grebowDawnSafingApproach2016} and reduced duty cycles\cite{ohAnalysisSystemMargins}, which are effective operational strategies, but can be conservative and typically require extensive Monte Carlo simulations.

Several approaches have been developed to design trajectories that are more robust against MTEs. Early works focused on studying various system margins \cite{ohAnalysisSystemMargins,laipertAutomatedMissedThrustPropellant2015} such as propellant and duty cycle margins for given nominal trajectories. To incorporate such margins into mission design, missed thrust recovery margin (MTRM) has emerged as a popular metric that can help analyze and robustify trajectories. It is defined as the longest outage duration from which a spacecraft can still recover and reach its target\cite{raymanDesignFirstInterplanetary2002}. Reference~\citenum{olympioDESIGNINGROBUSTLOWTHRUST} used a stochastic optimization method for increased robustness against a single thruster outage, and related the terminal chance constraint to MTRM. The MTRM was also directly constrained in a deterministic trajectory optimization framework\cite{olympioDeterministicMethodSpace}. Later, "ghost" trajectory design\cite{mccartyMissedThrustAnalysis2020} and virtual swarm optimization techniques\cite{venigallaMultiObjectiveLowThrustTrajectory2022} were proposed to directly optimize the MTRM along the trajectory. Recent works have also explored bilevel approaches\cite{sinhaBiLevelOptimalControl2025}, more efficient initial guess generation\cite{sinhaInitialGuessGeneration2025a}, learning-based approaches\cite{rubinsztejnNeuralNetworkOptimal2020,zavoliReinforcementLearningRobust2021}, reachability theory\cite{nathersonReachabililtybasedMissedThrust2026} for missed thrust-aware trajectory design. Passive safety under MTEs or control loss has also been studied for rendezvous\cite{alkhenaniChanceConstrainedAutonomousSpacecraft,elangoSuccessiveConvexificationPassivelySafe2025}, small-body proximity operations\cite{SafeAutonomousGuidance}, and icy-moon transfers\cite{nuccioStatisticalTrajectoryOptimization2026}.

Most of these methods approach MTE robustness from a deterministic, worst-case perspective, which carries several inherent challenges. First, enforcing a robustness margin over the full trajectory can lead to more conservative solutions. Second, the explicit trade-off between robustness and primary mission objectives such as fuel consumption introduces additional tuning parameters into the design process. Third, simultaneously optimizing the nominal and recovery trajectories across multiple failure scenarios can lead to high-dimensional problems that scale poorly as mission complexity grows. These challenges are further compounded in multibody environments, where the chaotic sensitivity of the dynamics makes deterministic margin constraints increasingly difficult to satisfy.

On the other hand, stochastic optimal control-based approaches remain relatively 
unexplored in the context of missed thrust design, despite showing promising results\cite{ozakiStochasticDifferentialDynamic2018,ridderhofChanceconstrainedCovarianceControl2020,kumagaiRobustCislunarLowThrust,grecoRobustSpaceTrajectory2022} 
in countering other operational uncertainties such as navigation and maneuver 
execution errors. Existing stochastic frameworks\cite{ozakiStochasticDifferentialDynamic2018,
ridderhofChanceconstrainedCovarianceControl2020,kumagaiRobustCislunarLowThrust} typically assume Gaussian disturbances, which poorly capture the uncertainty introduced by thruster failures. Previous flight data have shown to agree well with exponential and Weibull distributions\cite{imkenModelingSpacecraftSafe2018}, further motivating distribution-based trajectory design frameworks. Previous works have used such distributions in trajectory design in various ways. Deterministic formulations consider either the mean trajectory only\cite{rubinsztejnDesigningTrajectoriesResilient2021} or use scenario-based sequential convex programming (SCP)\cite{sekineRobustTrajectoryDesign2026}. More recent works have considered such distributions in a probabilistic approach using the Arrow-Hurwicz stochastic algorithm\cite{carpentierLowThrustInterplanetaryTrajectories}, a scenario-based differential dynamic programming framework\cite{shibukawaDDPBasedStochasticTrajectory}, and a missed thrust sensitivity metric-based trajectory search combined with finite horizon linear quadratic regulator for rendezvous applications\cite{nemirovskyTrajectoryDesignControl2026}. While these works study the MTE-resilient trajectory design problem through varied approaches, they either suffer from computational scalability challenges or consider only one MTE, and have not been demonstrated in chaotic nonlinear environments such as cislunar trajectory design. To the best of the authors' knowledge, no prior work has addressed the MTE-resilient trajectory design problem within a scalable stochastic optimal control framework, one whose computational complexity does not grow with the number of MTEs, that accommodates nonlinear spacecraft dynamics, and that controls both the mean and the dispersion of the spacecraft state.

To address this gap, we propose an SCP-based trajectory optimization framework that models the spacecraft as a stochastic hybrid system. Given that the exponential distribution well represents the MTE onset and recovery times\cite{imkenModelingSpacecraftSafe2018}, we model the switching as a Markov chain. Such systems are referred to as Markov Jump Systems (MJS)\cite{costaDiscreteTimeMarkovJump2005}, and are a subclass of stochastic hybrid systems. Despite its versatility, MJS applications in astrodynamics have been mostly limited to rendezvous\cite{shrivastavaFaultTolerantFeedbackControl2023}, stationkeeping\cite{shrivastavaRobustHaloOrbit2026} and fault-tolerant attitude control\cite{fuModelReferenceTracking2017}, with no prior work that utilizes MJS for robust trajectory optimization. 

Therefore, in this work, we utilize the Markov jump system framework to develop robust trajectory optimization frameworks. To this end, we leverage the recently developed chance-constrained covariance steering for Markov-jump linear systems\cite{shrivastavaChanceConstrainedCovarianceSteering2025,wangCovarianceSteeringDiscreteTime2026}, and combine it with a sequential convex programming-based covariance steering approach. To demonstrate our controller's performance, we include numerical simulations with multiple missed thrust events in the CR3BP dynamics. The main contributions of the paper are as follows: 1) We present the first application of covariance steering theory for nonlinear Markov jump systems using Sequential Convex Programming, 2) We formulate the missed thrust problem as a stochastic optimal control problem by modeling it as a Markov jump system, and leverage our framework to solve it under chance constraints, 3) We present a numerically stable and scalable framework to solve the nonlinear Markov jump covariance steering problem using the deviation-based lifted state formulation, and 4) We demonstrate some of the first work that robustifies trajectories in cislunar space against an arbitrary number of MTEs.


\subsection{Notation} $\mathbb{R}$, $\mathbb{N}$ and $\mathbb{Z}_{a:b}$ denote the set of real numbers, natural numbers and set of integers from $a$ to $b$, respectively. For a matrix $X$, $X^\top$, $\|X\|$, tr($X$) refers to its transpose, maximum singular value, and trace, respectively. For a vector $\mathbf{x}\in\mathbb{R}^n$, $\|\mathbf{x}\|$ indicates the $l_2$-norm of the vector. $I_n$ denotes the identity matrix of size $n\times n$. $\lambda_{\max}(X)$ is the largest eigenvalue of $X$. $E_r$ denotes the selection matrix that extracts the position from the state vector, and for a three-dimensional trajectory with $\mathbf{x}\in\mathbb{R}^6$, $E_r = [I_3\;\;\mathbf{0}_{3\times3}]$. For a scalar random variable $Z$ and $p\in(0,1)$, $\mathcal{Q}(Z;p)$ denotes its $p$-quantile.

\section{Problem Formulation}

In this section, we formulate the MTE-resilient trajectory optimization problem as a stochastic optimal control problem, and discuss the dynamics and statistical modeling of the MTEs.

\subsection{Equations of Motion}

We model the spacecraft dynamics using the Circular Restricted Three-Body Problem (CR3BP), which is a commonly used model in preliminary mission design. It considers the motion of a spacecraft, also referred as the secondary body, in the vicinity of the two primary bodies. It assumes that the spacecraft mass is negligible and does not affect the motion of the primaries, which move in a circular orbit about their barycenter. The spacecraft's motion is described relative to a rotating frame moving at the same angular velocity as the primaries. 

It is standard to study the CR3BP dynamics in non-dimensionalized units. The equations of motion for the spacecraft are given below\cite{connorhowellThreedimensionalPeriodicHalo1984}:
\begin{equation}
\label{eom}
    \begin{aligned}
\ddot{x} &= x+ 2\dot{y}  -\frac{(1-\mu)(x+\mu)}{d_1^3} - \frac{\mu(x-1+\mu)}{d_2^3} \\
\ddot{y} &=- 2\dot{x} +y  -\frac{(1-\mu)y}{d_1^3} - \frac{\mu y}{d_2^3} \\
\ddot{z} &= -\frac{(1-\mu)z}{d_1^3} - \frac{\mu z}{d_2^3}
\end{aligned}
\end{equation}
where  $\quad d_1=\sqrt{(x + \mu)^2 + y^2 + z^2},\quad d_2 = \sqrt{(x - 1 + \mu)^2 + y^2 + z^2},\quad \mu = \frac{m_2}{m_1 + m_2}$

For the state space form, we define the state vector $\mathbf{x}$ as $\mathbf{x}=[x,y,z,\dot{x},\dot{y},\dot{z}]^\top\in\mathbb{R}^6$. Similarly the control is denoted by $\mathbf{u}\in\mathbb{R}^3$. The state space model is given below:
\begin{equation}\label{statespace_cr3bp}
    \dot{\mathbf{x}} = \mathbf{f}(\mathbf{x}) + B \mathbf{u}
\end{equation}

where $\mathbf{f}(\mathbf{x})$ is based on the dynamics in \cref{eom} and $B=\begin{bmatrix}
        \mathbf{0}_{3\times3}&
        \mathbf{I}_3
    \end{bmatrix}^\top$. 
\subsection{MTE-Modeling via Markov Jump Systems}

To model MTEs, we utilize a stochastic hybrid system framework called Markov jump systems (MJS). Specifically, we consider a subcategory of MJS called sampled-data Markov jump system (SDMJS), in which the control and mode are discrete-time variables, and can only change at the finite sampling times.

The SDMJS is defined on the probability space $(\Omega,\mathcal{M},\mathbb{P})$, and represented as:

\begin{equation}\label{sdmjs}
    \mathcal{H} : \begin{cases} 
    \dot{\mathbf{x}}(t) = \mathbf{f}(\mathbf{x}(t),\, r_k) + {B}(r_k)\,\mathbf{u}(r_k,k), 
    & t \in [k\Delta_T,\,(k+1)\Delta_T) 
    \\[6pt] 
   \mathbb{P}(r_{k+1} = j \mid r_k = i) = p_{ij}, \quad i,j \in \mathcal{S}, 
    & t = (k+1)\Delta_T 
    \end{cases}
\end{equation}

If we define $\mathcal{S}$ as a singleton set, i.e. $r_k =1,\,\forall k$, then \cref{sdmjs} is equivalent to the standard nonlinear dynamics equation for single-mode systems, with zero-order hold (ZOH) control. Incorporating the mode dynamics generalizes the dynamics to allow for sudden switches in the dynamics function $\mathbf{f}$, or the control matrix ${B}$. This additional parameter $r_k$ is referred to as the active mode at sampling time $k\Delta_T$, and it can take values from the discrete set $\mathcal{S} = \{1,2,\cdots N\}$ for an $N$-mode Markov jump system. The switching between these modes is governed by a homogeneous Markov chain with transition probability $p_{ij}$ where $p_{ij} = \mathbb{P}(r_{k+1}=j|r_k=i)$. We also define the transition probability matrix $\mathcal{P}$ whose ($i$,$j$)$^\mathrm{th}$ entry is $p_{ij}$. The mode probability distribution vector is given by $\bm{\rho}_k=[\rho_k(1)\ldots\rho_k(N)]$. For a given initial mode distribution $\bm{\rho}_0$, $\bm{\rho}_k$ can be calculated as $\bm{\rho}_k=\bm{\rho}_0\mathcal{P}^k$.

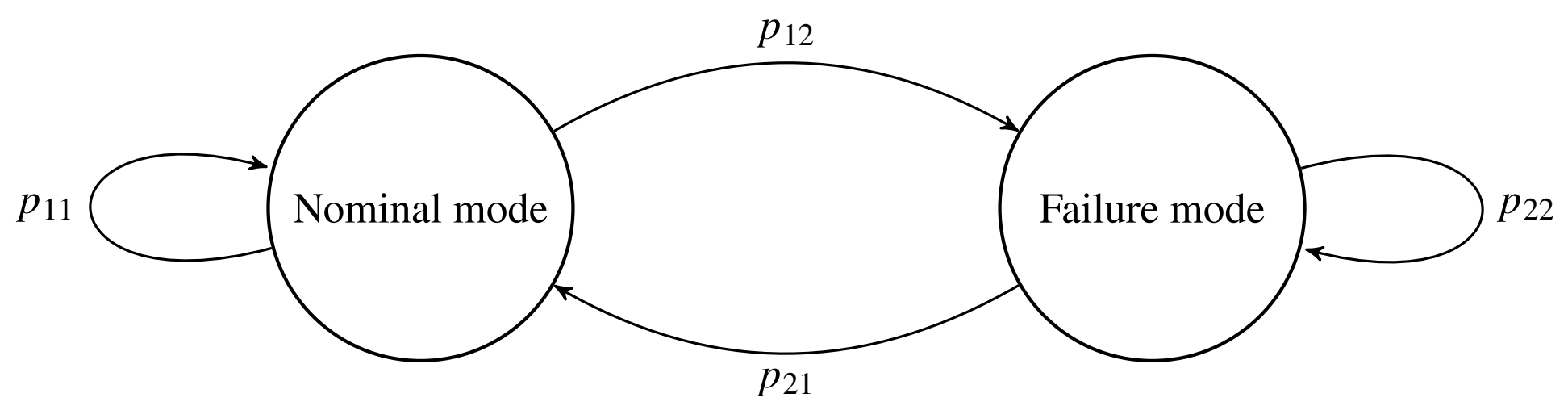
\begin{figure}[tbp!]
    \centering
    \begin{tikzpicture}[->, >=stealth', auto, semithick, node distance=6cm]
    \tikzstyle{every state}=[fill=white,draw=black,thick,text=black,scale=1,minimum size=2.5cm]

    \node[state] (A) {Nominal mode};
    \node[state] (B) [right of=A] {Failure mode};

    \path
    (A) edge[loop left] node{$p_{11}$} (A)
        edge[bend left] node{$p_{12}$} (B)
    (B) edge[loop right] node{$p_{22}$} (B)
        edge[bend left] node{$p_{21}$} (A);
    \end{tikzpicture}
    \caption{A schematic showing the Markov chain transition between nominal and failure modes}
    \label{fig:framework1}
\end{figure}

In this study, we consider a dual mode MJS description, with nominal and failure modes, but the framework can be extended to consider more modes, such as partial control failure (one of the thrusters stop working), or to model sudden changes in the dynamics. For the dual mode framework, we model MTEs as a failure mode, and the nominal mode as the mode when all thrusters are working. \Cref{fig:framework1} shows a Markov chain with two modes: mode 1 is labeled nominal mode and mode 2 is labeled failure mode. The associated transition probability with the arrows is also included. For instance, $p_{12}$ denotes the probability of transition from mode 1, i.e. nominal mode to mode 2, i.e. failure mode. Under the given modeling framework, the dynamics $\mathbf{f}(\mathbf{x},r_k)$ is the same for both modes, and the control matrix ${B}(r_k)$ is zero for the failure mode, which models the MTE, i.e. 

\begin{equation}\label{f_mjls}
    \mathbf{f}(\mathbf{x},1)= \mathbf{f}(\mathbf{x},2) = \mathbf{f}(\mathbf{x}),\quad B_k(r_k) = \begin{cases}
        \begin{bmatrix}
            \mathbf{0}_{3\times3}\\
            \mathbf{I}_3
        \end{bmatrix} & r_k=1\\
        \mathbf{0}_{6\times3} & r_k=2
    \end{cases}
\end{equation}

The transition probabilities $p_{ij}$ are derived from the MTE statistics. We characterize the MTEs by two parameters, $\text{TTF}_{90}$ and $\text{TTR}_{95}$, defined as the 90th percentile failure time and the 95th percentile recovery time, respectively. Therefore, once the spacecraft is in nominal mode, it will face a MTE within $\text{TTF}_{90}$ with at least $90\%$ probability, and once an MTE begins, recovery occurs within $\text{TTR}_{95}$ with $95\%$ probability. Since we model the MTE using an exponential distribution, the corresponding rate parameters are given by:
\begin{equation}\label{exp_rates}
    1-e^{-\lambda_F\,\text{TTF}_{90}} = 0.9,\qquad 1-e^{-\lambda_R\,\text{TTR}_{95}} = 0.95
    \quad\Longrightarrow\quad
    \lambda_F = \frac{\ln 10}{\text{TTF}_{90}},\qquad \lambda_R = \frac{\ln 20}{\text{TTR}_{95}}
\end{equation}
For the SDMJS defined in \cref{sdmjs} with time grid $\Delta_T$, the corresponding transition probability matrix is given by:
\begin{equation}\label{tpm_exp}
    \mathcal{P} = \begin{bmatrix}
        p_{11} & p_{12}\\
        p_{21} & p_{22}
    \end{bmatrix} = \begin{bmatrix}
        e^{-\lambda_F\Delta_T} & 1-e^{-\lambda_F\Delta_T}\\
        1-e^{-\lambda_R\Delta_T} & e^{-\lambda_R\Delta_T}
    \end{bmatrix}
\end{equation}
 
\begin{remark}
    Spacecrafts typically have a fixed sampling time for the control, and the thruster magnitude is rarely changed continuously. Moreover, most direct methods assume zero order hold (ZOH) for control, and the control is constant within each sampling interval. Therefore, the SDMJS framework is a justified assumption for spacecraft trajectory optimization.
\end{remark}

\subsection{Stochastic Optimal Control Formulation}

We can now define the stochastic optimal control problem we aim to solve in this work:

\begin{problem}\label{prob1} Solve:
    \begin{subequations}
    \begin{flalign}
        &\min_{\mathbf{x},\mathbf{u}} J = \mathcal{Q}\left(\int_{t_0}^{t_f} \|\mathbf{u}(t)\| dt; 0.99\right) \\
        &\text{s.t.}\;\;\;\forall k \in \mathbb{Z}_{0:T-1}, \notag\\
        &\quad \mathcal{H} : \begin{cases}
        \dot{\mathbf{x}}(t) = \mathbf{f}(\mathbf{x}(t),\,r_k) + {B}(r_k)\,\mathbf{u}(r_k,k), & t \in [k\Delta_T,\,(k+1)\Delta_T) \\[6pt]
        \mathbb{P}(r_{k+1} = j \mid r_k = i) = p_{ij},\quad i,j \in \mathcal{S}, & t = (k+1)\Delta_T
        \end{cases} && \clabel{Nonlinear SDMJS dynamics} \label{ocp_nldynamics}\\
        &\quad \mathbf{x}(0) \sim \mathcal{N}(\overline{\bm{\mu}}_0,\overline{\Sigma}_0), \; \mathbb{E}[\mathbf{x}(t_f)]=\bm{\mu}_f, \; cov(\mathbf{x}(t_f))\preceq\Sigma_f && \clabel{Boundary conditions} \label{bc}\\
        &\quad \mathbb{P}(\|\mathbf{u}(1,k)\| \leq u_{\max})\geq 1 - \epsilon_u && \clabel{Control chance constraint} \label{cc}\\
        &\quad \mathbb{P}\bigl(\|E_r\mathbf{x}_k-\mathbf{r}_{m}\| \geq R_{\mathrm{safe}}\bigr) \geq 1 - \epsilon_{\mathrm{ko}} && \clabel{Keep-out zone chance constraint} \label{cc_ko}\\
        &\quad cov(\mathbf{x}_k) \preceq \Sigma_{\max} && \clabel{Mid-course covariance bound} \label{cov_cap}
    \end{flalign}
    \end{subequations}
\end{problem}

\begin{assumption}
    The distribution for the initial mode distribution $\bm{\rho}_0$ and initial state $\mathbf{x}_0$ are independent.
\end{assumption}

\begin{assumption}\label{assum:rho_positive}
    The mode probability is nonzero for every mode at every time step, i.e.\ $\rho_k(i)>0$ $\forall i\in\mathcal{S},\,k\in\mathbb{Z}_{0:T}$.
\end{assumption}

\Cref{prob1} considers a $\Delta v_{99}$ cost minimization, subject to nonlinear SDMJS dynamics. We also assume the initial state is Gaussian distributed, and the final state is constrained to be within a target covariance of the target state. We also consider path constraints in the form of chance constraints: \Cref{cc} is the control bound chance constraint and \Cref{cc_ko} is a {keep-out zone} chance constraint. $E_r\triangleq[I_3\;\;\mathbf{0}_{3\times3}]$ denotes the selector matrix extracting the position from the state, $\mathbf{r}_{m}$ is the center of one of the primaries, such as the Moon, and $R_{\mathrm{safe}}$ is the prescribed safety radius. For cislunar transfers with complex lunar flybys, such keep-out zone constraints are critical to ensure low probability of impact with the lunar surface under various uncertainties. Finally, \cref{cov_cap} is a mid-course covariance bound that bounds the state dispersion at every node by a prescribed $\Sigma_{\max}$ and is useful to limit covariance growth throughout the trajectory.

{\Cref{prob1} is a computationally challenging problem for multiple reasons. Firstly, standard direct methods cannot deal with stochastic decision variables. Scenario-based methods are one approach to switching dynamics, but are not scalable to higher dimensions, or longer time horizons. Secondly, the nonlinear MJS dynamics are not only difficult to linearize about, but also introduce additional non-Gaussianity in the dynamics, rendering common stochastic control approaches difficult to apply. Even in the linearized case, the Markov jump dynamics lead to an exponential number of Gaussians needed to represent the state distribution, and the nonlinear dynamics already introduce non-Gaussianity. Therefore, we need to develop a more efficient framework to solve this problem.}

\section{Solution using Covariance Steering for Markov Jump Systems}

In this section, we describe our solution framework for \cref{prob1}. We leverage covariance steering theory for Markov jump Linear Systems (MJLS) and combine it with a SCP-based framework to deal with nonconvex constraints. 

\subsection{Linearization and Discretization of Nonlinear MJS}\label{sec:linearization}

To discretize and linearize the nonlinear MJS dynamics, we first note an analogy between the mean MJS dynamics and single-mode nonlinear dynamics. Taking the mean of the state propagation equation in \cref{sdmjs} gives:
\begin{equation}\label{sdmjs_mean}
    \mathbb{E}[\dot{\mathbf{x}}(t)] = \mathbb{E}[\mathbf{f}(\mathbf{x}(t),\, r_k)] + \mathbb{E}[{B}(r_k)\mathbf{u}(r_k,k)]
\end{equation}
Given the prior definition of $\mathbf{f}(\mathbf{x}(t),\, r_k)$ and $B(r_k)$ in \cref{f_mjls}, we first take the expectation over the current mode $r_k$ to get $\mathbb{E}[\mathbf{f}(\mathbf{x}(t),\, r_k)]=\mathbb{E}[\mathbf{f}(\mathbf{x}(t))]$ and $\mathbb{E}[{B}(r_k)\mathbf{u}(r_k,k)]=\rho_k(r_k=1){B}(r_k=1)\bar{\mathbf{u}}(r_k=1,k)$. A first-order approximation of the mean of the nonlinear MJS dynamics can be expressed as:
\begin{equation}\label{sdmjs_mean_2}
    \mathbb{E}[\dot{\mathbf{x}}(t)] =\mathbb{E}[\mathbf{f}(\mathbf{x}(t))]+ \rho_k(r_k=1){B}\bar{\mathbf{u}}(r_k=1,k)\approx \mathbf{f}(\mathbb{E}[\mathbf{x}(t)]) + \rho_k(r_k=1){B}\bar{\mathbf{u}}(r_k=1,k)
\end{equation}
Note that this is similar to a single-mode nonlinear dynamical system $\dot{\mathbf{x}}(t)=\mathbf{f}(\mathbf{x}(t)) +B\mathbf{u}$, with $\mathbf{u}$ scaled by the nominal mode probability $\rho_k(r_k=1)$. Therefore, the mean of the nonlinear MJS dynamics obeys similar dynamics as a single-mode nonlinear dynamical system with control scaled by $\rho_k(r_k)$. We leverage this to propose the following discretization and linearization procedure for the nonlinear SDMJS dynamics considered in \cref{prob1}:

\begin{enumerate}
    \item Given a reference MJLS trajectory $\{\bm{\mu}^{\mathrm{ref}}_k, \overline{\mathbf{u}}^{\,\mathrm{ref}}_k(i)\}_{k=0}^{T-1}$, $i\in\mathcal{S}$, define the corresponding single-mode reference trajectory as $\{\bm{\mu}^{\mathrm{ref}}_k, \rho_k(r_k=1)\overline{\mathbf{u}}^{\,\mathrm{ref}}_k(r_k=1)\}_{k=0}^{T-1}$.     
    \item Linearize and discretize the single-mode nonlinear dynamics as explained below:
    For a nonlinear dynamical system $\dot{\mathbf{x}}(t)=\mathbf{f}(\mathbf{x}(t)) +B\mathbf{u}$, the discrete-time linear system about the single-mode reference trajectory $\{\bm{\mu}^{\mathrm{ref}}_k, \rho_k(r_k=1)\overline{\mathbf{u}}^{\,\mathrm{ref}}_k(r_k=1)\}_{k=0}^{T-1}$ from Step 1 is given by:
    
    \begin{equation}
        \mathbf{x}_{k+1} = A_k \mathbf{x}_k + B_k \mathbf{u}_k + \mathbf{c}_k
    \end{equation}
    where $A_k$ and $B_k$ are given by\cite{oguriChanceConstrainedControlSafe2024}:
    \begin{equation}\label{lindisc}
    \begin{aligned}
    & A_k = \Phi\left(t_{k+1}, t_k\right), \quad B_k = 
    A_k \int_{t_k}^{t_{k+1}} \Phi^{-1}\left(t, t_k\right) B(t) \mathrm{d} t, \quad \mathbf{c}_k = 
    A_k \int_{t_k}^{t_{k+1}} \Phi^{-1}\left(t, t_k\right) \mathbf{c}(t) \mathrm{d} t
    \end{aligned}
    \end{equation}
    where $\Phi$ is the State Transition Matrix (STM), and it satisfies $ \frac{\mathrm{d}}{\mathrm{~d} t} \Phi\left(t, t_k\right)=A(t) \Phi\left(t, t_k\right), \quad \Phi\left(t_k, t_k\right)=I_{n_x}$.
    \item After obtaining the linearized matrices $A_k, B_k$, define the corresponding MJLS matrices as follows:
     $$A_k(r_k=1)=A_k(r_k=2)=A_k, B_k(r_k=1)=B_k, B_k(r_k=2)=\mathbf{0}_{n_x\times n_u}$$
\end{enumerate}
\begin{remark}
    We do not claim that this is an exact linearization of the original SDMJS. Under our MTE-based modeling, the
    mean of the SDMJS evolves similar to a single-mode nonlinear dynamical system with control scaled by the nominal mode probability $\rho_k(r_k=1)$. This allows us to apply the single-mode discretization and linearization process about the mean trajectory and define a corresponding MJLS. To our knowledge, an exact procedure for discretization and linearization of nonlinear Markov jump systems is unknown at the moment and is left for future work.
    \end{remark}
\subsection{Markov Jump Linear Systems}

The linearization performed in the previous section motivates us to discuss the Markov Jump Linear Systems (MJLS), which is a subclass of MJS wherein the subsystem dynamics for each mode are linear. A discrete-time MJLS, defined on the probability space $(\Omega,\mathcal{M},\mathbb{P})$ is given by:
\begin{equation}\label{MJLS}
    \mathbf{x}_{k+1} = A_k(r_k)\mathbf{x}_k + B_k(r_k)\mathbf{u}_k + \mathbf{c}_k(r_k) 
\end{equation}
The mode dynamics are defined by a time-homogeneous Markov chain, similar to nonlinear MJS, determined by the transition probability matrix $\mathcal{P}$. We now define the following quantities, which are standard in MJLS literature\cite{costaDiscreteTimeMarkovJump2005,shrivastavaChanceConstrainedCovarianceSteering2025}
\begin{equation}\label{mjls_moment_defs}
\begin{aligned}
    &\rho_k(i) \triangleq \mathbb{E} \left[\mathds{1}_{\{r_k=i\}}\right],
        & \mathbf{q}_k(i) \triangleq \mathbb{E} \left[\mathbf{x}_k \mathds{1}_{r_k=i} \right],   \\
        &\bm{\mu}_k\triangleq \mathbb{E}[\mathbf{x}_k]=\sum_{i\in\mathcal{S}} \mathbf{q}_k(i),
        & \Sigma_k\triangleq \mathbb{E} \left[(\mathbf{x}_k-\bm{\mu}_k) (\mathbf{x}_k-\bm{\mu}_k)^\top \right]\\
        &\overline{\mathbf{x}}_k(i)\triangleq\frac{\mathbf{q}_k(i)}{\rho_k(i)},& \hspace{-12mm}S_k(i) \triangleq \mathbb{E}[(\mathbf{x}-\overline{\mathbf{x}}_k(i))(\mathbf{x}-\overline{\mathbf{x}}_k(i))^\top\mathds{1}_{r_k=i}]
\end{aligned}
\end{equation}

To study the evolution of MJLS dynamics, instead of dealing with the state mean and covariance, the MJLS dynamics motivate us to work with $\mathbf{q}_k(i),\overline{\mathbf{x}}_k(i)$ and $S_k(i)$, which are related to the system's mean $\bm{\mu}_k$ and covariance $\Sigma_k$ as:
\begin{subequations}
    \begin{align}
    \bm{\mu}_k &= \sum_{i\in\mathcal{S}}\mathbf{q}_k(i) , \;\; \overline{\mathbf{x}}_k(i) =\frac{\mathbf{q}_k(i)}{\rho_k(i)}\label{mean_def}\\
    \Sigma_k &= \sum_{i\in\mathcal{S}} \left(S_k(i)  + \rho_k(i)(\overline{\mathbf{x}}_k(i)-\bm{\mu}_k)(\overline{\mathbf{x}}_k(i)-\bm{\mu}_k)^\top\right) \label{cov_def}
\end{align}
\end{subequations} 
We consider a mode-dependent affine feedback control policy, which is standard in covariance steering literature \cite{liuOptimalCovarianceSteering2024,rapakouliasDiscreteTimeOptimalCovariance}. The feedback component controls the dispersion in the single-mode system, and we extend that control law to mode dependent control for MJLS as follows:
\begin{equation}\label{controllaw}
    {\mathbf{u}_k(r_k=i) = \overline{\mathbf{u}}_k(r_k=i) + K_k(i) \left(\mathbf{x}_k-\overline{\mathbf{x}}_k(i)\right)}
\end{equation}
The MJLS mean and covariance propagation is given by \cite{shrivastavaChanceConstrainedCovarianceSteering2025}:
\begin{align}
    &\mathbf{q}_{k+1}(j)={\sum_{i\in\mathcal{S}}p_{ij}\left({A}_k(i)\mathbf{q}_k(i) +\rho_k(i)(B_k(i)\overline{\mathbf{u}}_k(i) + \mathbf{c}_k)\right)}
    \label{meanprop1}\\
    & S_{k+1}(j)=\sum_{i\in\mathcal{S}}p_{ij}\Bigl[\bigl(A_k(i)+B_k(i)K_k(i)\bigr)S_k(i)\bigl(A_k(i)+B_k(i)K_k(i)\bigr)^\top+\rho_k(i)\bigl(B_k(i)\overline{\mathbf{u}}_k(i)\overline{\mathbf{u}}_k^\top(i)B_k^\top(i)\notag\\
    &\quad+B_k(i)\overline{\mathbf{u}}_k(i)\mathbf{c}_k^\top+B_k(i)\overline{\mathbf{u}}_k(i)\overline{\mathbf{x}}_k^\top(i)A_k^\top(i)+\mathbf{c}_k\mathbf{c}_k^\top+\mathbf{c}_k\overline{\mathbf{u}}_k^\top(i)B_k^\top(i)+\mathbf{c}_k\overline{\mathbf{x}}_k^\top(i)A_k^\top(i)+A_k(i)\overline{\mathbf{x}}_k(i)\overline{\mathbf{u}}_k^\top(i)B_k^\top(i)\notag\\
    &\quad+A_k(i)\overline{\mathbf{x}}_k(i)\mathbf{c}_k^\top+A_k(i)\overline{\mathbf{x}}_k(i)\overline{\mathbf{x}}_k^\top(i)A_k^\top(i)\bigr)\Bigr]-\rho_{k+1}(j)\overline{\mathbf{x}}_{k+1}(j)\overline{\mathbf{x}}_{k+1}^\top(j)\label{covprop1} 
\end{align}

We note that \cref{covprop1} is a nonlinear equation in both mean and covariance quantities, and cannot be directly used in a convex optimization framework. This is a pertinent issue in covariance steering for Markov Jump Systems\cite{shrivastavaChanceConstrainedCovarianceSteering2025,wangCovarianceSteeringDiscreteTime2026}, and we discuss our approach to tackle this in the following sections.

\subsection{Linearized Covariance Steering Subproblem for MJLS}

To solve \cref{prob1}, we first consider a simplified problem, wherein we linearize the nonlinear dynamics of \cref{ocp_nldynamics}. This is an approximation to the original problem, but when solved iteratively, can lead to a locally optimal solution to \cref{prob1}. This is the working principle for most sequential convex programming algorithms.
\begin{problem}\label{prob2} Solve:
    \begin{subequations}
    \begin{flalign}
        &\min_{\mathbf{x}_k,\mathbf{u}_k(i)} J = \mathcal{Q}\left(\sum_{k=0}^{T-1} \|\mathbf{u}_k\|\Delta T; 0.99\right) + P(\bm{\xi}) \\
        &\text{s.t.}\;\;\;\forall k \in \mathbb{Z}_{0:T-1}, \notag\\
        &\quad \mathbf{x}_{k+1} - A_k (r_k)\mathbf{x}_k - B_k(r_k) \mathbf{u}_k - \mathbf{c}_k(r_k)=\bm{\xi}_k && \clabel{Linearized MJLS dynamics} \label{mean_dyn_linearized_prob2} \\ 
        &\quad \mathbf{x}_0\sim\mathcal{N}(\overline{\bm{\mu}}_0,\overline{{\Sigma}}_0),  \; \mathbb{E}[\mathbf{x}_T]=\bm{\mu}_f, \; {\Sigma}_T\preceq{\Sigma}_f && \clabel{Boundary conditions} \label{bc2}\\
        &\quad \mathbb{P}(\|\mathbf{u}_k(1)\| \leq u_{\max}) \geq 1 - \epsilon_u && \clabel{Control chance constraint} \label{cc2}\\
        &\quad \mathbb{P}\bigl(\|E_r\mathbf{x}_k-\mathbf{r}_{m}\| \geq R_{\mathrm{safe}}\bigr) \geq 1 - \epsilon_{\mathrm{ko}} && \clabel{Keep-out zone chance constraint} \label{cc_ko2}\\
        &\quad \Sigma_k \preceq \Sigma_{\max} && \clabel{Mid-course covariance bound} \label{cov_cap2}
    \end{flalign}
    \end{subequations}
where $i\in\mathcal{S}$, $r_k$ forms a Markov chain, $\Delta T$ is the sampling interval, and $P(\bm{\xi})$ is a penalty function for slack variable  vector $\bm{\xi}$,  containing $\{\bm{\xi}_k,\; \forall k \in \mathbb{Z}_{0:T-1}\}$. The distribution $\bm{\rho}_0$ and $\mathbf{x}_0$ are independent.
\end{problem}

We add the penalty function for the slack $\bm{\xi}$ to relax the linearized dynamics, which is standard in SCP algorithms to tackle artificial infeasibility. 

\Cref{prob2} remains a computationally intractable problem, as \cref{mean_dyn_linearized_prob2,cc2,cc_ko2} are stochastic constraints that cannot directly be enforced in an optimization framework. For the stochastic MJLS dynamics, we instead steer the mean and covariance of the MJLS, and the following subsection develops a convex reformulation of the resulting covariance steering problem.

\subsection{Convex Reformulation of the MJLS Covariance Steering Problem}\label{dev_framework}

We introduce two major modifications over our previous work\cite{shrivastavaMASTOMISSEDTHRUSTAWARE2026} to solve \Cref{prob2}: First, we leverage a lifted-state formulation, inspired by a recent work\cite{wangCovarianceSteeringDiscreteTime2026}, which solves the unconstrained MJLS covariance steering problem in a unified framework by propagating the second moment of a lifted state vector. This allows us to solve the mean and covariance subproblems jointly, simplifying the two-shot alternating framework proposed previously\cite{shrivastavaMASTOMISSEDTHRUSTAWARE2026}. However, a nominal implementation of such lifted state-based frameworks leads to very ill-conditioned optimization problems in astrodynamics. The numerical issues in covariance steering formulations for astrodynamics applications are well-documented, requiring additional scaling efforts to improve numerical convergence\cite{kumagaiRobustCislunarLowThrust}. The lifted state formulation exacerbates this problem significantly since it deals with the second moment propagation, and enforcing covariance constraints becomes numerically ill-conditioned. This is caused by the much smaller covariance values compared to the mean quantities, which is typical for most spacecraft trajectory optimization problems. To solve this, we propose a deviation-based linearization approach, wherein we define the state vector as a deviation from the reference trajectory. The resulting deviation-based state has mean and covariance in similar order of magnitude, vastly improving numerical properties of the optimization problem. While we introduce it for MJLS covariance steering, such formulations should also benefit other covariance steering problems where the mean and covariance quantities are of very different scales. 

\subsubsection{Deviation-Based State Formulation}

Let $\{\bm{\mu}^{\mathrm{ref}}_k,\overline{\mathbf{u}}^{\,\mathrm{ref}}_k(i)\}_{k=0}^{T-1}$ be the reference MJLS trajectory about which the dynamics are linearized and discretized, as explained in \cref{sec:linearization}. We define the deviation variables with respect to this reference as: 
\begin{equation}\label{dev_vars} \
    \delta\mathbf{x}_k \triangleq \mathbf{x}_k - \bm{\mu}^{\mathrm{ref}}_k, \qquad
    \delta\mathbf{u}_k(i) \triangleq \mathbf{u}_k(i) - \overline{\mathbf{u}}^{\,\mathrm{ref}}_k(i)
\end{equation}
Here, $\bm{\mu}^{\mathrm{ref}}_k$ is the reference state mean and $\overline{\mathbf{u}}^{\,\mathrm{ref}}_k(i)$ is the reference nominal mode-dependent control. Under the current bimodal MJLS structure with a nominal and failure mode, the reference mean trajectory is propagated as: 
\begin{equation}\label{ref_dyn}
    \bm{\mu}^{\mathrm{ref}}_{k+1} = A_k\,\bm{\mu}^{\mathrm{ref}}_k + \rho_k(r_k=1)B_k\,\overline{\mathbf{u}}^{\,\mathrm{ref}}_k(r_k=1) + \mathbf{c}_k
\end{equation}

Subtracting \cref{ref_dyn} from the MJLS dynamics given in \cref{MJLS}, we get the  corresponding linear discrete-time system dynamics for the deviation-based variables as:
\begin{equation}\label{dev_dyn}
    \delta\mathbf{x}_{k+1} = A_k(r_k)\,\delta\mathbf{x}_k + B_k(r_k)\,\delta\mathbf{u}_k + \mathbf{d}_k(r_k) 
\end{equation}
where the term $\mathbf{d}_k(i)$ is defined as:
\begin{equation}\label{dev_drift}
    \mathbf{d}_k(i) \triangleq \bigl(A_k(i)-A_k\bigr)\bm{\mu}^{\mathrm{ref}}_k + B_k(i)\overline{\mathbf{u}}^{\,\mathrm{ref}}_k(i) - B_k\mathbf{u}^{\mathrm{ref}}_k + \mathbf{c}_k(i) - \mathbf{c}_k
    = A_k(i)\,\bm{\mu}^{\mathrm{ref}}_k + B_k(i)\,\overline{\mathbf{u}}^{\,\mathrm{ref}}_k(i) + \mathbf{c}_k(i) - \bm{\mu}^{\mathrm{ref}}_{k+1}
\end{equation}
In \cref{dev_drift}, the quantities without a mode argument $i$ denote the single-mode quantities used for linearization in \cref{sec:linearization}. For our bimodal MJLS framework (\cref{f_mjls}), we have $A_k(1)=A_k(2)=A_k$, $B_k(1)=B_k$, $B_k(2)=\mathbf{0}$, $\mathbf{c}_k(1)=\mathbf{c}_k(2)=\mathbf{c}_k$, and $\mathbf{u}^{\mathrm{ref}}_k=\rho_k(1)\,\overline{\mathbf{u}}^{\,\mathrm{ref}}_k(1)$, which gives
\begin{equation}\label{dev_drift_bimodal}
    \mathbf{d}_k(1) = \rho_k(2)\,B_k\,\overline{\mathbf{u}}^{\,\mathrm{ref}}_k(1), \qquad \mathbf{d}_k(2) = -\rho_k(1)\,B_k\,\overline{\mathbf{u}}^{\,\mathrm{ref}}_k(1) .
\end{equation}
The $\mathbf{d}_k(i)$ terms absorb the mismatch between the deviation-based MJLS and the original MJLS equation in \cref{MJLS} and ensure the dynamics remain consistent across the two formulations.  

The MJLS moment quantities are translated to deviation coordinates as
\begin{equation}\label{dev_moments}
    \delta\mathbf{q}_k(i) \triangleq \mathbb{E}\left[\delta\mathbf{x}_k\mathds{1}_{\{r_k=i\}}\right] = \mathbf{q}_k(i)-\rho_k(i)\,\bm{\mu}^{\mathrm{ref}}_k, \qquad
    \delta\overline{\mathbf{x}}_k(i) \triangleq \frac{\delta\mathbf{q}_k(i)}{\rho_k(i)} = \overline{\mathbf{x}}_k(i)-\bm{\mu}^{\mathrm{ref}}_k, \qquad
    \delta\bm{\mu}_k \triangleq \sum_{i\in\mathcal{S}}\delta\mathbf{q}_k(i) = \bm{\mu}_k - \bm{\mu}^{\mathrm{ref}}_k
\end{equation}
Since the partial covariance $S_k(i)$ depends on deviation from the mean, it is invariant under the coordinate shift:
\begin{equation}\label{dev_S_invariance}
    S_k(i) = \mathbb{E}\left[\bigl(\delta\mathbf{x}_k-\delta\overline{\mathbf{x}}_k(i)\bigr)\bigl(\delta\mathbf{x}_k-\delta\overline{\mathbf{x}}_k(i)\bigr)^\top\mathds{1}_{\{r_k=i\}}\right]
\end{equation}
which is identical to the $S_k(i)$ defined in \cref{mjls_moment_defs}. The control policy $\delta\mathbf{u}_k(i)$ in the deviation-based framework is given by: 
\begin{equation}\label{dev_controllaw}
    \delta\mathbf{u}_k(r_k=i) = \delta\overline{\mathbf{u}}_k(i) + K_k(i)\bigl(\delta\mathbf{x}_k - \delta\overline{\mathbf{x}}_k(i)\bigr), \qquad
    \overline{\mathbf{u}}_k(i) = \overline{\mathbf{u}}^{\,\mathrm{ref}}_k(i) + \delta\overline{\mathbf{u}}_k(i)
\end{equation}
Since the feedback acts on the deviation from the mean trajectory, it is invariant under this coordinate change, similar to $S_k(i)$. Given \cref{dev_controllaw}, we can recover the true control $\mathbf{u}_k(i)$ as:
\begin{equation}\label{dev_policy_absolute}
    \mathbf{u}_k(i) = \underbrace{\overline{\mathbf{u}}^{\,\mathrm{ref}}_k(i) + \delta\overline{\mathbf{u}}_k(i)}_{=\,\overline{\mathbf{u}}_k(i)} + \; K_k(i)\bigl(\mathbf{x}_k - \overline{\mathbf{x}}_k(i)\bigr)
\end{equation}
\subsubsection{Lifted State Formulation}
We now define the deviation-based lifted state $\check{\mathbf{x}}_k\triangleq[\delta\mathbf{x}_k^\top\;\;1]^\top\in\mathbb{R}^{n_x+1}$, so that \cref{dev_dyn} becomes:
\begin{equation}\label{devlifted_dyn}
    \check{\mathbf{x}}_{k+1} = \check{A}_k(r_k)\,\check{\mathbf{x}}_k + \check{B}_k(r_k)\,\delta\mathbf{u}_k 
\end{equation}
with $\check{A}_k(r_k)$ and $\check{B}_k(r_k)$ given by:  
\begin{equation}\label{devlifted_matrices}
    \check{A}_k(i) = \begin{bmatrix} A_k(i) & \mathbf{d}_k(i) \\ \mathbf{0}_{1\times n_x} & 1 \end{bmatrix}, \quad
    \check{B}_k(i) = \begin{bmatrix} B_k(i) \\ \mathbf{0}_{1\times n_u} \end{bmatrix} 
\end{equation}
We also define the selection operators $E_x \triangleq [\,I_{n_x}\;\;\mathbf{0}_{n_x\times 1}\,]\in\mathbb{R}^{n_x\times(n_x+1)}$, which extracts the state part of the lifted state as $\delta\mathbf{x}_k=E_x\check{\mathbf{x}}_k$, and $\mathbf{e} \triangleq [\,\mathbf{0}_{1\times n_x}\;\;1\,]^\top\in\mathbb{R}^{n_x+1}$, the selection vector for the unity term, $\mathbf{e}^\top\check{\mathbf{x}}_k=1$. The lifted moments are defined on the deviation variables as follows:
\begin{subequations}\label{devlifted_moments}
\begin{align}
    \check{S}_k(i) &\triangleq \mathbb{E}\left[\check{\mathbf{x}}_k\check{\mathbf{x}}_k^\top\mathds{1}_{\{r_k=i\}}\right] = \begin{bmatrix} S_k(i)+\rho_k(i)\,\delta\overline{\mathbf{x}}_k(i)\,\delta\overline{\mathbf{x}}_k^\top(i) & \delta\mathbf{q}_k(i) \\ \delta\mathbf{q}_k^\top(i) & \rho_k(i) \end{bmatrix} \label{devlifted_S_def}\\
    \check{U}_k(i) &\triangleq \mathbb{E}\left[\delta\mathbf{u}_k\,\check{\mathbf{x}}_k^\top\mathds{1}_{\{r_k=i\}}\right] \in \mathbb{R}^{n_u\times(n_x+1)}, \qquad
    \check{Y}_k(i) \triangleq \mathbb{E}\left[\delta\mathbf{u}_k\,\delta\mathbf{u}_k^\top\mathds{1}_{\{r_k=i\}}\right] \in \mathbb{R}^{n_u\times n_u} \label{devlifted_UY_def}
\end{align}
\end{subequations}
so that the single matrix variable $\check{S}_k(i)$ simultaneously encodes the mode probability $\rho_k(i)=\mathbf{e}^\top\check{S}_k(i)\,\mathbf{e}$, the first moment information $\delta\mathbf{q}_k(i)=E_x\check{S}_k(i)\,\mathbf{e}$, and the second moment information $E_x\check{S}_k(i)E_x^\top=S_k(i)+\rho_k(i)\,\delta\overline{\mathbf{x}}_k(i)\,\delta\overline{\mathbf{x}}_k^\top(i)$. Without the deviation-based reformulation, this term causes the optimization to be ill-conditioned since the covariance dynamics are incorporated through this second moment expression only, and the significantly larger scale of the mean terms dominate the second moment $E_x\check{S}_k(i)E_x^\top$.

The policy \cref{dev_controllaw} can be represented in the lifted state formulation as:
\begin{equation}\label{control_lifted_def}
    \delta\mathbf{u}_k=\check{K}_k(r_k)\check{\mathbf{x}}_k
\end{equation}
with $\check{K}_k(i)\triangleq\bigl[K_k(i),\;\delta\overline{\mathbf{u}}_k(i)-K_k(i)\,\delta\overline{\mathbf{x}}_k(i)\bigr]$. Using \Cref{control_lifted_def} in \Cref{devlifted_UY_def}, we get:
\begin{equation}\label{devlifted_policy_eq}
    \check{U}_k(i)=\check{K}_k(i)\check{S}_k(i),\quad\check{Y}_k(i)=\check{K}_k(i)\check{S}_k(i)\check{K}_k^\top(i)=\check{U}_k(i)\check{S}_k^{-1}(i)\check{U}_k^\top(i)
\end{equation}

For the control, we can derive the mode-dependent mean as $\delta\overline{\mathbf{u}}_k(i)=\check{U}_k(i)\,\mathbf{e}/\rho_k(i)$ and covariance as $cov(\mathbf{u}_k\,|\,r_k=i)=\check{Y}_k(i)/\rho_k(i)-\delta\overline{\mathbf{u}}_k(i)\,\delta\overline{\mathbf{u}}_k^\top(i)$, the latter being invariant under this coordinate change.

Under the lifted state formulation, the second moment matrix $\check{S}_k(i)$ encodes both mean and covariance information, and can be propagated as: 
\begin{align}
    \check{S}_{k+1}(j) &= \sum_{i\in\mathcal{S}}\mathbb{E}\bigl[\check{\mathbf{x}}_{k+1}\check{\mathbf{x}}_{k+1}^\top\mathds{1}_{\{r_{k+1}=j\}}\mathds{1}_{\{r_k=i\}}\bigr] \notag\\
    &= \sum_{i\in\mathcal{S}} p_{ij}\Bigl(\check{A}_k(i)\check{S}_k(i)\check{A}_k^\top(i) + \check{A}_k(i)\check{U}_k^\top(i)\check{B}_k^\top(i) + \check{B}_k(i)\check{U}_k(i)\check{A}_k^\top(i) 
    + \check{B}_k(i)\check{Y}_k(i)\check{B}_k^\top(i)\Bigr)\label{devlifted_prop} 
\end{align}
which is linear in the lifted deviation variables $(\check{S},\check{U},\check{Y})$. \Cref{devlifted_policy_eq} is nonconvex in the moment variables. We therefore relax it to the semidefinite inequality $\check{Y}_k(i)\succeq\check{U}_k(i)\check{S}_k^{-1}(i)\check{U}_k^\top(i)$, which, for $\check{S}_k(i)\succ 0$, is equivalent by the Schur complement to the linear matrix inequality:
\begin{equation}\label{devlifted_schur}
    \begin{bmatrix} \check{S}_k(i) & \check{U}_k^\top(i) \\ \check{U}_k(i) & \check{Y}_k(i) \end{bmatrix} \succeq 0
\end{equation}
This relaxation is lossless, as discussed in \cref{sec:lossless_cvx}.
\begin{remark}
    \Cref{devlifted_prop} can be expanded block by block to recover the mean and covariance propagation equations derived earlier in \Cref{meanprop1,covprop1}. The lifted framework allows the mean and covariance quantities to be propagated (and optimized) jointly in the covariance steering framework without the need for solving the mean and covariance subproblem separately\cite{shrivastavaChanceConstrainedCovarianceSteering2025}.
\end{remark}

While \cref{devlifted_prop} is an exact propagation of MJLS first and second moments, we relax it using a slack variable 
$\bm{\nu}_k(j)\in\mathbb{R}^{n_x}$ for the mean dynamics as it represents linearized dynamics in a SCP framework. The corresponding slack matrix $N_k(j)\triangleq\begin{bsmallmatrix}\mathbf{0}_{n_x\times n_x} & \bm{\nu}_k(j)\\ \bm{\nu}_k^\top(j) & 0\end{bsmallmatrix}$ is incorporated as: 
\begin{equation}\label{devlifted_prop_nu}
    \check{S}_{k+1}(j) = \sum_{i\in\mathcal{S}} p_{ij}\Bigl(\check{A}_k(i)\check{S}_k(i)\check{A}_k^\top(i) + \check{A}_k(i)\check{U}_k^\top(i)\check{B}_k^\top(i) + \check{B}_k(i)\check{U}_k(i)\check{A}_k^\top(i) + \check{B}_k(i)\check{Y}_k(i)\check{B}_k^\top(i)\Bigr) + N_k(j) 
\end{equation}
The slacks $\{\bm{\nu}_k(j)\}_{k\in\mathbb{Z}_{0:T-1},j\in\mathcal{S}}$ are penalized in the objective and vanish as the SCP converges.

\subsubsection{Control Chance Constraint}

The control chance constraint in \cref{prob2} is given by:
\begin{equation}\label{cc_control_mjls}
    \mathbb{P}(\|\mathbf{u}_k(r_k=i)\| \leq u_{\max}) \geq 1 - \epsilon_u
\end{equation}

We now derive a deterministic approximation of this chance constraint in lifted variables. Our previous work \cite{shrivastavaChanceConstrainedCovarianceSteering2025} proposed a deterministic bound using the multivariate Chebyshev inequality. However, we found it to be too conservative for astrodynamics applications, and we instead propose to use a Gaussian-based approximation to the control chance constraint. We note that the state and therefore the control are not Gaussian distributed, and even in the linearized case form an exponential number of Gaussians. But given that the MJS in this work has only two modes with high probability of staying in nominal mode, we can approximate the control by a Gaussian distribution. To achieve this,  we seek to minimize the Kullback-Leibler divergence \cite{kullbackInformationTheoryStatistics1997} between the MJLS distribution and the approximating Gaussian. It has been shown previously \cite{minkaExpectationPropagationApproximate2001} that this is achieved by performing moment matching, i.e. setting the first two moments of the approximating 
Gaussian to be the same as the original distribution which is being approximated. We leverage this useful property to approximate the MJLS control using Gaussian distributions. For a control law defined in \cref{controllaw}, we have $\mathbb{E}[\mathbf{u}_k(i)|r_k=i]=\overline{\mathbf{u}}_k(i)$ and $\mathbb{E}[(\mathbf{u}_k(i)-\overline{\mathbf{u}}_k(i))(\mathbf{u}_k(i)-\overline{\mathbf{u}}_k(i))^\top|r_k=i]=\check{Y}_k(i)/\rho_k(i)-\delta\overline{\mathbf{u}}_k(i)\,\delta\overline{\mathbf{u}}_k^\top(i)$. Therefore, $\mathbf{u}_k(i)$ can be approximated by $\mathcal{N}(\overline{\mathbf{u}}_k(i),\check{Y}_k(i)/\rho_k(i)-\delta\overline{\mathbf{u}}_k(i)\,\delta\overline{\mathbf{u}}_k^\top(i))$. A deterministic approximation of \cref{cc_control_mjls} using \cref{det_cc_norm} is given by:

\begin{lemma}[Lemma 3 in Reference~\citenum{oguriChanceConstrainedControlSafe2024}]\label{det_cc_norm}
    Let ${\mathbf{u}_k(r_k=i)} \sim \mathcal{N}(\bar{\mathbf{u}}_k(r_k=i), \check{Y}_k(i)/\rho_k(i)-\delta\overline{\mathbf{u}}_k(i)\,\delta\overline{\mathbf{u}}_k^\top(i)) \in \mathbb{R}^{n_u}$, $\epsilon \in (0, 1)$, and $u_{\max} > 0$. Then, $\mathbb{P}[\|{\mathbf{u}_k(r_k=i)}\| \leq u_{\max}] \geq 1 - \epsilon$ is implied by
\begin{equation}\label{det_cc_norm_eq}
\|\bar{\mathbf{u}}_k(r_k=i)\| + m_{\chi^2}(\epsilon, n_u) \sqrt{\lambda_{\max}\left(\frac{Y_k(r_k=i)}{\rho_k(i)}\right)} \leq u_{\max} 
\end{equation}
where $m_{\chi^2}(\epsilon,n_u)$ denotes the square root of the $(1-\epsilon)$-quantile  of the chi-squared distribution with $n_u$ degrees of freedom.
\end{lemma}

\noindent{Note that while \cref{det_cc_norm_eq} is linear in the mean variable $\bar{\mathbf{u}}_k(r_k=i)$, it is nonconvex in the covariance terms $Y_k(r_k=i)$.} In deviation coordinates the conditional control mean is $\overline{\mathbf{u}}_k(1)=\overline{\mathbf{u}}^{\,\mathrm{ref}}_k(1)+\delta\overline{\mathbf{u}}_k(1)$ and the conditional covariance is $\check{Y}_k(1)/\rho_k(1)-\delta\overline{\mathbf{u}}_k(1)\delta\overline{\mathbf{u}}_k^\top(1)$. We convexify \cref{det_cc_norm_eq} by rearranging to isolate the square root on the left-hand side, squaring both sides, and expanding the square to get: 
\begin{equation}\label{cc_u_square}
    \begin{aligned}
        &m^2_{\chi^2}(\epsilon_u,n_u)\,\lambda_{\max}\Bigl(\frac{\check{Y}_k(1)}{\rho_k(1)}-\delta\overline{\mathbf{u}}_k(1)\delta\overline{\mathbf{u}}_k^\top(1)\Bigr) + 2u_{\max}\bigl\|\overline{\mathbf{u}}^{\,\mathrm{ref}}_k(1)+\delta\overline{\mathbf{u}}_k(1)\bigr\| \\
        &\quad - \bigl\|\overline{\mathbf{u}}^{\,\mathrm{ref}}_k(1)\bigr\|^2 - 2\,\overline{\mathbf{u}}^{\,\mathrm{ref}\top}_k(1)\delta\overline{\mathbf{u}}_k(1) - \bigl\|\delta\overline{\mathbf{u}}_k(1)\bigr\|^2 - u^2_{\max} \leq 0 .
    \end{aligned}
\end{equation}
The $\delta\overline{\mathbf{u}}_k(1)\delta\overline{\mathbf{u}}_k^\top(1)$ term inside $\lambda_{\max}(\cdot)$ and $\bigl\|\delta\overline{\mathbf{u}}_k(1)\bigr\|^2$ terms make \cref{cc_u_square} a nonconvex constraint. However, we can drop them and it will remain a sufficient condition for \cref{cc_u_square} since they are nonpositive terms. To ease SCP convergence, we also relax the resulting constraint with a penalized nonnegative slack $\zeta^{\mathrm{cc}}_k\geq0$, which gives the following convex constraint, a sufficient condition for \cref{cc_u_square} when $\zeta^{\mathrm{cc}}_k=0$:
\begin{equation}\label{devlifted_cc}
    \frac{m^2_{\chi^2}(\epsilon_u,n_u)}{\rho_k(1)}\,\lambda_{\max}\bigl(\check{Y}_k(1)\bigr) \leq u^2_{\max} - 2u_{\max}\bigl\|\overline{\mathbf{u}}^{\,\mathrm{ref}}_k(1)+\delta\overline{\mathbf{u}}_k(1)\bigr\| + \bigl\|\overline{\mathbf{u}}^{\,\mathrm{ref}}_k(1)\bigr\|^2 + 2\,\overline{\mathbf{u}}^{\,\mathrm{ref}\top}_k(1)\delta\overline{\mathbf{u}}_k(1) + \zeta^{\mathrm{cc}}_k
\end{equation}
Squaring \cref{det_cc_norm_eq} is valid only when the RHS is positive, i.e. the mean control itself satisfies the bound, so we impose \cref{devlifted_cc} along with the following constraints:
\begin{equation}\label{devlifted_cc_valid}
    \bigl\|\overline{\mathbf{u}}^{\,\mathrm{ref}}_k(1)+\delta\overline{\mathbf{u}}_k(1)\bigr\| \leq u_{\max}, \qquad \zeta^{\mathrm{cc}}_k\geq0
\end{equation}
The conservatism of \cref{devlifted_cc} is of order $\|\delta\overline{\mathbf{u}}_k(1)\|^2$, which collapses to below tolerance values as the SCP converges. 

\subsubsection{Mid-course Covariance Bound}

To enforce the mid-course covariance bound from \cref{cov_cap2}, we first express the state covariance $\Sigma_k$ in terms of the lifted deviation-based state variable. We know $\sum_{i\in\mathcal{S}}E_x\check{S}_k(i)E_x^\top=\mathbb{E}\bigl[\delta\mathbf{x}_k\delta\mathbf{x}_k^\top\bigr]=\Sigma_k+\delta\bm{\mu}_k\delta\bm{\mu}_k^\top$. Therefore,
\begin{equation}\label{dev_M}
    \Sigma_k = \sum_{i\in\mathcal{S}}E_x\check{S}_k(i)E_x^\top - \delta\bm{\mu}_k\delta\bm{\mu}_k^\top, \qquad \delta\bm{\mu}_k = \sum_{i\in\mathcal{S}}E_x\check{S}_k(i)\,\mathbf{e}
\end{equation}
The mid-course covariance bound constraint can then be given by: 
\begin{equation}\label{mid_cov_dev}
    \sum_{i\in\mathcal{S}}E_x\check{S}_k(i)E_x^\top - \left(\sum_{i\in\mathcal{S}}E_x\check{S}_k(i)\,\mathbf{e}\right)\left(\sum_{i\in\mathcal{S}}E_x\check{S}_k(i)\,\mathbf{e}\right)^\top\preceq \Sigma_{\max}
\end{equation}
\Cref{mid_cov_dev} is a nonconvex constraint because of the quadratic term on the LHS. However, since it is negative semidefinite, we can drop it to get the sufficient condition: 
\begin{equation}\label{devlifted_cap}
    \sum_{i\in\mathcal{S}}E_x\check{S}_k(i)E_x^\top \preceq \Sigma_{\max}
\end{equation}
Similar to control chance constraints, the conservatism in \cref{devlifted_cap} is of order $\|\delta\bm{\mu}_k\|^2$, which vanishes as the SCP converges.

\subsubsection{Keep-out Zone Chance Constraint}\label{sec:ko_cc}

For the keep-out zone chance constraint \cref{cc_ko2}, we first consider the feasible set: 
\begin{equation}\label{ko_nonconvex}
    g_{\mathrm{ko}}(\mathbf{x}_k) \triangleq R_{\mathrm{safe}} - \|E_r\mathbf{x}_k-\mathbf{r}_{m}\| \leq 0
\end{equation}
\Cref{ko_nonconvex} is a nonconvex constraint, and we linearize it about the reference $\bm{\mu}^{\mathrm{ref}}_k$. Since $g_{\mathrm{ko}}(\mathbf{x}_k)$ is concave, the linearization would be a global over-approximator of $g_{\mathrm{ko}}(\mathbf{x}_k)$ for all $\mathbf{x}_k$
\begin{equation}\label{ko_linearized}
    g_{\mathrm{ko}}(\mathbf{x}_k) \leq g_{\mathrm{ko}}(\bm{\mu}^{\mathrm{ref}}_k) + \nabla g_{\mathrm{ko}}^\top(\bm{\mu}^{\mathrm{ref}}_k)\bigl(\mathbf{x}_k-\bm{\mu}^{\mathrm{ref}}_k\bigr), \qquad
    \nabla g_{\mathrm{ko}}(\bm{\mu}^{\mathrm{ref}}_k) = -E_r^\top\,\frac{E_r\bm{\mu}^{\mathrm{ref}}_k-\mathbf{r}_{m}}{\bigl\|E_r\bm{\mu}^{\mathrm{ref}}_k-\mathbf{r}_{m}\bigr\|}
\end{equation}
Therefore, the following is a sufficient condition for \cref{cc_ko2}:
\begin{equation}\label{cc_ko_mjls}
    \mathbb{P}\Bigl(g_{\mathrm{ko}}(\bm{\mu}^{\mathrm{ref}}_k) + \nabla g_{\mathrm{ko}}^\top(\bm{\mu}^{\mathrm{ref}}_k)\bigl(\mathbf{x}_k-\bm{\mu}^{\mathrm{ref}}_k\bigr) \leq 0\Bigr) \geq 1-\epsilon_{\mathrm{ko}}
\end{equation}
We note that the feasible set is now convex, and under the Gaussian approximation of the state $\mathbf{x}_k$ by $\mathcal{N}(\bm{\mu}_k,\Sigma_k)$, the deterministic reformulation of the keep-out zone chance
constraint is given by: 
\begin{lemma}\cite{onoProbabilisticPlanningContinuous2013}\label{det_cc_affine}
    Let $\mathbf{x} \sim \mathcal{N}(\bm{\mu}, \Sigma) \in \mathbb{R}^{n_x}$, $\mathbf{a}\in\mathbb{R}^{n_x}$, $b\in\mathbb{R}$, and $\epsilon \in (0, 0.5)$. Then, $\mathbb{P}[\mathbf{a}^\top\mathbf{x} \leq b] \geq 1 - \epsilon$ is equivalent to
\begin{equation}\label{det_cc_affine_eq}
    \mathbf{a}^\top\bm{\mu} + m_{\mathcal{N}}(\epsilon)\sqrt{\mathbf{a}^\top\Sigma\,\mathbf{a}} \leq b
\end{equation}
where $m_{\mathcal{N}}(\epsilon)$ denotes the $(1-\epsilon)$-quantile of the standard normal distribution.
\end{lemma}
Applying \cref{det_cc_affine} to \cref{cc_ko_mjls} with $\mathbf{a}=\nabla g_{\mathrm{ko}}(\bm{\mu}^{\mathrm{ref}}_k)$ and $\delta\bm{\mu}_k=\bm{\mu}_k-\bm{\mu}^{\mathrm{ref}}_k$ gives the deterministic reformulation:
\begin{equation}\label{ko_det}
    \nabla g_{\mathrm{ko}}^\top(\bm{\mu}^{\mathrm{ref}}_k)\delta\bm{\mu}_k+ m_{\mathcal{N}}(\epsilon_{\mathrm{ko}})\sqrt{\nabla g_{\mathrm{ko}}^\top(\bm{\mu}^{\mathrm{ref}}_k)\Sigma_k\,\nabla g_{\mathrm{ko}}(\bm{\mu}^{\mathrm{ref}}_k)} \leq -g_{\mathrm{ko}}(\bm{\mu}^{\mathrm{ref}}_k) 
\end{equation}
\Cref{ko_det} is a nonconvex constraint, so we rearrange and square both sides, which is valid when the right-hand side is nonnegative, to get the equivalent form, where the second inequality is the validity condition of the squaring:
\begin{equation}\label{ko_squared}
    m^2_{\mathcal{N}}(\epsilon_{\mathrm{ko}})\,\nabla g_{\mathrm{ko}}^\top(\bm{\mu}^{\mathrm{ref}}_k)\,\Sigma_k\,\nabla g_{\mathrm{ko}}(\bm{\mu}^{\mathrm{ref}}_k) \leq \bigl(g_{\mathrm{ko}}(\bm{\mu}^{\mathrm{ref}}_k) + \nabla g_{\mathrm{ko}}^\top(\bm{\mu}^{\mathrm{ref}}_k)\delta\bm{\mu}_k\bigr)^2, \qquad g_{\mathrm{ko}}(\bm{\mu}^{\mathrm{ref}}_k) + \nabla g_{\mathrm{ko}}^\top(\bm{\mu}^{\mathrm{ref}}_k)\delta\bm{\mu}_k \leq 0
\end{equation}
Substituting $\Sigma_k$ from \cref{dev_M} and rearranging the terms to left-hand side, the first inequality in \cref{ko_squared} is equivalent to
\begin{equation}\label{ko_squared_dev}
    \begin{aligned}
        &m^2_{\mathcal{N}}(\epsilon_{\mathrm{ko}})\sum_{i\in\mathcal{S}}\nabla g_{\mathrm{ko}}^\top(\bm{\mu}^{\mathrm{ref}}_k) E_x\check{S}_k(i)E_x^\top \nabla g_{\mathrm{ko}}(\bm{\mu}^{\mathrm{ref}}_k) - g_{\mathrm{ko}}(\bm{\mu}^{\mathrm{ref}}_k)^2 - 2\,g_{\mathrm{ko}}(\bm{\mu}^{\mathrm{ref}}_k)\,\nabla g_{\mathrm{ko}}^\top(\bm{\mu}^{\mathrm{ref}}_k)\delta\bm{\mu}_k \\
        &\quad - \bigl(1+m^2_{\mathcal{N}}(\epsilon_{\mathrm{ko}})\bigr)\bigl(\nabla g_{\mathrm{ko}}^\top(\bm{\mu}^{\mathrm{ref}}_k)\delta\bm{\mu}_k\bigr)^2 \leq 0
    \end{aligned}
\end{equation}
The last term in \cref{ko_squared_dev} makes it a nonconvex constraint. However, as done in previous sections, we ignore this term since it is a nonpositive term, and the resulting constraint would be sufficient to satisfy \cref{ko_det}. We also relax the resulting constraint and the validity condition by nonnegative slacks $\zeta^{\mathrm{ko},1}_k$ and $\zeta^{\mathrm{ko},2}_k$, respectively: 
\begin{subequations}\label{devlifted_ko}
\begin{align}
    m^2_{\mathcal{N}}(\epsilon_{\mathrm{ko}})\sum_{i\in\mathcal{S}}\nabla g_{\mathrm{ko}}^\top(\bm{\mu}^{\mathrm{ref}}_k) E_x\check{S}_k(i)E_x^\top \nabla g_{\mathrm{ko}}(\bm{\mu}^{\mathrm{ref}}_k) &\leq g_{\mathrm{ko}}(\bm{\mu}^{\mathrm{ref}}_k)^2 + 2\,g_{\mathrm{ko}}(\bm{\mu}^{\mathrm{ref}}_k)\,\nabla g_{\mathrm{ko}}^\top(\bm{\mu}^{\mathrm{ref}}_k)\delta\bm{\mu}_k + \zeta^{\mathrm{ko},1}_k \label{devlifted_ko_main}\\
    g_{\mathrm{ko}}(\bm{\mu}^{\mathrm{ref}}_k) + \nabla g_{\mathrm{ko}}^\top(\bm{\mu}^{\mathrm{ref}}_k)\delta\bm{\mu}_k &\leq \zeta^{\mathrm{ko},2}_k, \qquad \zeta^{\mathrm{ko},1}_k,\,\zeta^{\mathrm{ko},2}_k\geq0 \label{devlifted_ko_valid}
\end{align}
\end{subequations}
where $\nabla g_{\mathrm{ko}}$ is given in \cref{ko_linearized}. The conservatism introduced by ignoring the nonconvex term is of the order $\|\delta\bm{\mu}_k\|^2$, which vanishes as the SCP converges. 

\subsubsection{Boundary Conditions in Deviation Variables}
With $\delta\bm{\mu}_0\triangleq\overline{\bm{\mu}}_0-\bm{\mu}^{\mathrm{ref}}_0$, the second moment initial condition is given as:
\begin{equation}\label{devlifted_bc0}
    \check{S}_0(i) = \rho_0(i)\begin{bmatrix} \overline{\Sigma}_0 + \delta\bm{\mu}_0\delta\bm{\mu}_0^\top & \delta\bm{\mu}_0 \\ \delta\bm{\mu}_0^\top & 1 \end{bmatrix}
\end{equation}
Similarly, let $\delta\bm{\mu}_f\triangleq\bm{\mu}_f-\bm{\mu}^{\mathrm{ref}}_T$, then the terminal conditions are given as:
\begin{equation}\label{devlifted_bcf}
    \sum_{i\in\mathcal{S}}E_x\check{S}_T(i)\,\mathbf{e} = \delta\bm{\mu}_f, \qquad
    \sum_{i\in\mathcal{S}}E_x\check{S}_T(i)E_x^\top \preceq \Sigma_f + \delta\bm{\mu}_f\delta\bm{\mu}_f^\top, \qquad
    \check{S}_T(i)\succeq 0
\end{equation}
which are all convex  in $\check{S}_T(i)$. The first two are jointly equivalent to $\mathbb{E}[\mathbf{x}_T]=\bm{\mu}_f$ and $\Sigma_T\preceq\Sigma_f$. The last condition ensures that the terminal lifted moment is a valid second moment.

\subsubsection{Objective Function}\label{sec:dev_objective}

We wish to optimize the $\Delta v_{99}$ cost in \cref{prob1}. Since the deterministic formulation for chance constraints used a Gaussian approximation for control, we use the same approximation for a deterministic upper bound of the $\Delta v_{99}$. For a Gaussian control policy, a deterministic upper bound for $\Delta v_{99}$ was given in Ref~\cite{oguriChanceConstrainedControlSafe2024}. However, it involves a nonconvex term for the feedback cost. In this work, we use an approximate expression for the $\Delta v_{99}$ cost function: 

\begin{equation}\label{dev_objective}
    J = \sum_{i\in\mathcal{S}}\sum_{k=0}^{T-1}\Bigl[\rho_k(i)\,\bigl\|\overline{\mathbf{u}}_k(i)\bigr\| + m^2_{\chi^2}(0.01,n_u)\,\rho_k(i)\,\mathrm{tr}\bigl(cov(\mathbf{u}_k\,|\,r_k=i)\bigr)\Bigr]
\end{equation}

\Cref{dev_objective} utilizes the trace of the control covariance instead of its maximum eigenvalue. However, they have a similar effect of controlling the dispersion, and have previously been used in other covariance steering frameworks\cite{kumagaiSequentialChanceConstrainedCovariance2024}. In the current lifted deviation-based state formulation, the control policy can be expressed as:
\begin{equation}\label{dev_u_gauss}
    \mathbf{u}_k\,|\,r_k=i\sim\mathcal{N}\biggl(\overline{\mathbf{u}}^{\,\mathrm{ref}}_k(i)+\frac{\check{U}_k(i)\,\mathbf{e}}{\rho_k(i)},\;\frac{\check{Y}_k(i)}{\rho_k(i)}-\frac{\check{U}_k(i)\,\mathbf{e}\,\mathbf{e}^\top\check{U}_k^\top(i)}{\rho_k^2(i)}\biggr)
\end{equation}
\Cref{dev_u_gauss} results in a nonconvex objective function because of the quadratic term in the covariance. Since it is nonpositive, we drop it, and its conservatism will vanish as the SCP converges. Therefore, we have the convex upper bound on $\Delta v_{99}$ cost function: 
\begin{equation}\label{dev_objective_lifted}
    \check{J}\bigl(\check{U},\check{Y}\bigr) = \sum_{i\in\mathcal{S}}\sum_{k=0}^{T-1}\Bigl[\bigl\|\rho_k(i)\,\overline{\mathbf{u}}^{\,\mathrm{ref}}_k(i)+\check{U}_k(i)\,\mathbf{e}\bigr\| + m^2_{\chi^2}(0.01,n_u)\,\mathrm{tr}\bigl(\check{Y}_k(i)\bigr)\Bigr]
\end{equation}

\begin{remark}
    We note that the deviation-based formulation allows us to drop the nonconvex nonpositive terms in a nonconvex constraint, without introducing additional conservatism at optimality. The deviation-based framework ensures that those terms vanish as the SCP converges, implying that the sufficient conditions used in this work do not introduce any conservatism at optimality, while avoiding the need for the linearization-based convexification steps standard in sequential convex programming frameworks\cite{kumagaiRobustCislunarLowThrust,oguriSuccessiveConvexificationFeasibility}.
\end{remark}

\subsubsection{Losslessness of the Policy Relaxation}\label{sec:lossless_cvx}

In \cref{devlifted_schur}, we propose a relaxation of the nonconvex constraint $\check{Y}_k(i)=\check{U}_k(i)\check{S}_k^{-1}(i)\check{U}_k^\top(i)$. Such relaxations have been shown to be lossless at optimality for both single-mode covariance steering problems \cite{liuOptimalCovarianceSteering2024,kumagaiRobustCislunarLowThrust} and MJLS covariance steering problems \cite{shrivastavaChanceConstrainedCovarianceSteering2025,wangCovarianceSteeringDiscreteTime2026}, as long as the objective function \cref{dev_objective_lifted} is strictly increasing in $\check{Y}_k(i)$, i.e. $\frac{\partial J}{\partial \check{Y}_k(i)}\succ 0$. Under the current $\Delta v_{99}$ cost function in \cref{dev_objective_lifted}, we satisfy this condition and therefore the relaxation is lossless at optimality for each convex subproblem.

\subsubsection{Convex Subproblem}
We have now proposed convex deterministic reformulations for all nonconvex, stochastic constraints in \cref{prob2}. Therefore, the convex subproblem is given by: 

\begin{problem}\label{prob_devlifted} Given a reference trajectory $\{\bm{\mu}^{\mathrm{ref}}_k, \overline{\mathbf{u}}^{\,\mathrm{ref}}_k(i)\}$, solve:
    \begin{subequations}\label{ocp_devlifted}
    \begin{flalign}
    &\min_{\substack{\check{S}_k(i),\,\check{U}_k(i),\,\check{Y}_k(i)\\ \bm{\nu},\,\bm{\zeta}}}\;\; \check{J}\bigl(\check{U},\check{Y}\bigr) + P_{\mathrm{cc}}(\bm{\zeta}) + P_{\bm{\nu}}(\bm{\nu}) && \clabel{$\Delta v_{99}$ cost (\cref{dev_objective_lifted}) with slack penalty} \label{eq:devlifted_obj}\\
    &\text{s.t.}\;\;\;\forall k \in \mathbb{Z}_{0:T-1},\; i,j\in\mathcal{S}, \nonumber\\
    &\quad \text{\Cref{devlifted_bc0,devlifted_bcf}} && \clabel{Boundary conditions} \label{eq:devlifted_bc}\\
    &\quad \text{\Cref{devlifted_prop_nu}}, \qquad \forall k\in\mathbb{Z}_{0:T-1},\; j\in\mathcal{S} && \clabel{Lifted moment propagation} \label{eq:devlifted_dyn}\\
    &\quad \text{\Cref{devlifted_schur}} && \clabel{Lossless relaxation of $Y_k(i)$} \label{eq:devlifted_lmi}\\
    &\quad \text{\Cref{devlifted_cc,devlifted_cc_valid}} && \clabel{Control chance constraint} \label{eq:devlifted_cc}\\
    &\quad \text{\Cref{devlifted_cap}}, \qquad \forall k\in\mathbb{Z}_{1:T-1} && \clabel{Mid-course covariance bound} \label{eq:devlifted_cap}\\
    &\quad \text{\Cref{devlifted_ko}}, \qquad \forall k\in\mathbb{Z}_{1:T-1} && \clabel{Keep-out zone chance constraint} \label{eq:devlifted_ko}\\
    &\quad \frac{1}{\kappa_x}\bigl\|\delta\bm{\mu}_k\bigr\|_\infty \leq \Delta_{\mathrm{TR}}, \qquad \frac{1}{\kappa_u}\bigl\|\check{U}_k(1)\,\mathbf{e}\bigr\|_\infty \leq \Delta_{\mathrm{TR}} && \clabel{Trust region constraint} \label{eq:devlifted_tr}
    \end{flalign}
    \end{subequations}
where $\bm{\zeta}$ stacks all chance-constraint slacks, $\bm{\nu}$ stacks the mean-dynamics slacks $\bm{\nu}_k(j)$, $P_{\mathrm{cc}}$, $P_{\bm{\nu}}$ are the corresponding penalty functions defined in \cref{sec:scp_scheme}, and $\kappa_x,\kappa_u>0$ are user-specified scaling constants for trust region constraints.
\end{problem}
Once \cref{prob_devlifted} is solved, the controller of \cref{controllaw} is recovered from the lifted variables via
\begin{equation}\label{devlifted_gain_recovery}
    \check{K}_k(i)=\check{U}_k(i)\check{S}_k^{-1}(i)=\bigl[K_k(i),\;\delta\overline{\mathbf{u}}_k(i)-K_k(i)\,\delta\overline{\mathbf{x}}_k(i)\bigr]
\end{equation} and the true control policy can be computed as $\overline{\mathbf{u}}_k(i)=\overline{\mathbf{u}}^{\,\mathrm{ref}}_k(i)+\delta\overline{\mathbf{u}}_k(i)$ and $\overline{\mathbf{x}}_k(i)=\bm{\mu}^{\mathrm{ref}}_k+\delta\overline{\mathbf{x}}_k(i)$ with $\delta\overline{\mathbf{x}}_k(i)=E_x\check{S}_k(i)\,\mathbf{e}/\rho_k(i)$. 

\subsection{SCP framework for Nonlinear Markov Jump Covariance Steering}\label{SCP_framework}

We now describe the SCP framework used to solve \cref{prob1}. We employ the SCvx$\ast$ framework \cite{oguriSuccessiveConvexificationFeasibility} to update the reference dynamics and SCP parameters. Each SCvx$\ast$ iteration solves the semidefinite program given in \cref{prob_devlifted}, and the algorithm is summarized in \Cref{alg:masto}. We call our algorithm MASTO, which stands for \emph{Missed-thrust Aware Stochastic Trajectory Optimization}. We discuss the details on initial guess generation, step acceptance and trust region update rules, and numerical scaling below.

\begin{algorithm}[H]
    \caption{The MASTO Algorithm}
    \label{alg:masto}
    \setstretch{1}
    \begin{algorithmic}[1]
        \State Initialize the reference trajectory $\{\bm{\mu}^{\mathrm{ref}}_k,\overline{\mathbf{u}}^{\,\mathrm{ref}}_k(i)\}$ (\Cref{sec:init_guess}) and the SCP parameters
        \While{iterations $\leq$ max iterations}
            \State Linearize about $\{\bm{\mu}^{\mathrm{ref}}_k,\overline{\mathbf{u}}^{\,\mathrm{ref}}_k(i)\}$ if it was updated (\Cref{sec:linearization} and \cref{ko_linearized})
            \State Solve \Cref{prob_devlifted}; compute metrics $\chi$, $\Delta J$, $\Delta L$ and $\rho$
            \If{converged}
                \State \Return the solution by recovering $K_k(i)$, $\overline{\mathbf{u}}_k(i)$ \Comment{\Cref{devlifted_gain_recovery}}
            \EndIf
            \If{$\rho\in[1-\eta_0,1+\eta_0]$} \Comment{step acceptance}
                \State Update the reference trajectory $\{\bm{\mu}^{\mathrm{ref}}_k,\overline{\mathbf{u}}^{\,\mathrm{ref}}_k(i)\}$
                \If{$|\Delta J|<\delta$} 
                    \State Update the Lagrange multipliers and SCP parameters $\bm{\lambda},\bm{\psi}$, $w$ and $\delta$ \Comment{Ref.~\citenum{oguriSuccessiveConvexificationFeasibility}}
                \EndIf
            \EndIf
            \State Update $\Delta_{\mathrm{TR}}$ \Comment{\Cref{tr_update}}
        \EndWhile
    \end{algorithmic}
\end{algorithm} 

\subsubsection{Initial Guess Trajectory}\label{sec:init_guess}

Most SCP-based spacecraft trajectory optimization works use a deterministic open-loop minimum-fuel or minimum energy trajectory for the initial guess\cite{kumagaiRobustCislunarLowThrust,oguriConvexBlockCholeskyApproach2026}. However, for this work, we generate a new initial guess with scaled control norm bound $u_{\max,k}=\rho_k(r_k=1) u_{\max}$. We use a deterministic SCvx$\ast$ framework to get this initial guess and use the deterministic open-loop minimum-fuel transfer as initialization for this SCP run. The augmented control bound better represents the actual control authority available to the Markov Jump System, and is therefore a better initial guess for the SCP framework.

\subsubsection{SCP scheme}\label{sec:scp_scheme}

For the SCP scheme, we adopt the augmented-Lagrangian cost, step acceptance and trust-region update rule proposed in Reference~\citenum{kumagaiRobustCislunarLowThrust}, which is based on the SCvx$\ast$ framework\cite{oguriSuccessiveConvexificationFeasibility}. Therefore, the penalty terms of \cref{prob_devlifted} are
\begin{equation}\label{devlifted_penalties}
    P_{\mathrm{cc}}(\bm{\zeta}) = \bm{\psi}^\top\bm{\zeta} + \sqrt{w}\,\|\bm{\zeta}\|_1 + \frac{w}{2}\|\bm{\zeta}\|^2, \qquad
    P_{\bm{\nu}}(\bm{\nu}) = \bm{\lambda}^\top\bm{\nu} + \sqrt{w}\,\|\bm{\nu}\|_1 + \frac{w}{2}\|\bm{\nu}\|^2
\end{equation}
where the $l_1$ norm terms are added to improve convergence. We note that our framework inherits the same ``inexact linearization'' challenge discussed in Reference~\citenum{kumagaiRobustCislunarLowThrust}, as the MJLS covariance propagation equation in each convex subproblem is based on the reference state dynamics and thus ignores the second-order effects. This can negatively affect the SCP convergence, and the SCvx$\ast$ convergence guarantees\cite{oguriSuccessiveConvexificationFeasibility} do not carry over to our framework. However, addressing it would require the propagation of second order state tensors, which would be computationally intensive, without significant SCP performance improvements. 

Therefore, we use the modified rule for step acceptance and trust region updates, which was previously shown to perform well for such problems\cite{kumagaiRobustCislunarLowThrust}. To determine step acceptance, we use the parameter $\rho=\Delta J/\Delta L$, which is the ratio of the actual reduction $\Delta J$ of the nonlinear augmented cost to the reduction $\Delta L$ predicted by \cref{prob_devlifted} (see Reference~\citenum{kumagaiRobustCislunarLowThrust} for details). For a given set of predefined SCP parameters $1\geq\eta_0>\eta_1> \eta_2 > 0$, $\alpha_1, \alpha_2 > 1$, and $0<\Delta_{\mathrm{TR,min}}< \Delta_{\mathrm{TR,max}}$, the step is accepted if $\rho \in [1-\eta_0,\, 1+\eta_0]$ and the trust region update rule is given by:
\begin{equation}\label{tr_update}
    \Delta_{\mathrm{TR}} \leftarrow 
    \left\{
    \begin{array}{ll}
        \min\{\alpha_2\Delta_{\mathrm{TR}},\, \Delta_{\mathrm{TR,max}}\} & \text{if } \rho \in [1-\eta_2,\, 1+\eta_2] \\
        \Delta_{\mathrm{TR}} & \text{elseif } \rho \in [1-\eta_1,\, 1+\eta_1] \\
        \max\{\Delta_{\mathrm{TR}}/\alpha_1,\, \Delta_{\mathrm{TR,min}}\} & \text{otherwise}
    \end{array}
    \right.
\end{equation}
When an accepted step also satisfies the condition $|\Delta J|<\delta$, the multipliers $\bm{\lambda},\bm{\psi},w$ and SCP parameters $\gamma$ are updated. We refer the reader to Algorithm~1 of Reference~\citenum{oguriSuccessiveConvexificationFeasibility} for more details about the SCvx$\ast$ framework.

\subsubsection{Numerical Scaling Strategy}

The deviation-based state formulation improves the scaling of the optimization problem since the mean and covariance are now of similar magnitude. However, some other scaling concerns remain. Depending on the choice of the units, the position and velocity can be of very different scales in some astrodynamics problems. Further, the lifted state formulation has an $\mathcal{O}(1)$ term which propagates the mode probability, while other terms are $\mathcal{O}(\|\delta \mathbf{x}\|^2)$ whose exact scale depends on the choice of units. To simplify the scaling and make it more robust to the choice of units used in the optimization process, we propose a static matrix scaling, given by the following transformation of coefficient matrices and boundary conditions: 
\begin{equation}\label{eq:static_scaling}
    \check{A}_k(i) \leftarrow \check{D}\,\check{A}_k(i)\,\check{D}^{-1}, \qquad
    \check{B}_k(i) \leftarrow s_u\,\check{D}\,\check{B}_k(i), \qquad
    \Sigma \leftarrow D_x\,\Sigma\,D_x, \qquad
    \delta\bm{\mu} \leftarrow D_x\,\delta\bm{\mu},
\end{equation}
for $\Sigma\in\{\overline{\Sigma}_0,\Sigma_f,\Sigma_{\max}\}$ and $\delta\bm{\mu}\in\{\delta\bm{\mu}_0,\delta\bm{\mu}_f\}$. Here, $\check{D} = \operatorname{blkdiag}(D_x, 1)$, where $D_x = \operatorname{blkdiag}(\sigma_{r_f}^{-1}I_3,\,\sigma_{v_f}^{-1}I_3)$ and $s_u = u_{\max}$. $\sigma_{r_f}$ and $\sigma_{v_f}$ correspond to the terminal covariance bound, given by $\Sigma_f = \operatorname{blkdiag}(\sigma^2_{r_f}I_3,\,\sigma^2_{v_f}I_3)$. \Cref{eq:static_scaling} corresponds to a scaling change in variables as $\delta\mathbf{x}\leftarrow D_x\,\delta\mathbf{x}$, $\delta\mathbf{u}\leftarrow\delta\mathbf{u}/s_u$, making them all roughly $\mathcal{O}(1)$. For the trust region constraint scaling constants, we use $\kappa_u=1$ and $\kappa_x=\|D_x(\bm{\mu}_f-\overline{\bm{\mu}}_0)\|_\infty$. 

\section{Numerical Example}\label{sec:results}

In this section, we demonstrate the resilience of MASTO trajectories against statistical MTEs through nonlinear Monte Carlo simulations. We consider two three-dimensional transfers: (a) a L2 Halo to L1 Halo transfer in the Earth-Moon CR3BP, and (b) a complex multi-revolution L2 Near Rectilinear Halo Orbit (NRHO) to L1 Halo transfer in the Earth-Moon CR3BP. We first describe the simulation set-up below:

\subsection{Simulation Set-up}
\subsubsection{SCP Framework Parameters}

We use the following parameters for the SCvx$\ast$ framework used in this work: 

\begin{table}[h!]
\centering
\caption{SCvx$\ast$ parameters}
\label{tab:scvx_params}
\begin{tabular}{ll}
\hline\hline
Parameter                                                    & Value                \\ \hline
$\{ \epsilon_{\text{opt}}, \epsilon_{\text{feas}} \}$         & $\{ 10^{-3}, 5\times10^{-4} \}$ \\
$\{ \eta_0, \eta_1, \eta_2 \}$                               & $\{ 1, 0.90, 0.25 \}$   \\
$\{ \alpha_1, \alpha_2 \}$                                   & $\{ 2, 2 \}$         \\
$\beta$                                                      & 2                    \\
$\gamma$                                                     & 0.5                  \\
$\{ \Delta_{\text{TR,min}}, \Delta_{\text{TR,max}}, \Delta_{\text{TR}}^{(1)} \}$ & $\{ 10^{-5}, 1, 0.5 \}$ \\
$\{ w^{(1)}, w_{\text{max}} \}$                               & $\{ 10^4, 10^6 \}$   \\ \hline\hline
\end{tabular}
\end{table}

We solve all optimization problems with CVXPY~\cite{diamondCVXPYPythonEmbeddedModeling2016}, using Clarabel~\cite{goulartClarabelInteriorpointSolver2024} as the primary solver and MOSEK~\cite{apsMOSEKOptimizerAPI} as a fallback. In our tests, Clarabel was faster and solved the semidefinite programs to tighter tolerances, but it was also more brittle and failed to return a solution in some cases. In comparison, MOSEK was more robust but yielded less accurate solutions. We therefore use MOSEK only when Clarabel fails. All optimization problems and nonlinear Monte Carlo simulations were run on a 2023 MacBook Pro with 16~GB of RAM.

\subsubsection{Initial Guess Trajectory Generation}

In both cases, we use a previously generated initial guess trajectory using the minimum-fuel deterministic control policy. The details of optimizing and computing the deterministic minimum-fuel trajectory for these cases can be found in Refs.~\citenum{sidhoumPontryaginBellmanDifferential2026} and \citenum{kumagaiRobustCislunarLowThrust} for the L2 Halo to L1 Halo case and the L2 NRHO to L1 Halo case, respectively. In this first step of the MASTO framework, we solve the minimum-fuel deterministic optimization problem (using the SCvx$\ast$ framework with the same parameters as \Cref{tab:scvx_params}) with mode-probability scaled $u_{\max}$. For numerical simplicity, we set the initial mode probability to be the same as the stationary distribution of the corresponding Markov chain, which allows for a constant $u_{\max}$ constraint in this problem. Note that this is not a hard requirement, and any initial mode probability distribution can be used, although the convergence of the SCP framework may be affected. 

\subsubsection{Uncertainty Modeling and Path Constraints}

We sample the initial state dispersion from a Gaussian distribution with covariance ${\Sigma}_0 = \text{blkdiag}(\sigma^2_{r_0} I_3, \sigma^2_{v_0}I_3)$, where \{$\sigma_{r_0}$, $\sigma_{v_0}$\} denote the standard deviations of the position and velocity at the departure orbit. At the final time, we impose a terminal covariance constraint with ${\Sigma}_f = \text{blkdiag}(\sigma^2_{r_f} I_3, \sigma^2_{v_f}I_3)$, where \{$\sigma_{r_f}$, $\sigma_{v_f}$\} are the standard deviations of the position and velocity for target final covariance. The mid-course covariance bound is also imposed as ${\Sigma}_{\max} = \text{blkdiag}(\sigma^2_{r_{\max}} I_3, \sigma^2_{v_{\max}}I_3)$ throughout the trajectory.  

In this work, we want to focus on how the trajectory is affected by introducing robustness to MTEs, and we therefore exclude other sources of uncertainty. Thus, the mode-switching Markov chain and the dispersed initial condition are the only sources of uncertainty. Nonetheless, the MASTO framework can also extend to consider other sources of uncertainty such as maneuver execution error and dynamics modeling error, modeled as additive Gaussian noise. The Markov jump-based switching also allows for mode-dependent uncertainty. For instance, we can have increased maneuver execution error in a mode with partial thruster failure. 

We also impose a lunar surface clearance constraint, which is a nonconvex keep-out zone constraint and is convexified as explained in \Cref{sec:ko_cc}. The keep-out zone is defined at $100$ km altitude above the lunar surface, and is especially critical in the complex L2 NRHO to L1 Halo transfer case because of its close flyby. We enforce the lunar surface keep-out zone constraint with probability $1-\epsilon_{\mathrm{ko}}$ where $\epsilon_{\mathrm{ko}}=10^{-4}$. For the control chance constraints, we set $\epsilon_u=0.05$, i.e. the control bound should be satisfied with at least $95\%$ probability. In the Monte Carlo simulations, we cap the executed control by $u_{\max}$, in case the commanded control exceeds the upper bound. 

\subsubsection{Single-mode Comparison Framework}

To assess the robustness of MASTO trajectories to MTEs, we compare its performance to a non-MTE aware trajectory optimization framework. While the deterministic open-loop policy is a traditional benchmark, it is an unfair comparison as it does not account for any correction maneuvers to control the dispersion. On the other hand, closed-loop control strategies such as the full block covariance steering framework\cite{kumagaiRobustCislunarLowThrust} have emerged as a promising framework for robust spacecraft trajectory optimization problems under uncertainty. Therefore, in this work, we compare the MASTO framework with a simplified version of closed-loop single-mode covariance steering framework, which only accounts for the uncertainty introduced by the initial dispersion. We use \Cref{alg:masto} to solve the single-mode covariance steering problem, and define the corresponding Markov jump system as a single-mode system, with the same SCP parameters as \Cref{tab:scvx_params}. For the single-mode problem, we additionally impose the covariance floor $\Sigma_k\succeq\Sigma_{\min}$ at every node, with $\Sigma_{\min}=\operatorname{blkdiag}\bigl((1\,\mathrm{km})^2I_3,\,(10\,\mathrm{cm/s})^2I_3\bigr)$, which is a linear matrix inequality (LMI) in the lifted variables. Without it, the single-mode controller can collapse the planned covariance to a small value which makes it difficult to optimize, especially in the absence of other uncertainties. The MASTO framework generalizes the single-mode covariance steering framework as we can define the single-mode system as a non-switching single-mode MJS . We note that our comparison benchmark does not consider other uncertainties such as navigation uncertainty which is tackled in Reference \cite{kumagaiRobustCislunarLowThrust}.

\subsection{L2 Halo to L1 Halo Transfer in Earth-Moon CR3BP}

We first consider the L2 Halo to L1 Halo transfer in the Earth-Moon CR3BP. \Cref{tab:hh_ic} lists the initial conditions for the departure and arrival orbit, which we obtain after performing differential correction\cite{spreenTrajectoryDesignTargeting2021} to get a periodic orbit. The simulation parameters for this case with the MTE statistics are included in \Cref{tab:hh_params}. For the given $\text{TTF}_{90}$ of $30$ days and $40$ days of total time-of-flight, the spacecraft will experience at least $1$ MTE with more than $90\%$ probability during the transfer. We run \Cref{alg:masto} to get the MTE-resilient MASTO trajectory. It converges in approximately 5 minutes with $21$ iterations ($\approx15$ s per iteration on average).
\begin{table}[h!]
    \centering
    \caption{Initial conditions for L2 Halo to L1 Halo transfer ($y=\dot{x}=\dot{z}=0$)}
    \label{tab:hh_ic}
    \begin{tabular}{lrrrr}
    \hline\hline
    \multicolumn{1}{c}{Orbit} & \multicolumn{1}{c}{$x$ (ND)} & \multicolumn{1}{c}{$z$ (ND)} & \multicolumn{1}{c}{$\dot{y}$ (ND)} & \multicolumn{1}{c}{Period (days)} \\ \hline
    L2 Halo departure orbit  & 1.16               & $-0.124716843616057$ & $-0.208673496954019$ & $\approx14.3$ \\
    L1 Halo arrival orbit & 0.85               & $0.175465054221552$  & $0.262898036867705$  & $\approx11.1$ \\ \hline\hline
    \end{tabular}
\end{table}
\begin{table}[htbp!]
    \centering
    \caption{Parameters for the L2 Halo to L1 Halo transfer}
    \label{tab:hh_params}
    \begin{tabular}{lrl|lrl}
    \hline\hline
    Parameter & Value & Unit & Parameter & Value & Unit \\ \hline
    Time-of-flight           & 40   & days     & Discretization nodes & 101 & -- \\
    $u_{\max}$   & 0.3  & mm/s$^2$ & $\text{TTF}_{90}$    & 30  & days \\
    $\sigma_{r_0}$           & 5    & km       & $\text{TTR}_{95}$    & 0.5 & days \\
    $\sigma_{v_0}$           & 0.1  & m/s      & \multirow{4}{*}{$\mathcal{P}$} & \multirow{4}{*}{$\begin{bmatrix}0.96980 & 0.03020\\0.90900 & 0.09100\end{bmatrix}$} & \multirow{4}{*}{--} \\
    $\sigma_{r_f}$           & 100  & km       & & & \\
    $\sigma_{v_f}$           & 10   & m/s      & & & \\
    $\sigma_{r,\max}$        & 500  & km       & & & \\
    $\sigma_{v,\max}$        & 15   & m/s      & $\bm{\rho}_0$ & $[0.968,\ 0.032]^\top$ & -- \\ \hline\hline
    \end{tabular}
\end{table}

\begin{figure}[htbp!]
    \centering
    \begin{subfigure}[b]{0.48\textwidth}
        \centering
        \includegraphics[height=0.78\linewidth]{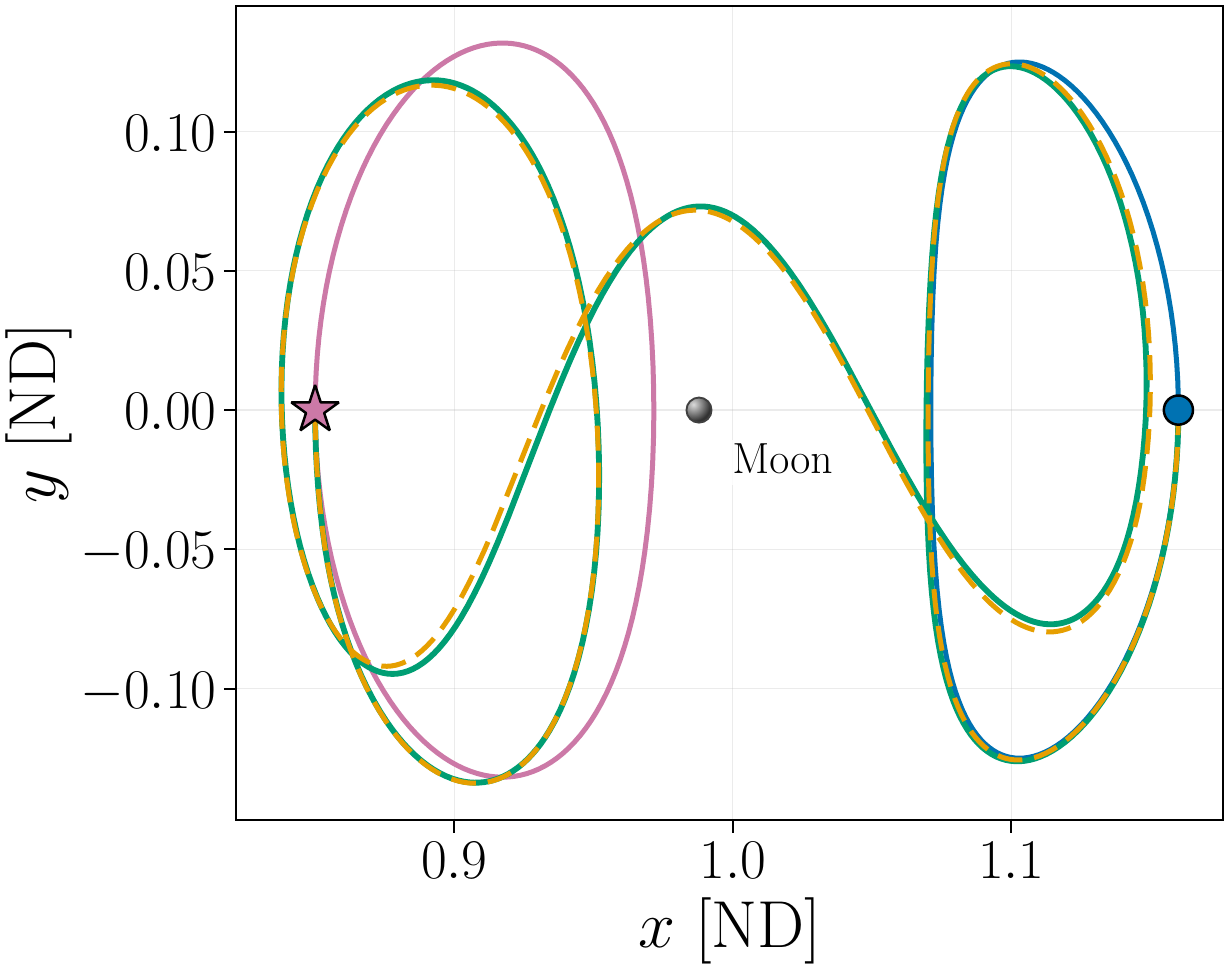}
        \caption{$xy$-projection.}
        \label{fig:hh_ovm_xy}
    \end{subfigure}
    \hfill
    \begin{subfigure}[b]{0.48\textwidth}
        \centering
        \includegraphics[height=0.78\linewidth]{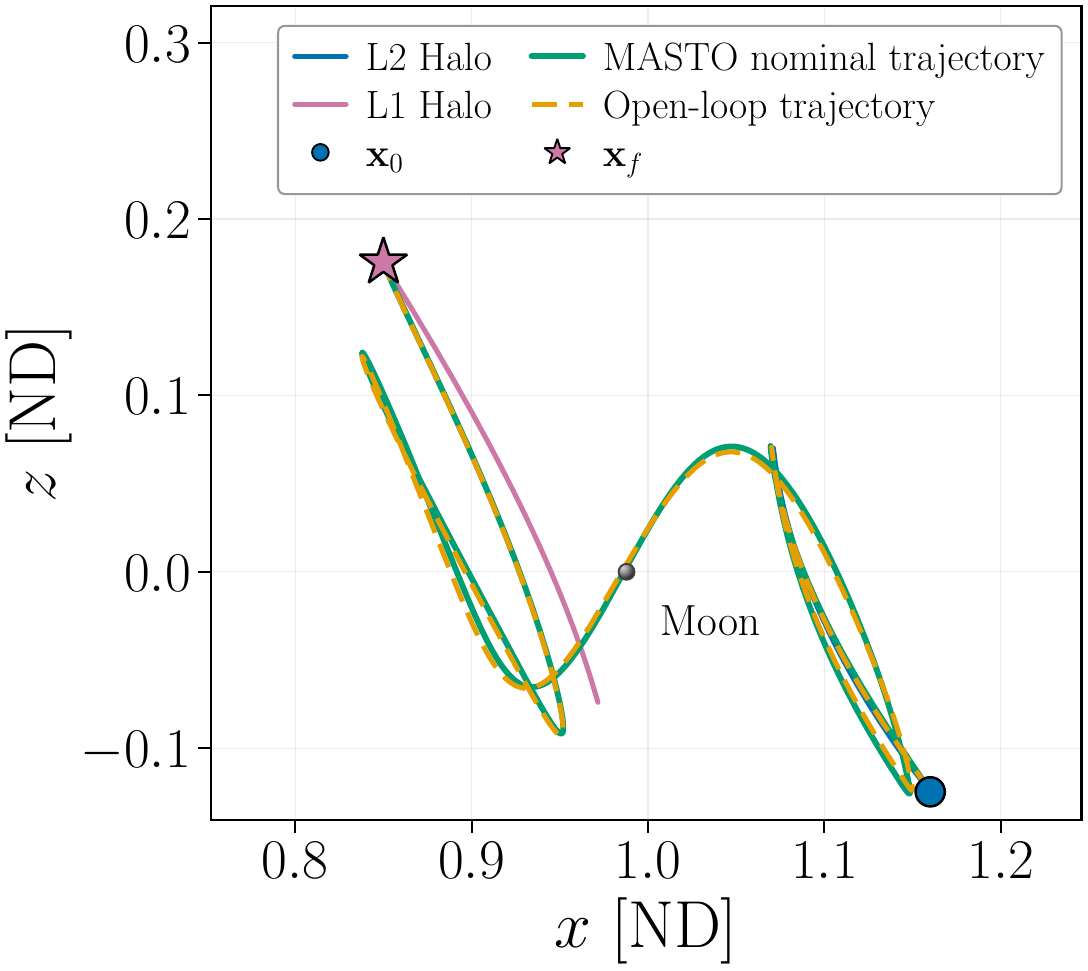}
        \caption{$xz$-projection.}
        \label{fig:hh_ovm_xz}
    \end{subfigure}

    \vspace{0.6em}

    \begin{subfigure}[b]{0.48\textwidth}
        \centering
        \includegraphics[height=0.84\linewidth]{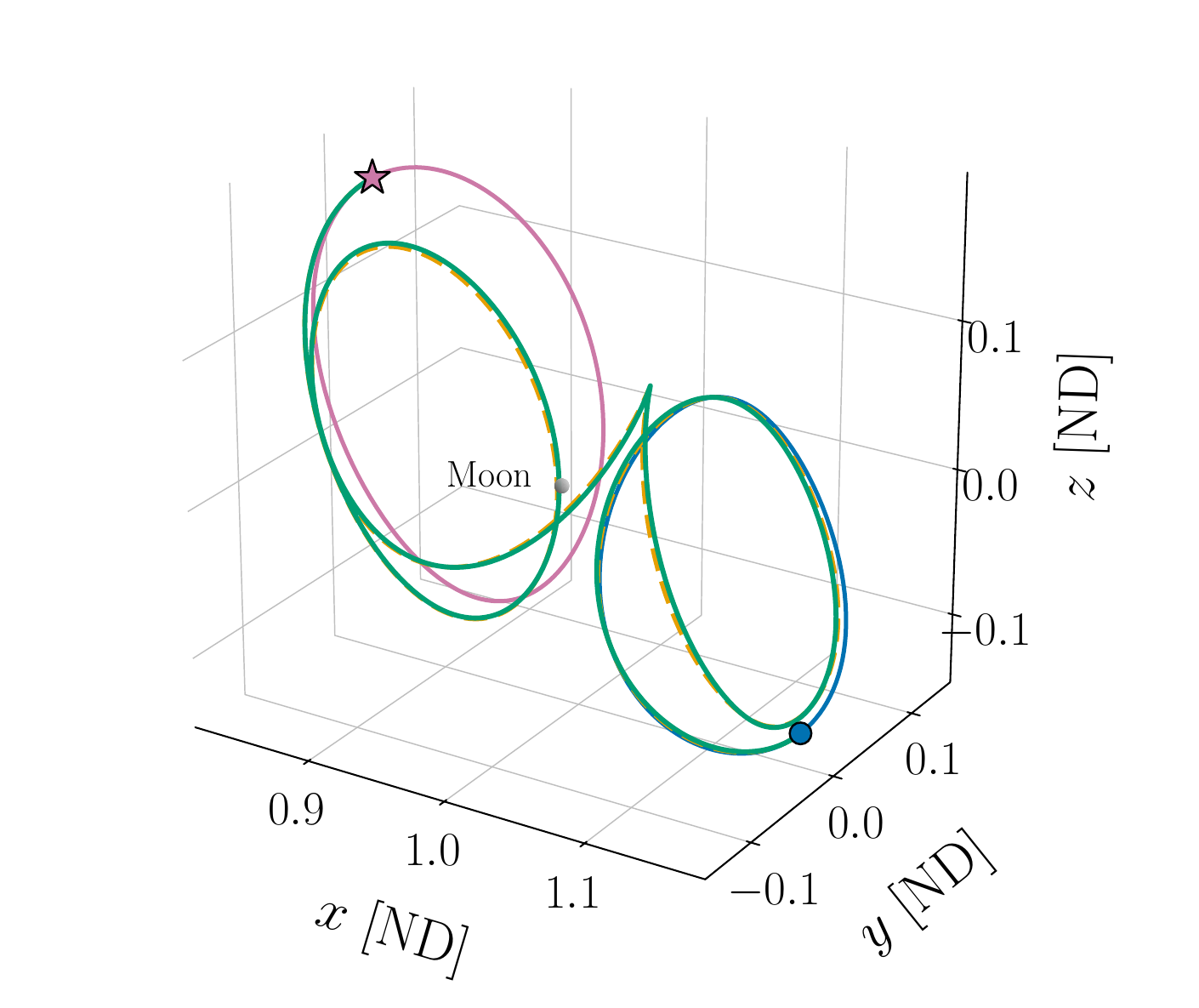}
        \caption{Three-dimensional view.}
        \label{fig:hh_ovm_3d}
    \end{subfigure}
    \hfill
    \begin{subfigure}[b]{0.48\textwidth}
        \centering
        \includegraphics[height=0.84\linewidth]{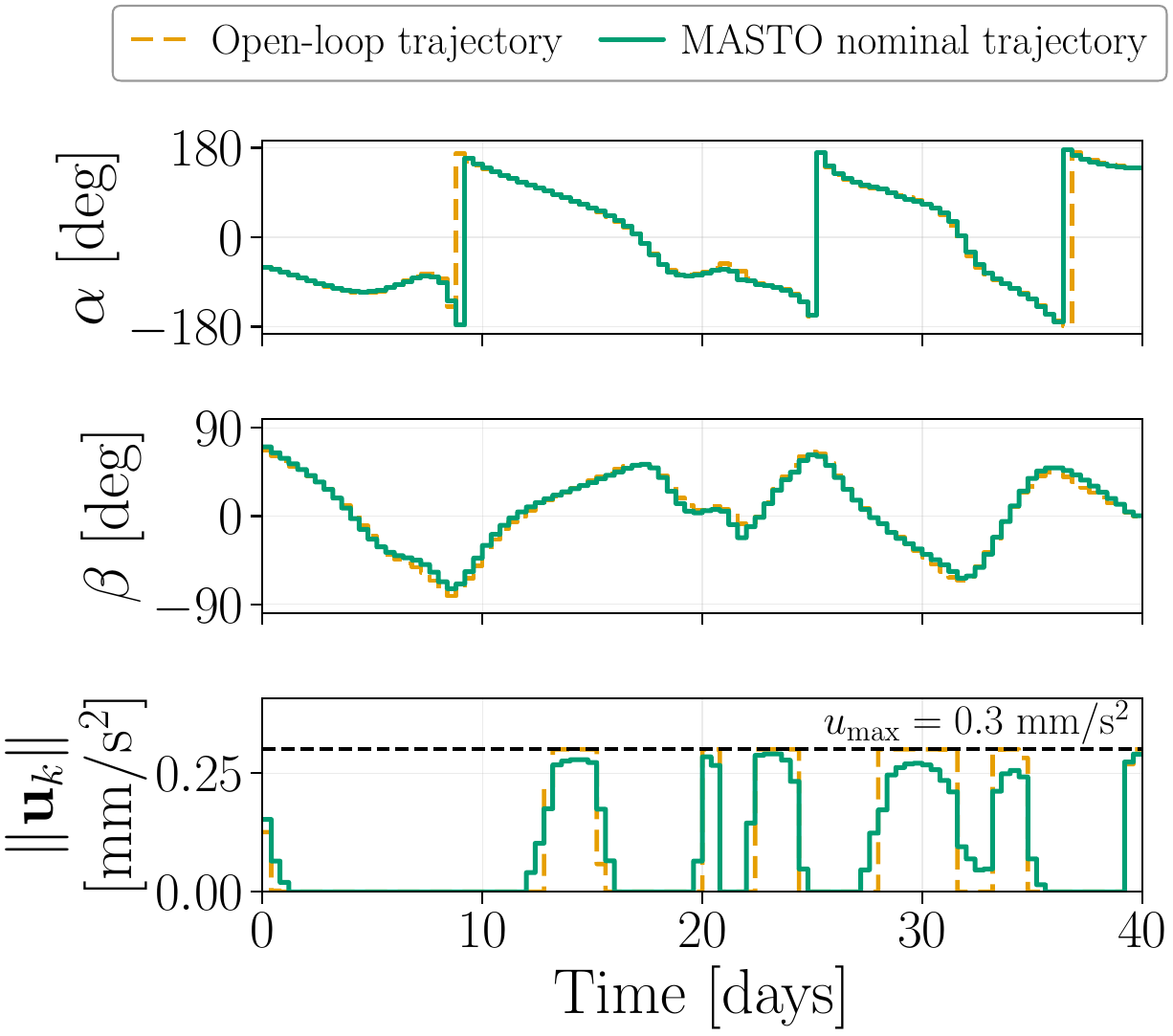}
        \caption{Control angles and magnitude.}
        \label{fig:hh_ovm_ctrl}
    \end{subfigure}
    \caption{Nominal MASTO against the open-loop minimum-fuel solution, L2 Halo to L1 Halo transfer.}
    \label{fig:hh_ovm}
\end{figure}

We first compare the MASTO nominal trajectory to the deterministic open-loop reference in \Cref{fig:hh_ovm}. \Cref{fig:hh_ovm_xy}, \Cref{fig:hh_ovm_xz}, and \Cref{fig:hh_ovm_3d} compare the $xy$, $xz$ projections, and the three-dimensional trajectory respectively, while \Cref{fig:hh_ovm_ctrl} compares the respective control profiles, with angles $\alpha$ and $\beta$ defined as
\begin{equation}\label{control_angles} 
    \alpha_k = \operatorname{atan2}\bigl(u_{y,k},\,u_{x,k}\bigr) \in (-180^\circ,\,180^\circ], \qquad
    \beta_k = \arcsin\biggl(\frac{u_{z,k}}{\|\mathbf{u}_k\|}\biggr) \in [-90^\circ,\,90^\circ]
\end{equation}
where $\mathbf{u}_k=[u_{x,k}\;\;u_{y,k}\;\;u_{z,k}]^\top$.
We note that the two trajectories look quite similar in the position space, but the control profile comparison offers better clarity into the differences between the two trajectories. Specifically, we note that the MASTO nominal trajectory introduces some margin in the control profile for most thrusting arcs, which allows for corrective maneuvers while respecting control chance constraints. While introducing these margins, the MASTO nominal control profile also tends to have longer thrusting arcs, to compensate for the new margins.

\begin{figure}[htbp!]
    \centering
    \begin{subfigure}[b]{0.48\textwidth}
        \centering
        \includegraphics[height=0.67\linewidth]{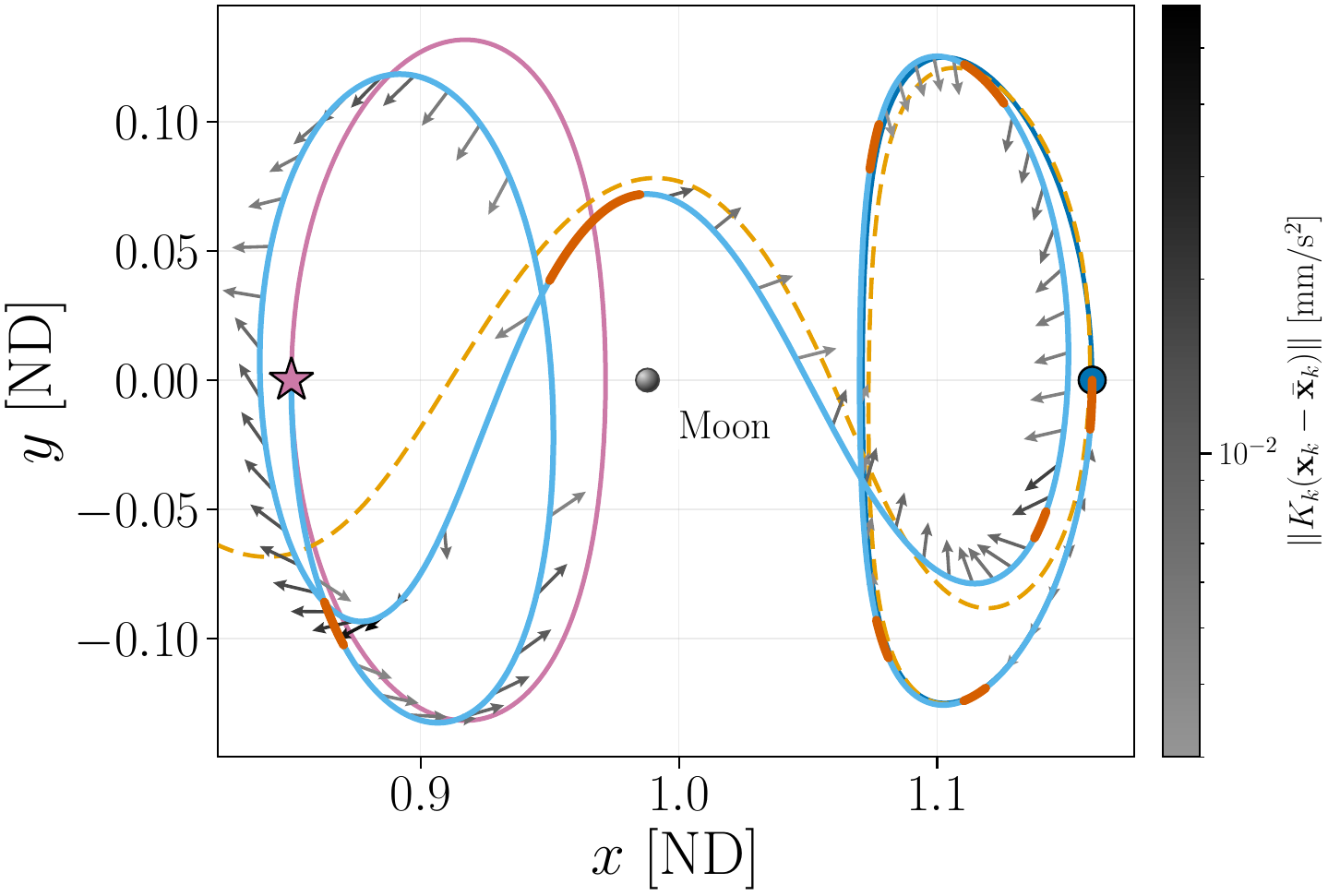}
        \caption{$xy$-projection.}
        \label{fig:hh_smp_xy}
    \end{subfigure}
    \hfill
    \begin{subfigure}[b]{0.48\textwidth}
        \centering
        \includegraphics[height=0.67\linewidth]{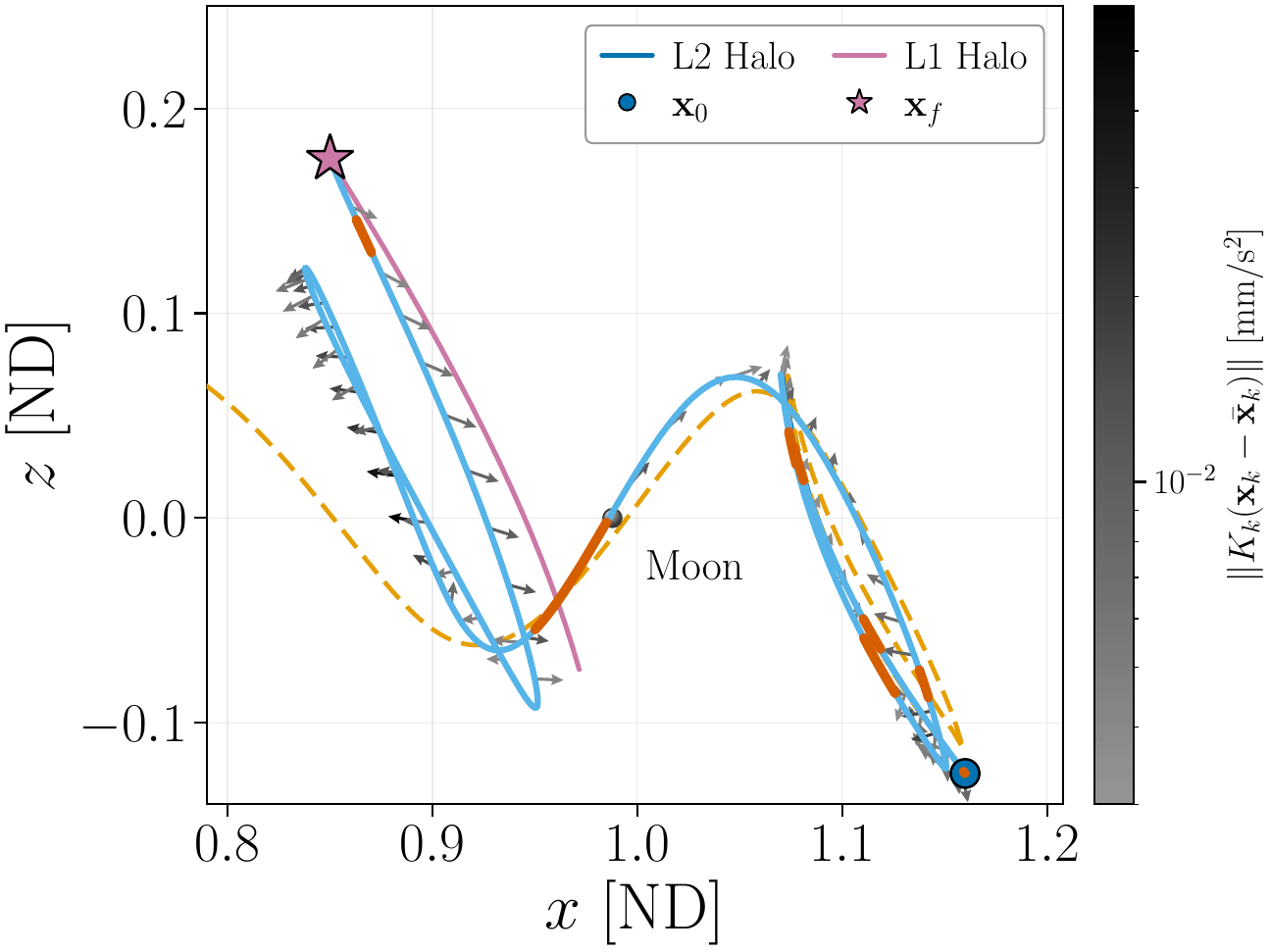}
        \caption{$xz$-projection.}
        \label{fig:hh_smp_xz}
    \end{subfigure}

    \vspace{0.6em}

    \begin{subfigure}[b]{0.48\textwidth}
        \centering
        \includegraphics[height=0.84\linewidth]{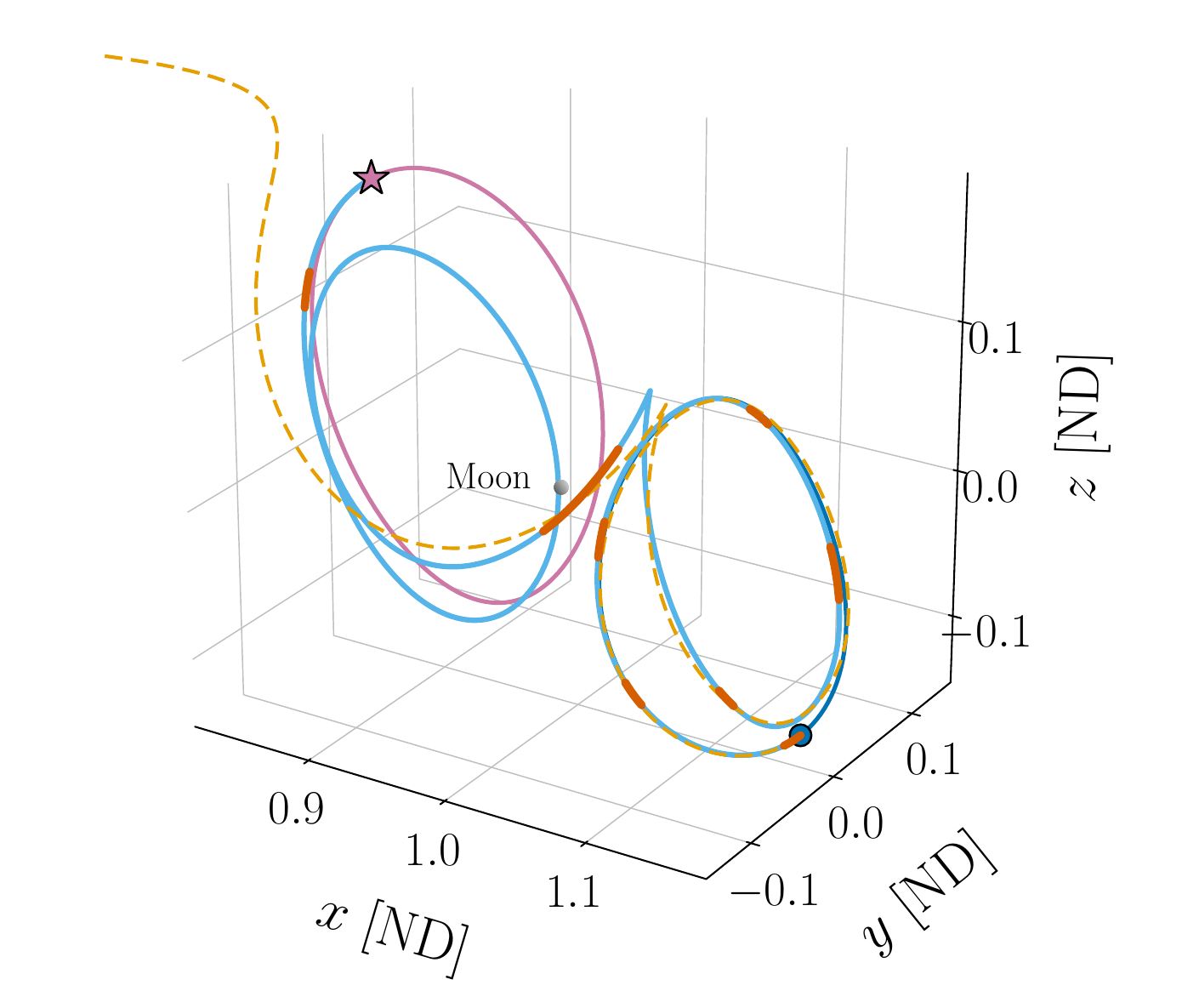}
        \caption{Three-dimensional view.}
        \label{fig:hh_smp_3d}
    \end{subfigure}
    \hfill
    \begin{subfigure}[b]{0.48\textwidth}
        \centering
        \includegraphics[height=0.84\linewidth]{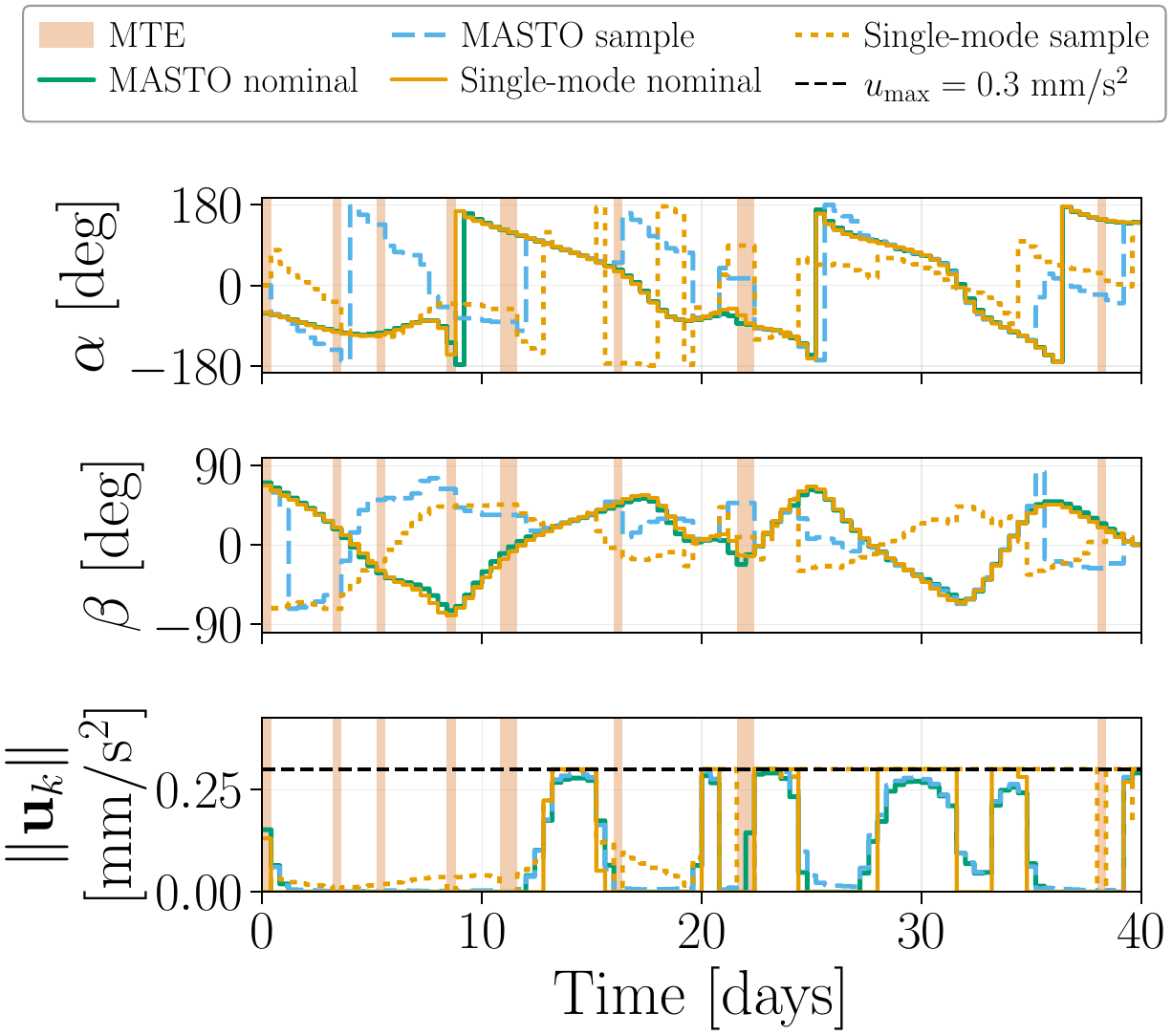}
        \caption{Control angles and magnitude, outages shaded.}
        \label{fig:hh_smp_ctrl}
    \end{subfigure}
    \caption{Monte Carlo sample with eight MTEs: MASTO (blue, outage arcs red) against the single-mode controller under the same outages (dashed orange).}
    \label{fig:hh_sample}
\end{figure}

We now study how the MASTO nominal trajectory performs under several MTEs, and how it compares against a non-MTE aware stochastic optimal control method (single-mode covariance steering). We perform nonlinear Monte Carlo simulations with 250 samples, and generate the MTEs according to the given distribution in \Cref{tab:hh_params}. We first isolate one sample: it faces eight MTEs, which correspond to a total of four days spent in failure mode. Two of these MTEs coincide with the MASTO nominal thrusting arc, leading to around twenty hours of lost thrust time. We perform this simulation under both the MASTO policy and the single-mode covariance steering policy, and the corresponding comparison is shown in \Cref{fig:hh_sample}. \Cref{fig:hh_smp_xy} and \Cref{fig:hh_smp_xz} show the $xy$ and $xz$ projections, respectively, where the blue and red show the MASTO trajectory denoting the nominal and failure mode respectively. The arrows indicate the feedback action with the colorbar corresponding to the magnitude of the corrective maneuver. The orange dashed line indicates the same simulation performed with the single-mode policy, which diverges under this sequence of MTEs. \Cref{fig:hh_smp_3d} shows the trajectory for the two cases, highlighting how the single-mode controller fails to bring the spacecraft to the desired final orbit's vicinity in the given time under the control norm constraints. We note that the feedback tries to recover the lost thrusting time in the single-mode case, as shown in the control history plot in \Cref{fig:hh_smp_ctrl}, but still diverges. We also note that in the nominal case, only the first MTE affects the single-mode thrusting arc, while the MASTO nominal thrusting arcs are affected by two MTEs. However, since the MASTO nominal trajectory is designed anticipating future MTEs, it is able to recover the spacecraft despite multiple MTEs, unlike the single-mode case where a single MTE proved to be catastrophic.

\begin{figure}[htbp!]
    \centering
    \begin{subfigure}[b]{0.48\textwidth}
        \centering
        \includegraphics[height=0.78\linewidth]{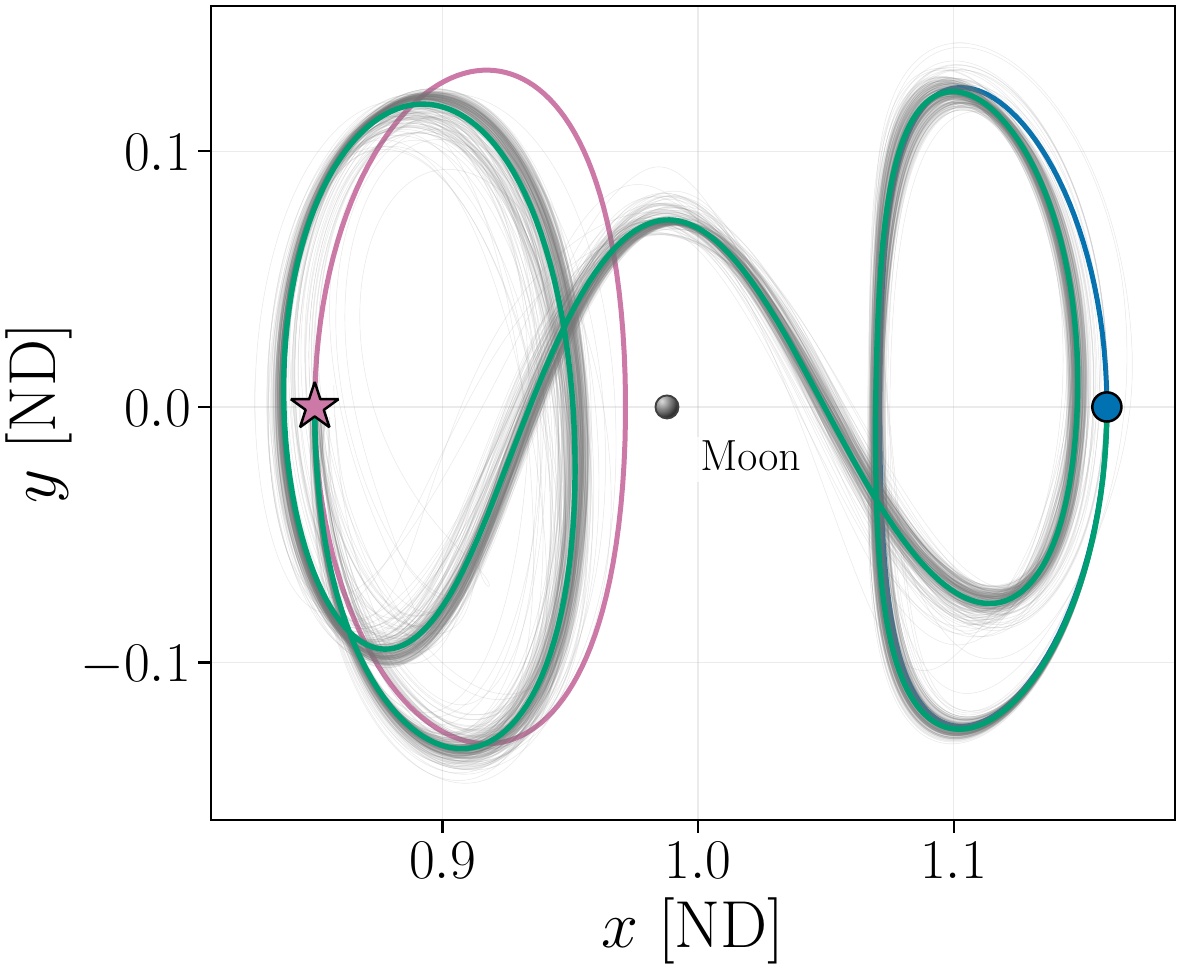}
        \caption{$xy$-projection.}
        \label{fig:hh_mc_xy}
    \end{subfigure}
    \hfill
    \begin{subfigure}[b]{0.48\textwidth}
        \centering
        \includegraphics[height=0.78\linewidth]{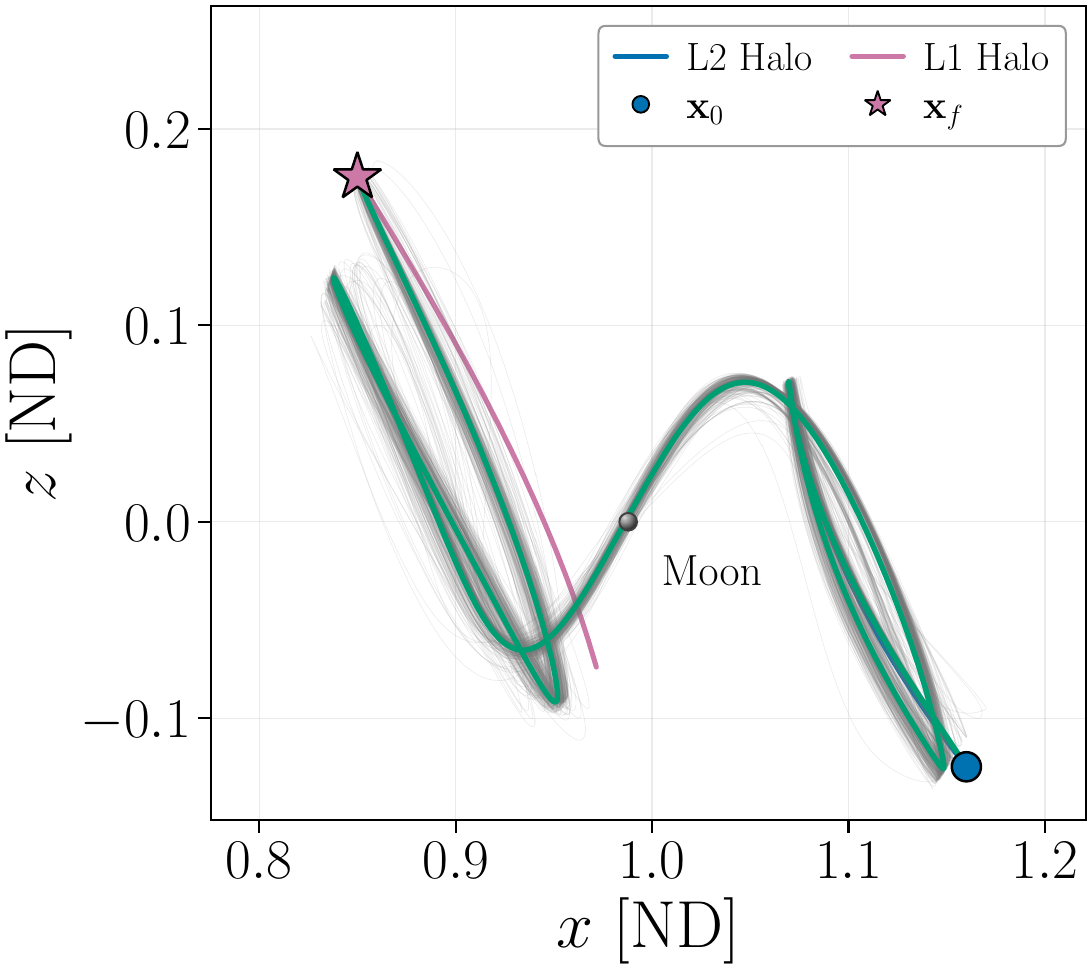}
        \caption{$xz$-projection.}
        \label{fig:hh_mc_xz}
    \end{subfigure}

    \vspace{0.6em}

    \begin{subfigure}[b]{0.48\textwidth}
        \centering
        \includegraphics[height=0.84\linewidth]{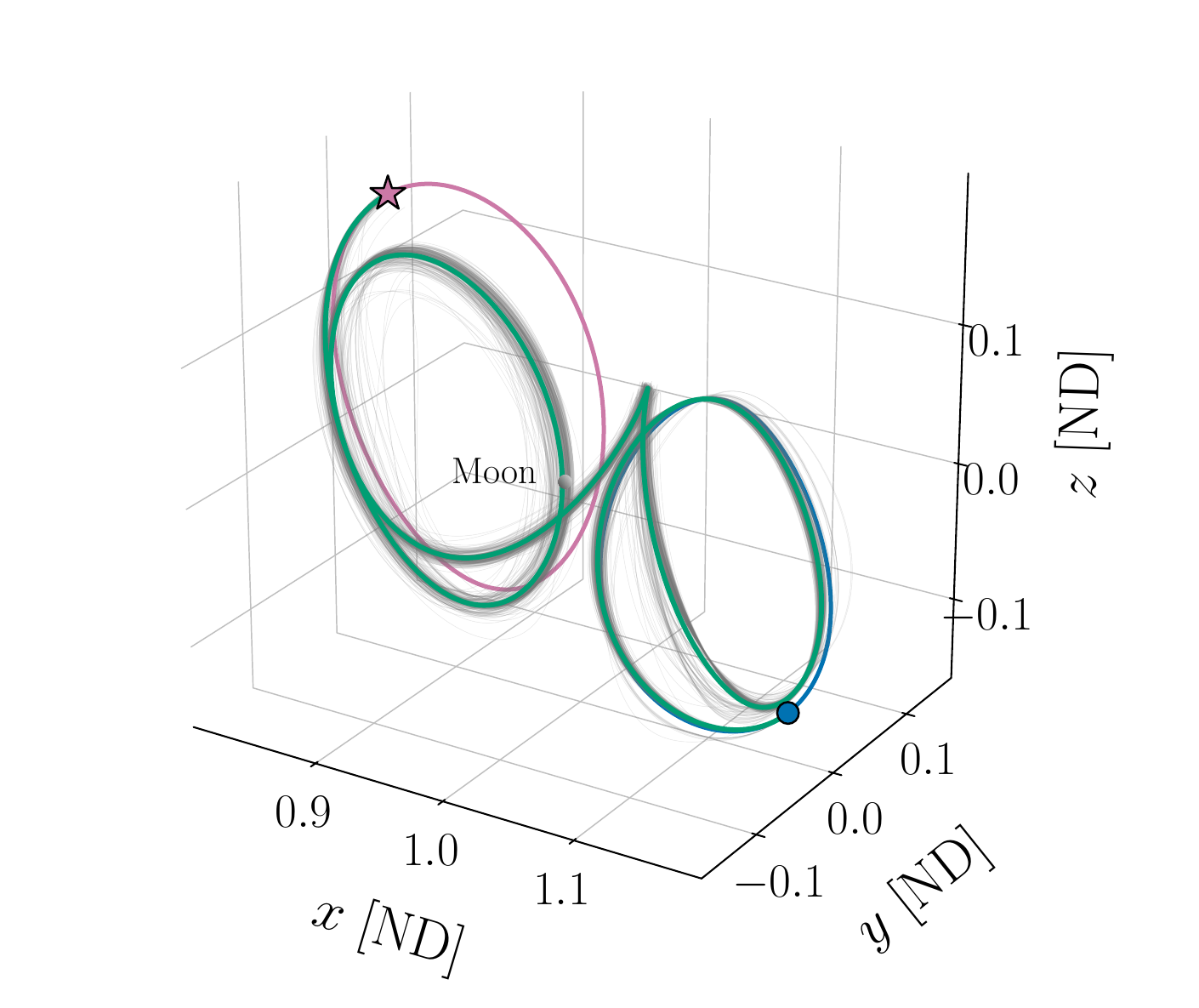}
        \caption{Three-dimensional view.}
        \label{fig:hh_mc_3d}
    \end{subfigure}
    \hfill
    \begin{subfigure}[b]{0.48\textwidth}
        \centering
        \includegraphics[height=0.84\linewidth]{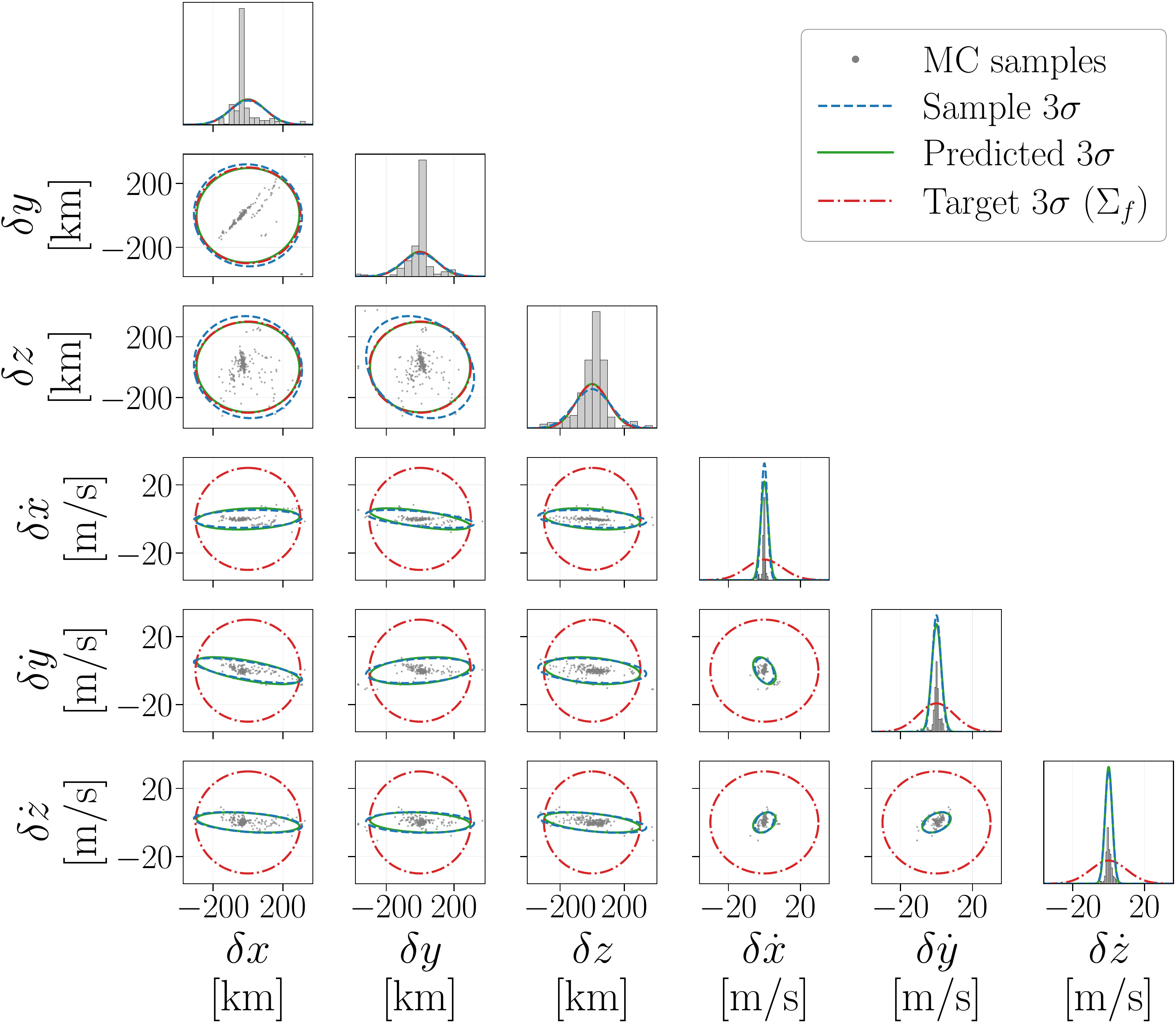}
        \caption{Terminal state dispersion.}
        \label{fig:hh_mc_disp}
    \end{subfigure}
    \caption{Monte Carlo results for the L2 Halo to L1 Halo transfer; deviations amplified $5\times$ for trajectory plots.}
    \label{fig:hh_mc}
\end{figure}

\begin{figure}[htbp!]
    \centering
    \begin{subfigure}[b]{0.48\textwidth}
        \centering
        \includegraphics[height=0.74\linewidth]{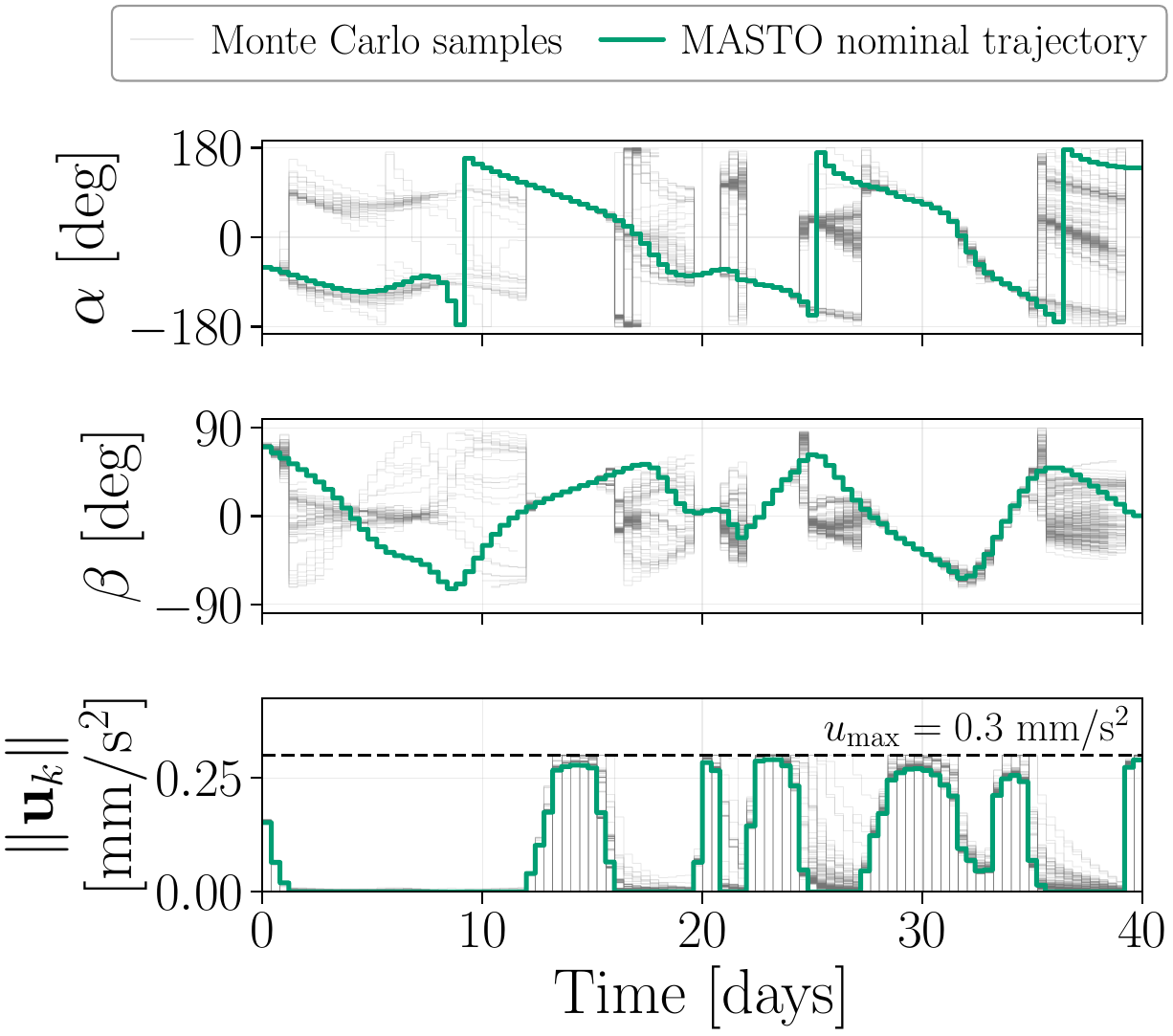}
        \caption{Monte Carlo control angles and magnitude.}
        \label{fig:hh_mc_ctrl}
    \end{subfigure}
    \hfill
    \begin{subfigure}[b]{0.48\textwidth}
        \centering
        \includegraphics[height=0.74\linewidth]{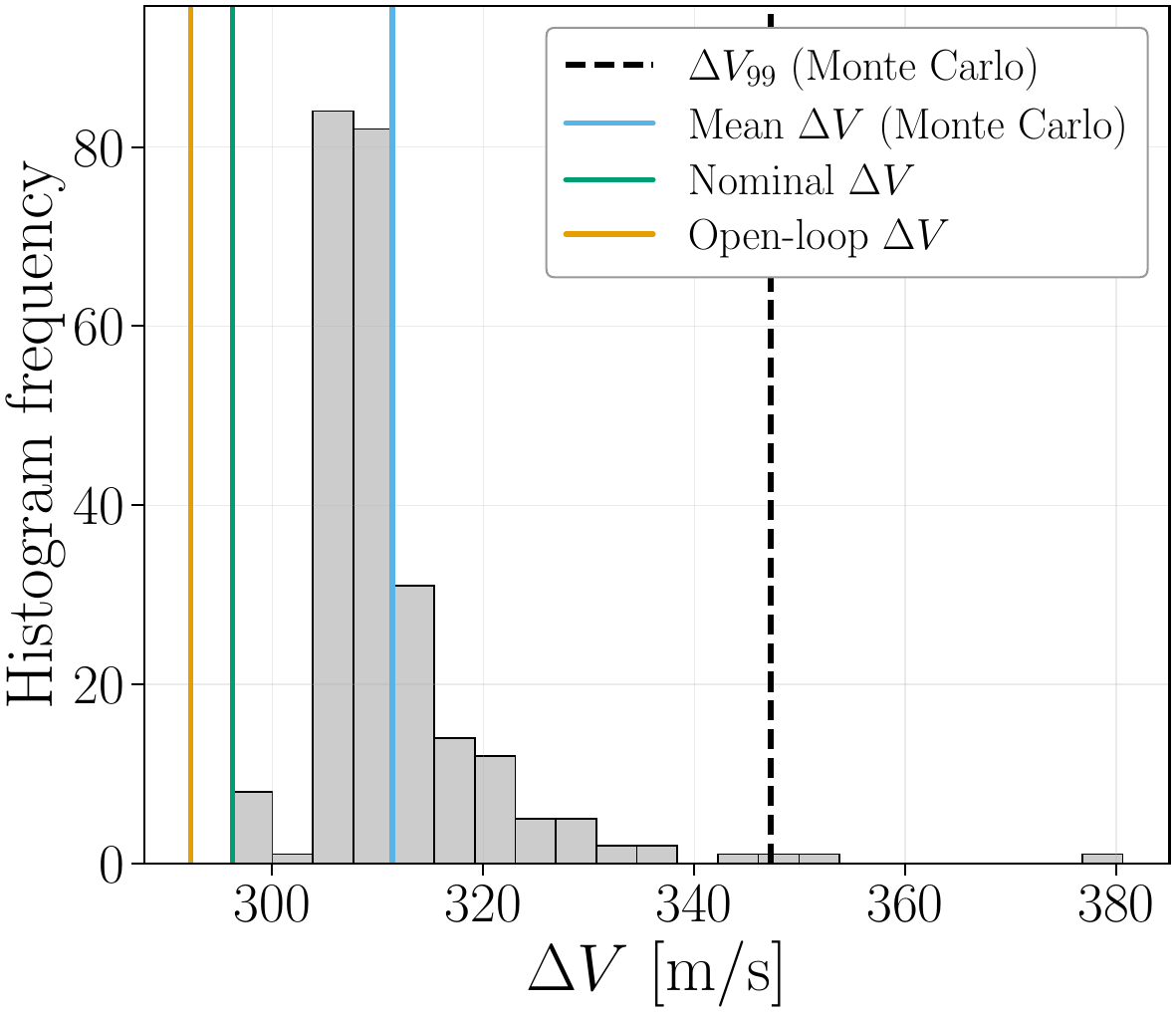}
        \caption{$\Delta v_{99}$ histogram.}
        \label{fig:hh_dv_hist}
    \end{subfigure}
    \caption{Monte Carlo control histories and $\Delta v$ distribution for the L2 Halo to L1 Halo transfer.}
    \label{fig:hh_ctrl_dv}
\end{figure}

Further, all Monte Carlo sample trajectories under the MASTO policy are shown in \Cref{fig:hh_mc}, with \Cref{fig:hh_mc_xy}, \Cref{fig:hh_mc_xz} and \Cref{fig:hh_mc_3d} showing the $xy$, $xz$ projections and the three-dimensional trajectory respectively. We amplify the deviations by five times for visual clarity in the trajectory plots. \Cref{fig:hh_mc_disp} shows the terminal state dispersion, which shows good agreement between the sample and predicted terminal covariance, while satisfying the terminal covariance bound. In position space, $92.0\%$ of samples fall within $3\sigma_{r_f}$, meeting the prescribed bounds well. The terminal distribution exhibits some non-Gaussian behavior, and can be explained by both the nonlinear dynamics and the MTE-uncertainty model which introduce non-Gaussianity. Addressing this issue would require the use of higher-order uncertainty quantification techniques, such as the conjugate unscented transformation\cite{qiNonGaussianDistributionSteering2025} in the MASTO framework, which is left for future work. The corresponding control history plot for the Monte Carlo simulations is shown in \Cref{fig:hh_mc_ctrl}, with the $\Delta v_{99}$ histogram shown in \Cref{fig:hh_dv_hist}. We note that the nominal $\Delta v=296.3$ m/s from MASTO framework, which represents the open-loop component of the control policy, is only slightly higher than the deterministic reference's $292.3$ m/s $\Delta v$ consumption, while introducing resilience to MTEs. This implies that the cost of adding resilience to MTEs is minimal in the nominal trajectory for this case. In the Monte Carlo simulations with arbitrary number of MTEs, the MASTO framework yields a $\Delta v_{99}=347.2$ m/s which is $18.7\%$ higher than the open-loop reference. This is similar to our previous results\cite{shrivastavaMASTOMISSEDTHRUSTAWARE2026} which showed similar percentage increase in $\Delta v_{99}$ for other transfers, with marginal increases in nominal $\Delta v$ compared to the open-loop reference.
\begin{figure}[htbp!]
    \centering
    \includegraphics[width=0.6\textwidth]{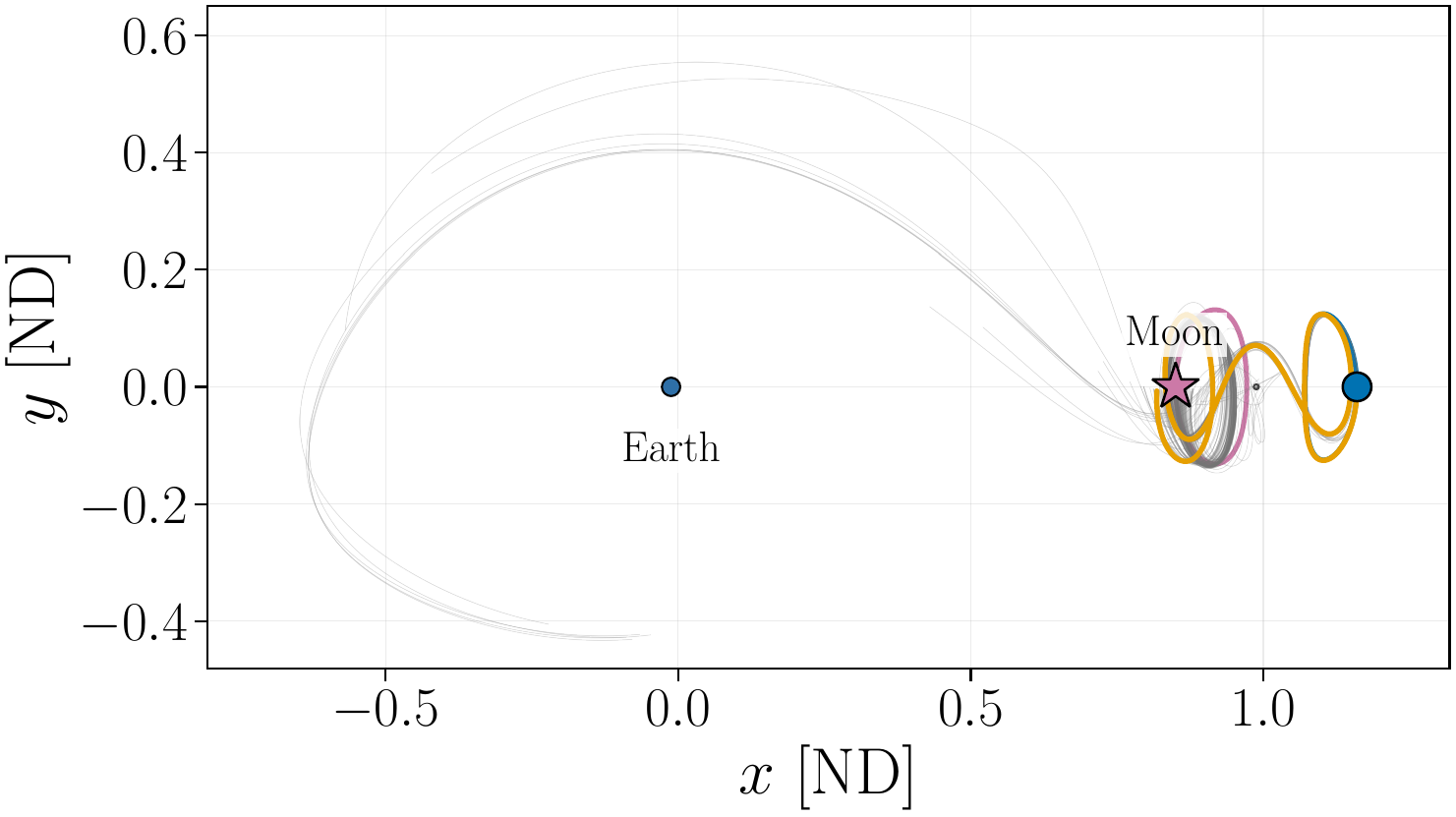}
    \caption{Single-mode Monte Carlo simulations ($xy$-projection) with MTEs for the L2 Halo to L1 Halo transfer.}
    \label{fig:hh_mc_singlemode_xy}
\end{figure}

We also perform Monte Carlo simulations with the single-mode covariance steering policy, with the same MTE samples simulated in \Cref{fig:hh_mc} for the MASTO Monte Carlo simulations. \Cref{fig:hh_mc_singlemode_xy} shows the $xy$-projection of the Monte Carlo simulations, with $133$ of the total $250$ samples diverging (the final position exceeds $20\sigma_{r_f}$). In comparison, the MASTO trajectories had zero divergences, and also satisfied the lunar surface clearance constraint in all samples.

\subsection{L2 NRHO to L1 Halo Transfer in Earth-Moon CR3BP}

\begin{figure}[htbp!]
    \centering
    \begin{subfigure}[b]{0.48\textwidth}
        \centering
        \includegraphics[height=0.80\linewidth]{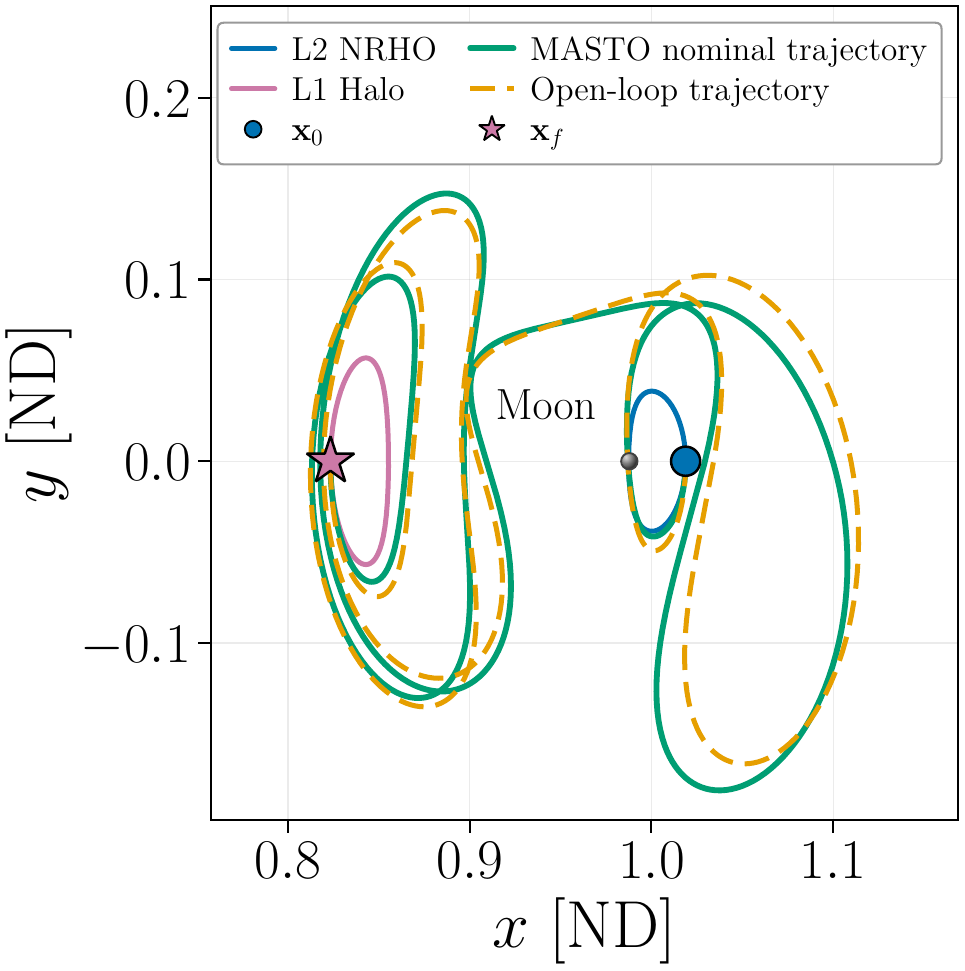}
        \caption{$xy$-projection.}
        \label{fig:nh_ovm_xy}
    \end{subfigure}
    \hfill
    \begin{subfigure}[b]{0.48\textwidth}
        \centering
        \includegraphics[height=0.80\linewidth]{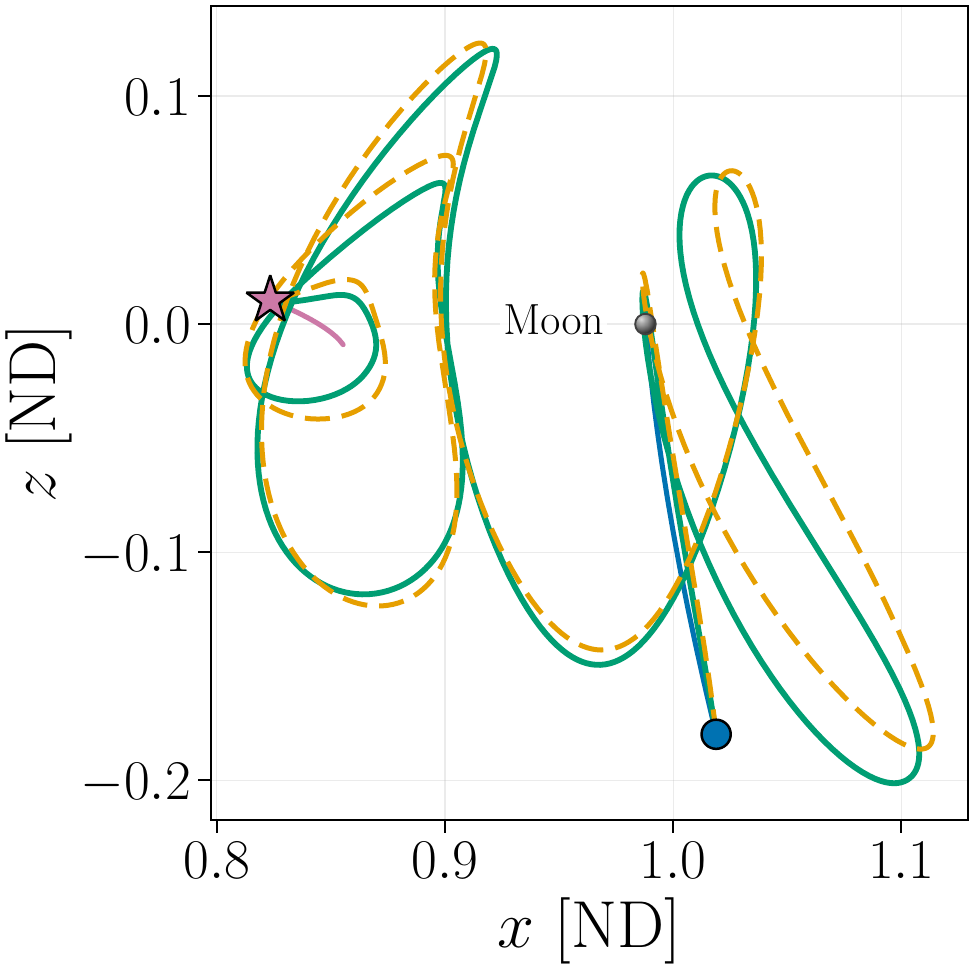}
        \caption{$xz$-projection.}
        \label{fig:nh_ovm_xz}
    \end{subfigure}

    \vspace{0.6em}

    \begin{subfigure}[b]{0.48\textwidth}
        \centering
        \includegraphics[height=0.84\linewidth]{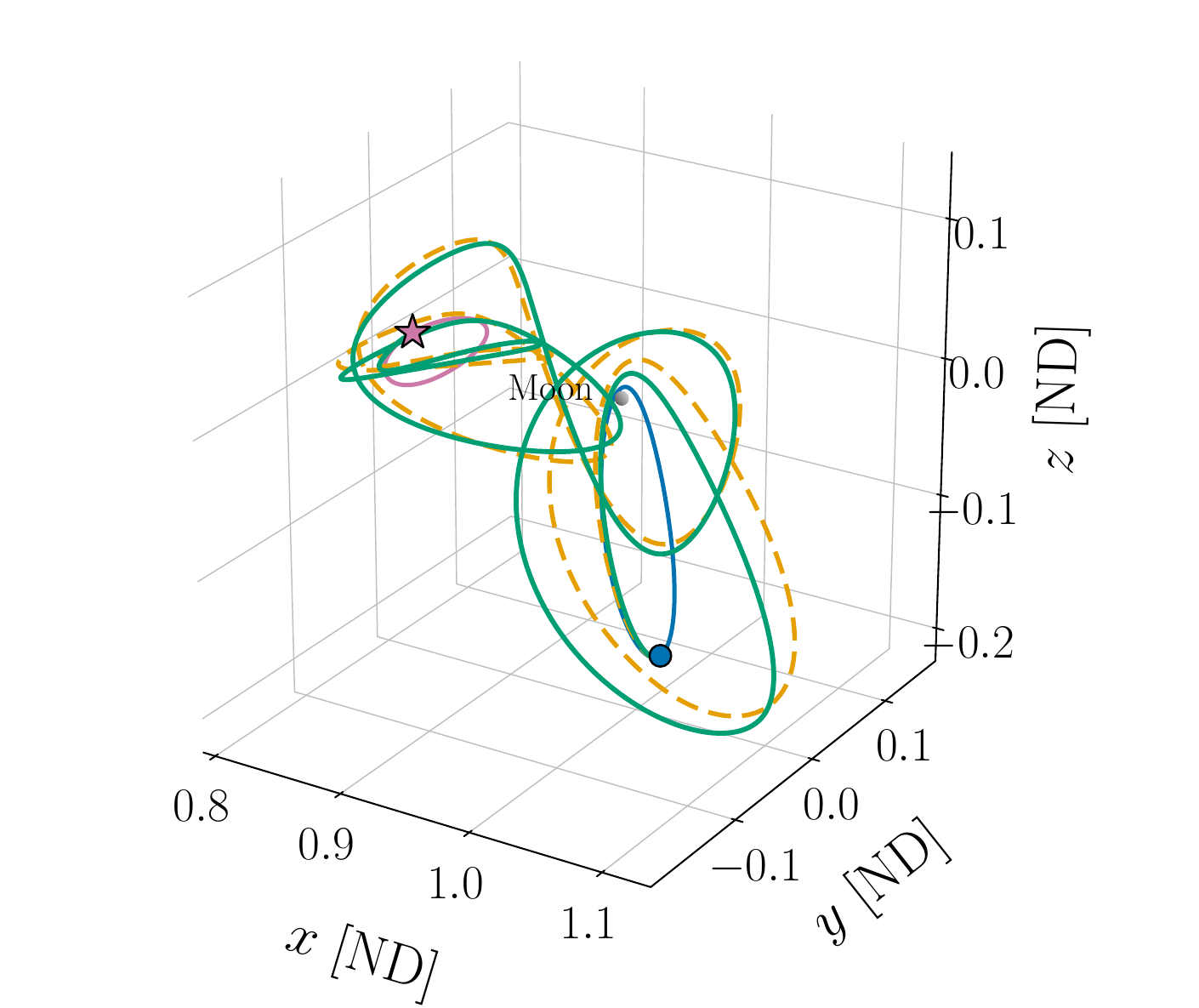}
        \caption{Three-dimensional view.}
        \label{fig:nh_ovm_3d}
    \end{subfigure}
    \hfill
    \begin{subfigure}[b]{0.48\textwidth}
        \centering
        \includegraphics[height=0.84\linewidth]{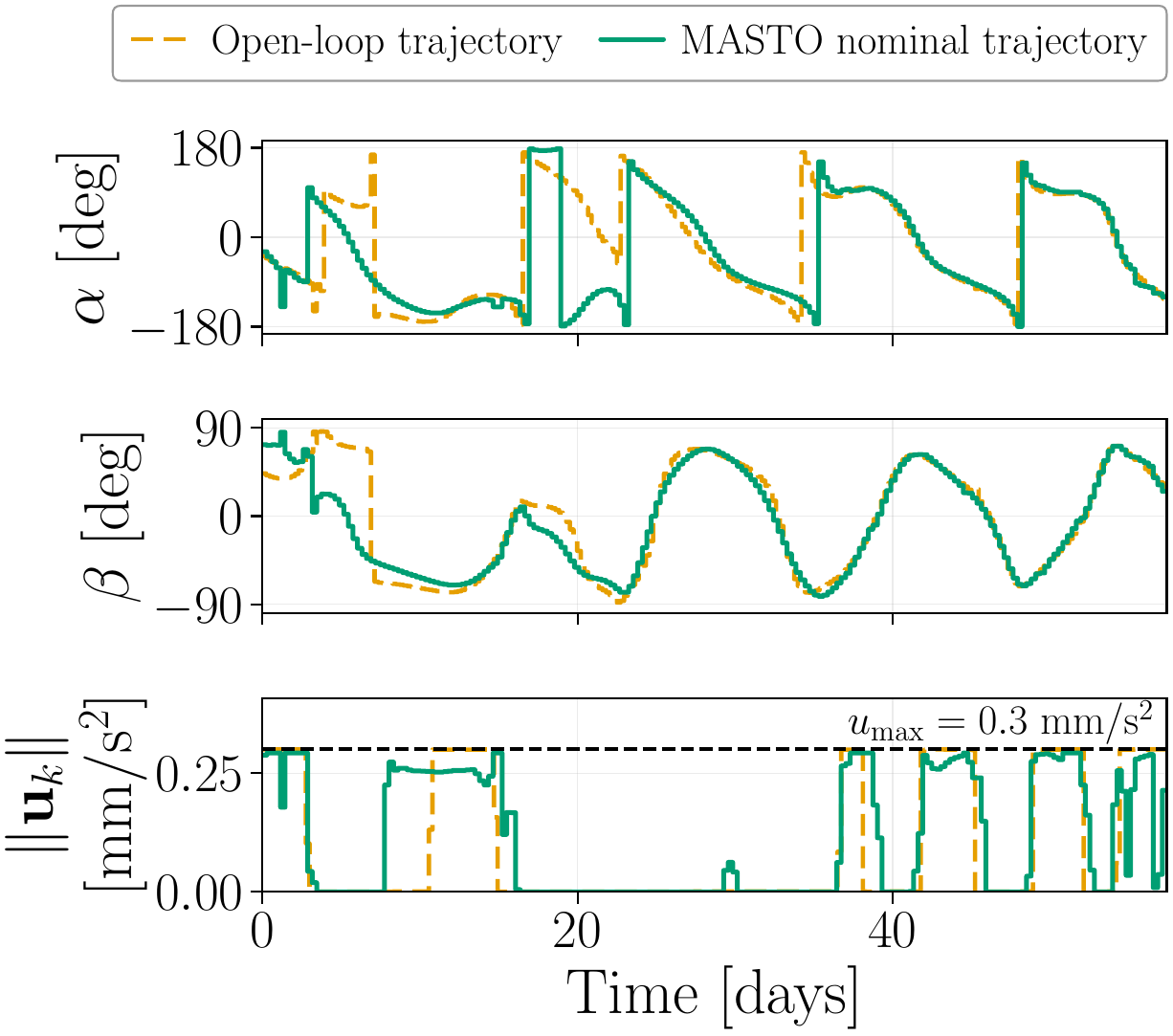}
        \caption{Control angles and magnitude.}
        \label{fig:nh_ovm_ctrl}
    \end{subfigure}
    \caption{Nominal MASTO against the open-loop minimum-fuel solution, L2 NRHO to L1 Halo transfer.}
    \label{fig:nh_ovm}
\end{figure}

We next consider a more complex multi-revolution transfer from an L2 NRHO to an L1 Halo orbit in the Earth-Moon CR3BP. \Cref{tab:nh_ic} lists the initial conditions for the departure and arrival orbits, and \Cref{tab:nh_params} lists the simulation parameters. We enforce the same lunar surface clearance chance constraint with $100$ km altitude. We note that the reference deterministic transfer in this case (shown in dashed lines in \Cref{fig:nh_ovm}) has a long coasting arc in the middle which relies on proper execution of the first two thrusting arcs. Under the current MTE model listed in \Cref{tab:nh_params}, the spacecraft is likely to face at least one MTE in the first half of the transfer trajectory with $90\%$ probability. This can be catastrophic, especially considering the close lunar flyby. Further, the second half of the transfer has several thrusting arcs that are critical to ensure the spacecraft reaches the L1 Halo orbit which has very different geometry and inclination compared to the departure orbit. 
\begin{table}[h!]
\centering
\caption{Initial conditions for L2 NRHO to L1 Halo transfer ($y=\dot{x}=\dot{z}=0$)}
\label{tab:nh_ic}
\begin{tabular}{lrrrr}
\hline\hline
\multicolumn{1}{c}{Orbit} & \multicolumn{1}{c}{$x$ (ND)} & \multicolumn{1}{c}{$z$ (ND)} & \multicolumn{1}{c}{$\dot{y}$ (ND)} & \multicolumn{1}{c}{Period (days)} \\ \hline
L2 NRHO  & $1.018826173554963$ & $-0.179797844569828$ & $-0.096189089845127$ & $\approx6.4$ \\
L1 Halo  & $0.823383959653906$ & $0.010388134109586$  & $0.128105259453086$  & $\approx11.9$ \\ \hline\hline
\end{tabular}
\end{table}
\begin{table}[htbp!]
\centering
\caption{Parameters for the L2 NRHO to L1 Halo transfer}
\label{tab:nh_params}
\begin{tabular}{lrl|lrl}
\hline\hline
Parameter & Value & Unit & Parameter & Value & Unit \\ \hline
Time-of-flight           & 57.4 & days     & Discretization nodes & 201 & -- \\
$u_{\max}$   & 0.3  & mm/s$^2$ & $\text{TTF}_{90}$    & 30  & days \\
$\sigma_{r_0}$           & 5    & km       & $\text{TTR}_{95}$    & 0.5 & days \\
$\sigma_{v_0}$           & 0.1  & m/s      & \multirow{4}{*}{$\mathcal{P}$} & \multirow{4}{*}{$\begin{bmatrix}0.97820 & 0.02180\\0.82090 & 0.17910\end{bmatrix}$} & \multirow{4}{*}{--} \\
$\sigma_{r_f}$           & 100  & km       & & & \\
$\sigma_{v_f}$           & 10   & m/s      & & & \\
$\sigma_{r,\max}$        & 500  & km       & & & \\
$\sigma_{v,\max}$        & 15   & m/s      & $\bm{\rho}_0$ & $[0.974,\ 0.026]^\top$ & -- \\ \hline\hline
\end{tabular}
\end{table}

We run \Cref{alg:masto} to get the MTE-resilient MASTO trajectory for this transfer. It converges in $51$ SCP iterations, taking $38$ minutes in total ($\approx45$ s per iteration on average) under the same computational setup. The complexity of this transfer with its close lunar flyby, its sensitivity to MTEs and the larger problem size contribute to the longer convergence time.

We first compare how the MASTO nominal trajectory differs from the deterministic open-loop reference trajectory. \Cref{fig:nh_ovm_xy}, \Cref{fig:nh_ovm_xz},  and \Cref{fig:nh_ovm_3d} compare the $xy$, $xz$ projections, and the three-dimensional trajectories respectively, while \Cref{fig:nh_ovm_ctrl} compares the respective control profiles. Unlike the previous example, the MASTO nominal trajectory in this case has clearly different trajectory geometry as seen in \Cref{fig:nh_ovm_xy} and \Cref{fig:nh_ovm_xz}. \Cref{fig:nh_ovm_3d} shows that the MASTO nominal trajectory takes a longer turn coming out of the departure NRHO, and shows significantly different approach geometry. The control profile also differs substantially from the open-loop reference, as the MASTO nominal trajectory introduces some margins while spreading the thrusting arcs over a longer time duration. 

\begin{figure}[htbp!]
    \centering
    \begin{subfigure}[b]{0.48\textwidth}
        \centering
        \includegraphics[height=0.80\linewidth]{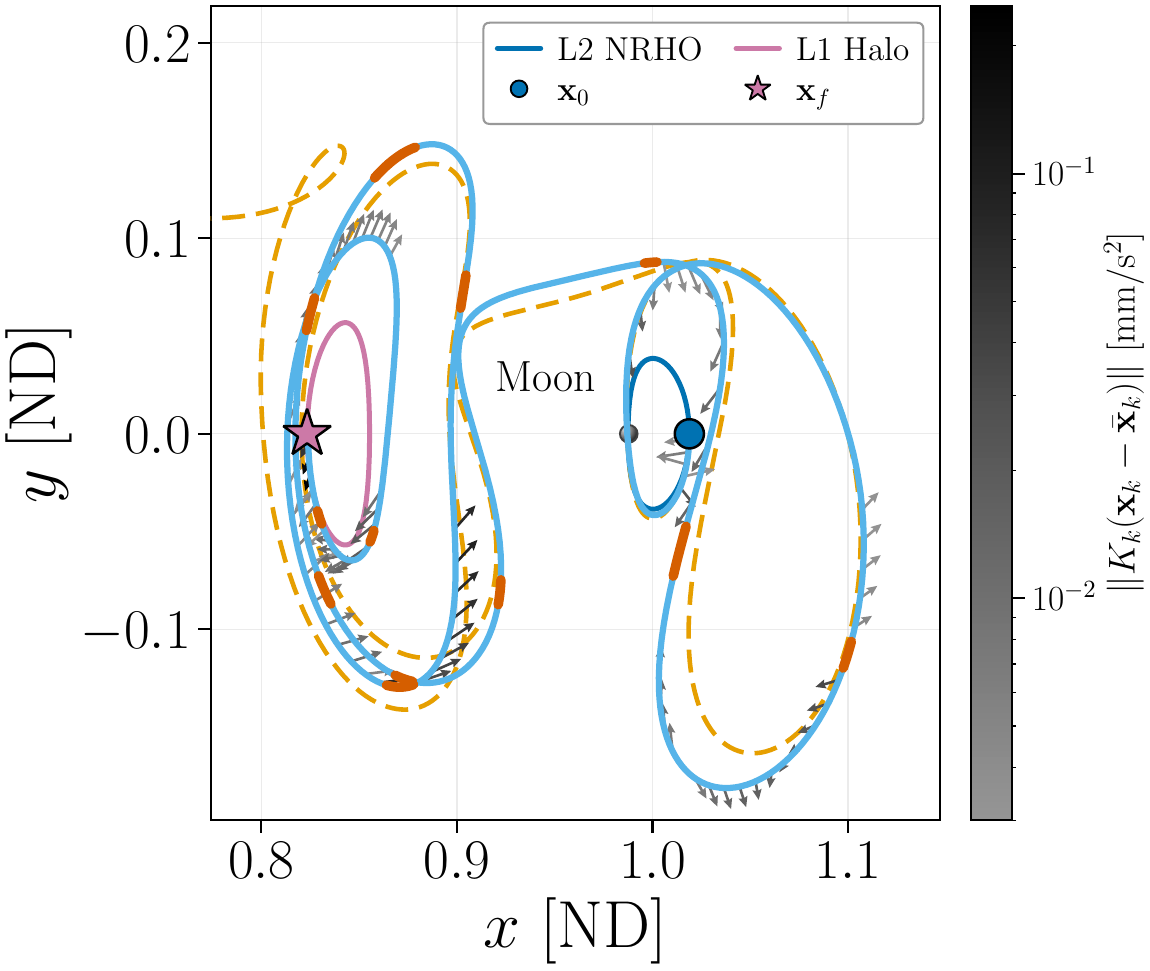}
        \caption{$xy$-projection.}
        \label{fig:nh_smp_xy}
    \end{subfigure}
    \hfill
    \begin{subfigure}[b]{0.48\textwidth}
        \centering
        \includegraphics[height=0.80\linewidth]{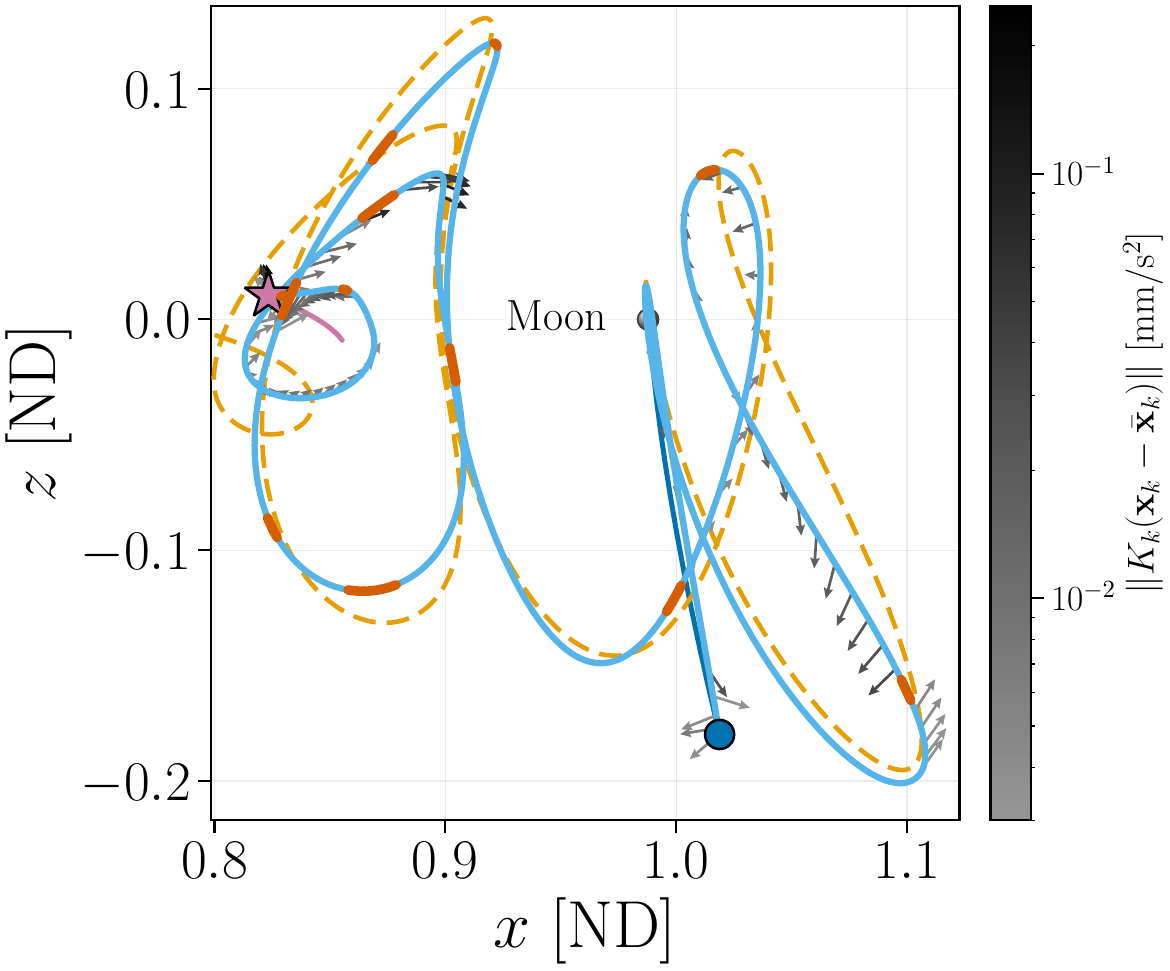}
        \caption{$xz$-projection.}
        \label{fig:nh_smp_xz}
    \end{subfigure}

    \vspace{0.6em}

    \begin{subfigure}[b]{0.48\textwidth}
        \centering
        \includegraphics[height=0.84\linewidth]{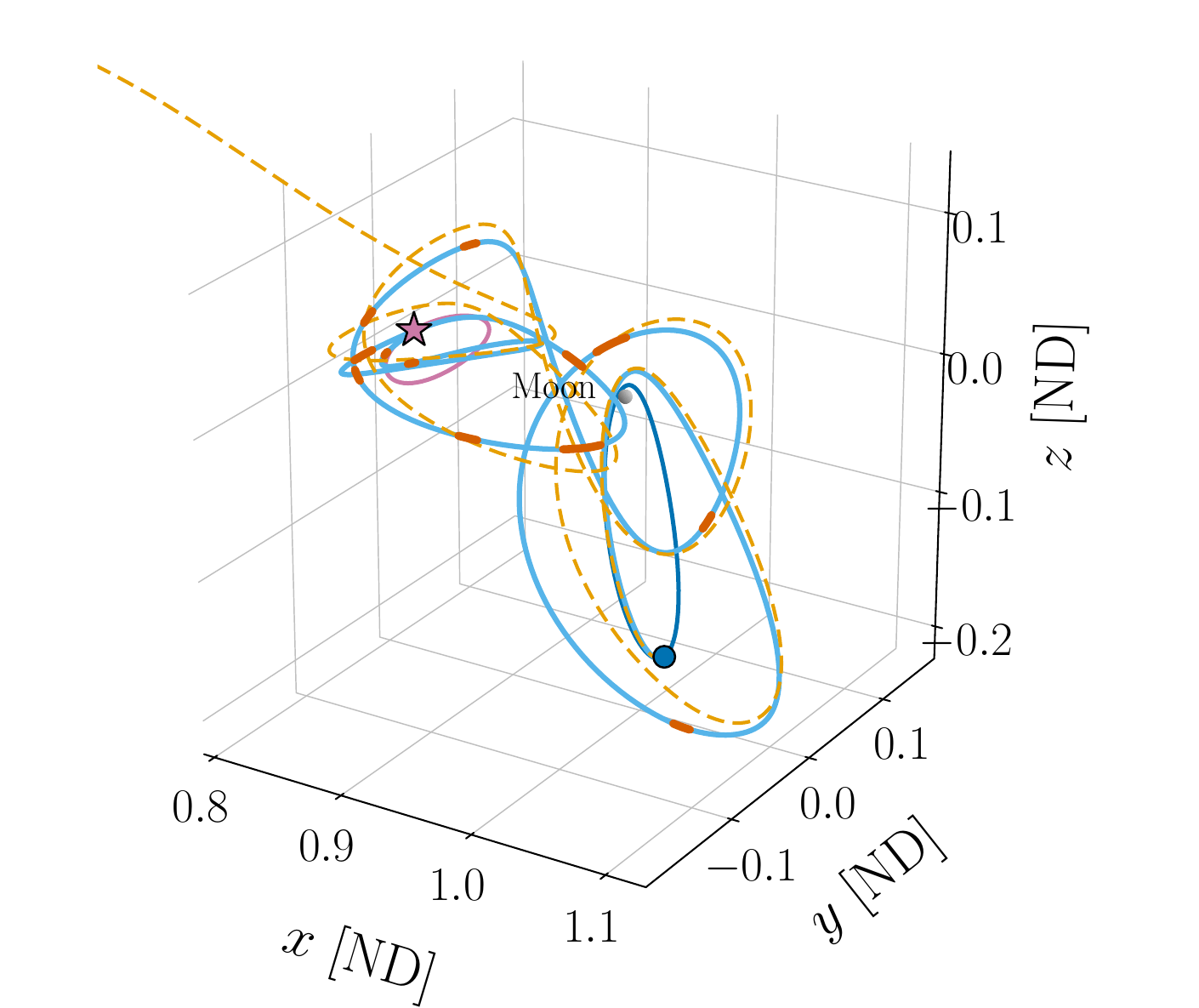}
        \caption{Three-dimensional view.}
        \label{fig:nh_smp_3d}
    \end{subfigure}
    \hfill
    \begin{subfigure}[b]{0.48\textwidth}
        \centering
        \includegraphics[height=0.84\linewidth]{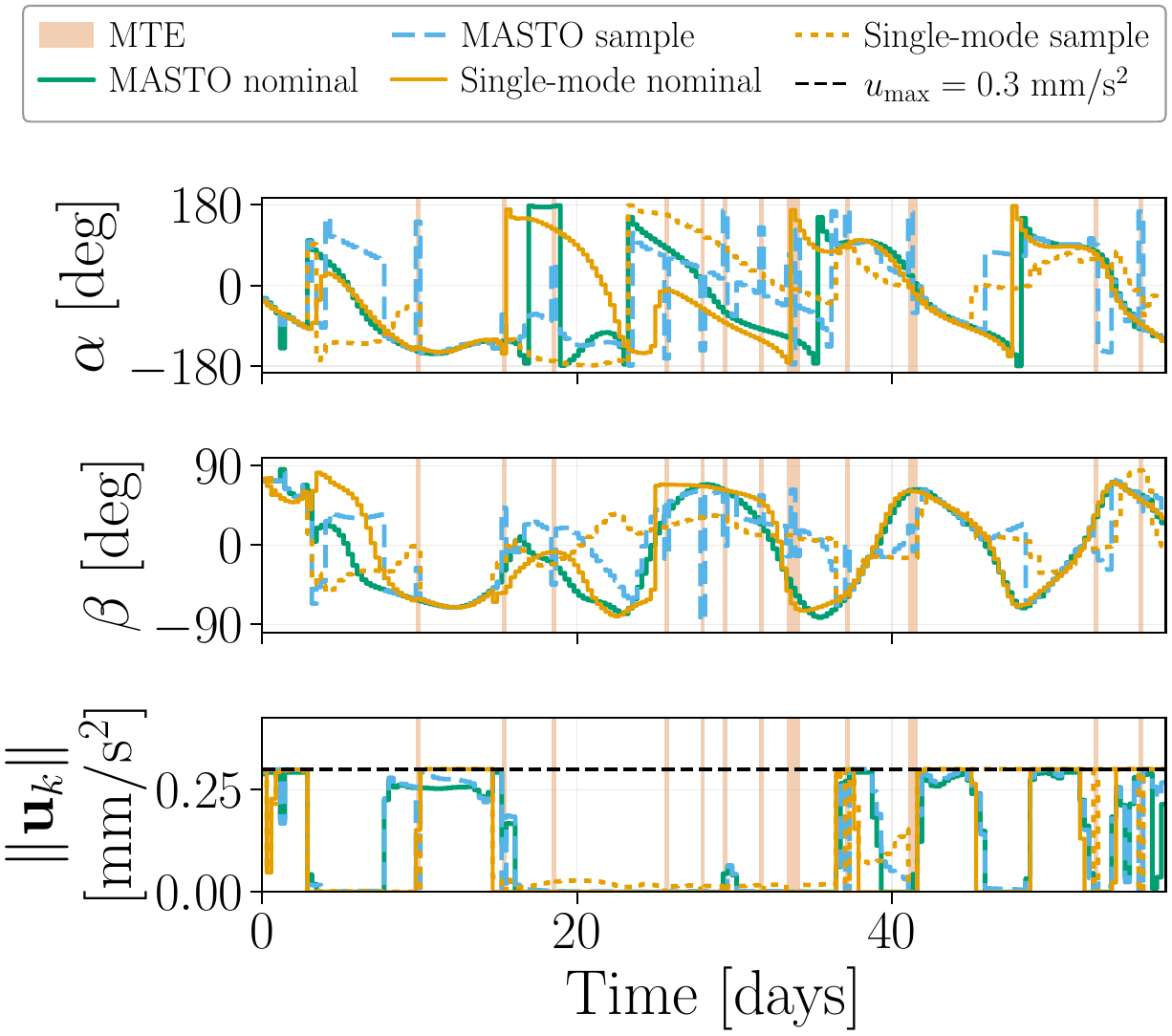}
        \caption{Control angles and magnitude, outages shaded.}
        \label{fig:nh_smp_ctrl}
    \end{subfigure}
    \caption{Monte Carlo sample with twelve MTEs: MASTO (blue, outage arcs red) against the single-mode controller under the same outages (dashed orange).}
    \label{fig:nh_sample}
\end{figure}

We now study how the MASTO nominal trajectory performs under several MTEs for this transfer, and how it compares against a non-MTE aware stochastic optimal control method (single-mode covariance steering). We perform nonlinear Monte Carlo simulations, with the MTEs generated according to the given distribution in \Cref{tab:nh_params}. As in the previous example, we first isolate one sample to compare the two trajectories in detail. This sample faces twelve MTEs, six of which coincide with the MASTO nominal thrusting arcs, while three coincide with the single-mode nominal thrusting arcs. We perform this simulation under both the MASTO policy and the non-MTE aware single-mode covariance steering policy, and the corresponding comparison is shown in \Cref{fig:nh_sample}. \Cref{fig:nh_smp_xy} and \Cref{fig:nh_smp_xz} show the $xy$ and $xz$ projections, respectively, {where the blue and red show the MASTO trajectory denoting the nominal and failure mode respectively. The shaded arrows in \Cref{fig:nh_smp_xy,fig:nh_smp_xz} indicate the corrective feedback action taken by the MASTO policy, with darker shades indicating higher feedback magnitude. The orange dashed line indicates the non-MTE aware single-mode trajectory under the same sequence of MTEs, which fails to recover the spacecraft under this sequence of MTEs, while the MASTO policy recovers.} \Cref{fig:nh_smp_3d} {shows} the three-dimensional trajectory, and \Cref{fig:nh_smp_ctrl} the control angles and magnitude with MTEs shaded, for both the MASTO and single-mode policies. This sample demonstrates how the MASTO policy is able to recover the spacecraft despite multiple MTEs, through its feedback action, while the non-MTE aware method fails to recover. 

\begin{figure}[htbp!]
    \centering
    \begin{subfigure}[b]{0.48\textwidth}
        \centering
        \includegraphics[height=0.80\linewidth]{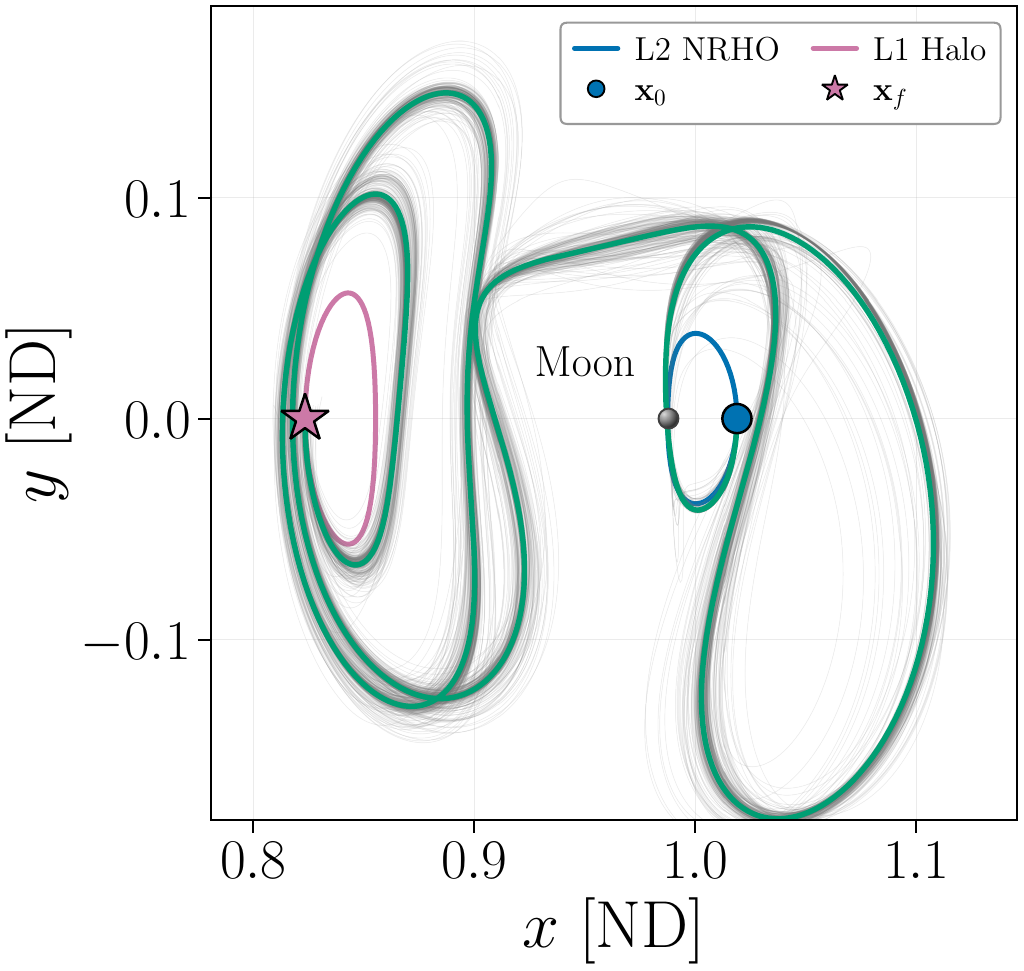}
        \caption{$xy$-projection.}
        \label{fig:nh_mc_xy}
    \end{subfigure}
    \hfill
    \begin{subfigure}[b]{0.48\textwidth}
        \centering
        \includegraphics[height=0.80\linewidth]{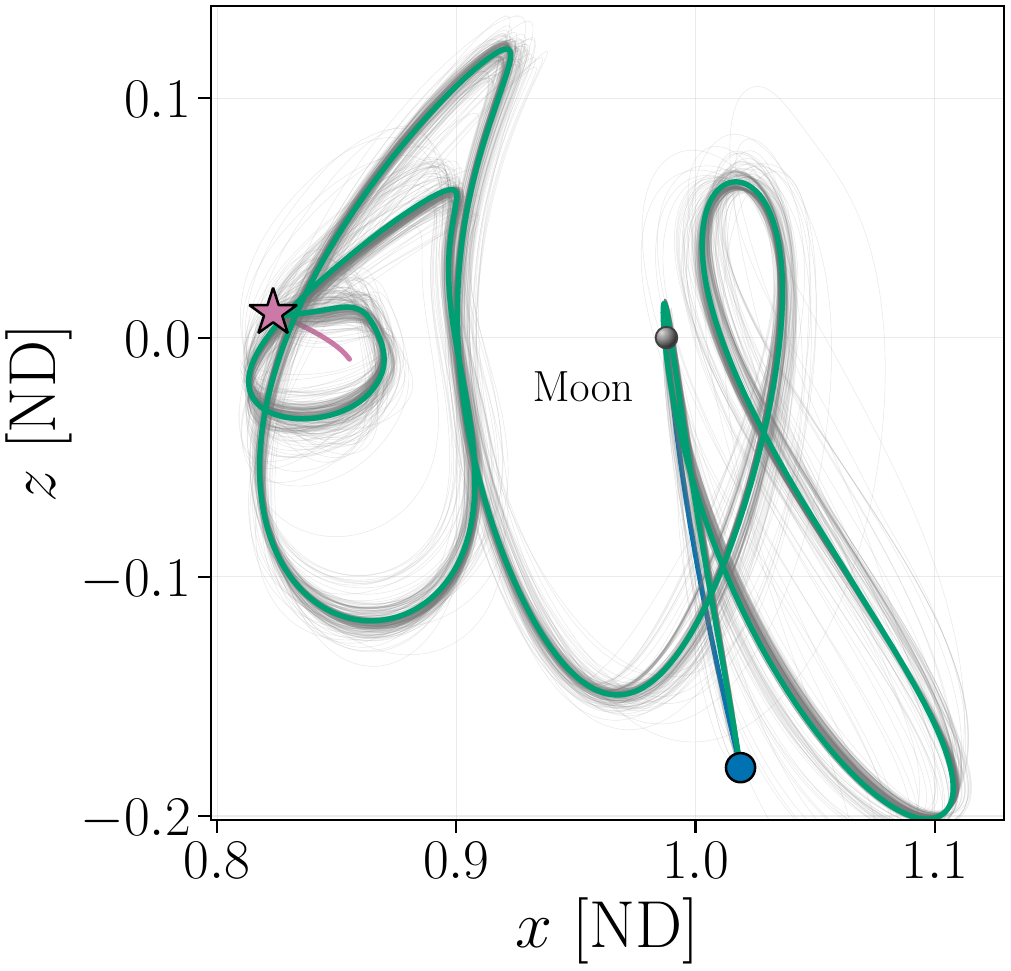}
        \caption{$xz$-projection.}
        \label{fig:nh_mc_xz}
    \end{subfigure}

    \vspace{0.6em}

    \begin{subfigure}[b]{0.48\textwidth}
        \centering
        \includegraphics[height=0.84\linewidth]{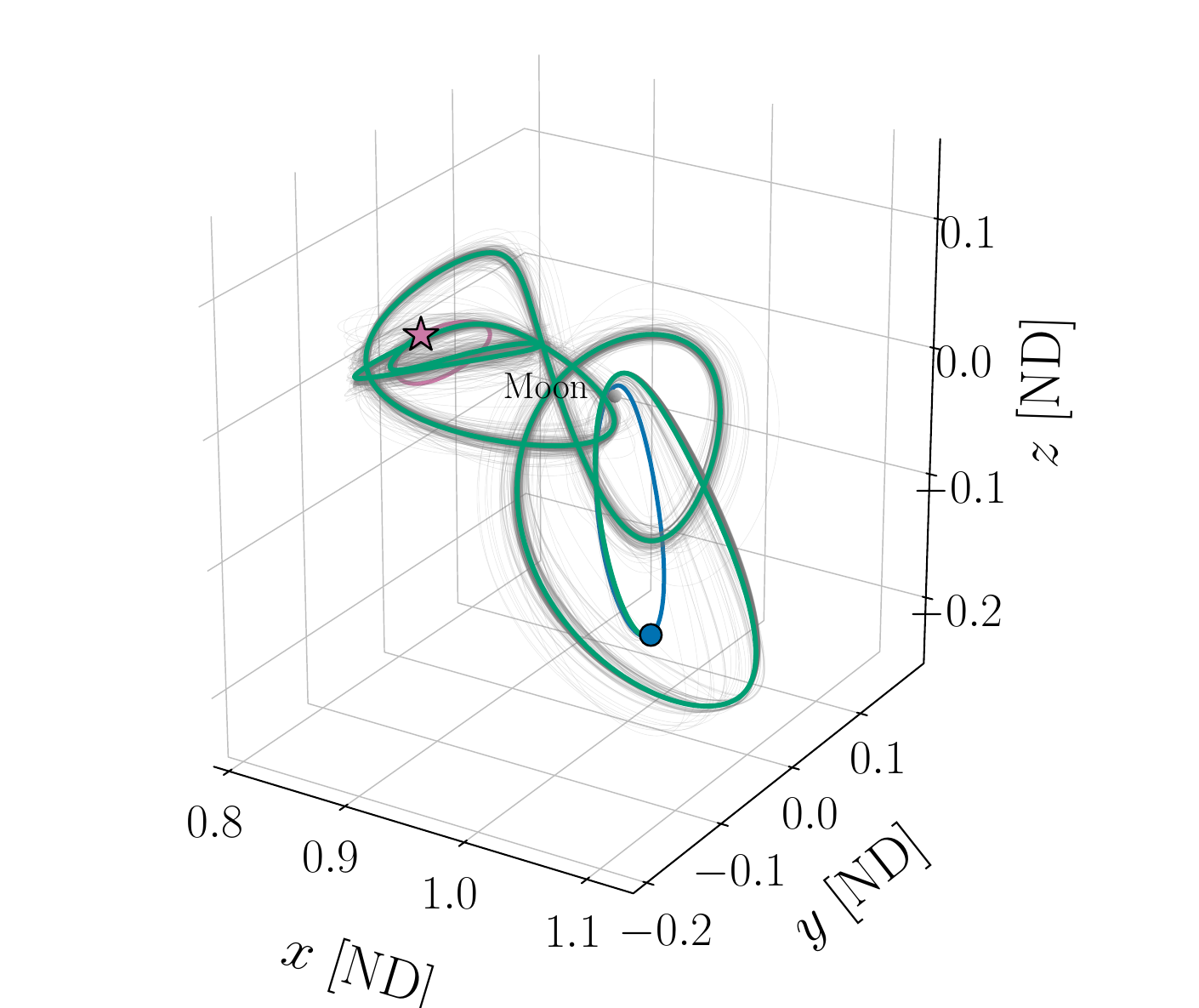}
        \caption{Three-dimensional view.}
        \label{fig:nh_mc_3d}
    \end{subfigure}
    \hfill
    \begin{subfigure}[b]{0.48\textwidth}
        \centering
        \includegraphics[height=0.84\linewidth]{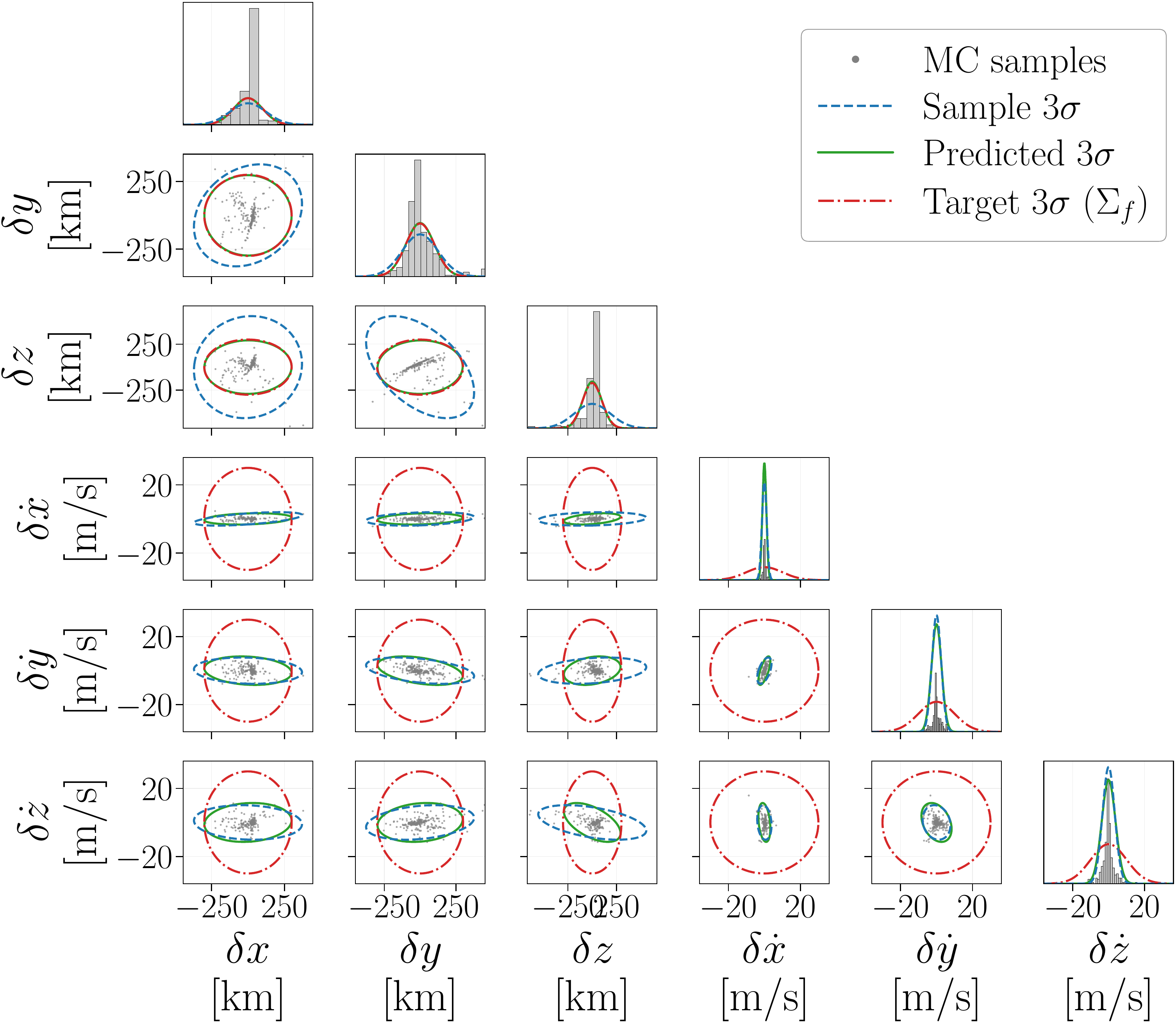}
        \caption{Terminal state dispersion.}
        \label{fig:nh_mc_disp}
    \end{subfigure}
    \caption{Monte Carlo results for the L2 NRHO to L1 Halo transfer; deviations amplified $5\times$ for trajectory plots.}
    \label{fig:nh_mc}
\end{figure}

\begin{figure}[htbp!]
    \centering
    \begin{subfigure}[b]{0.48\textwidth}
        \centering
        \includegraphics[height=0.74\linewidth]{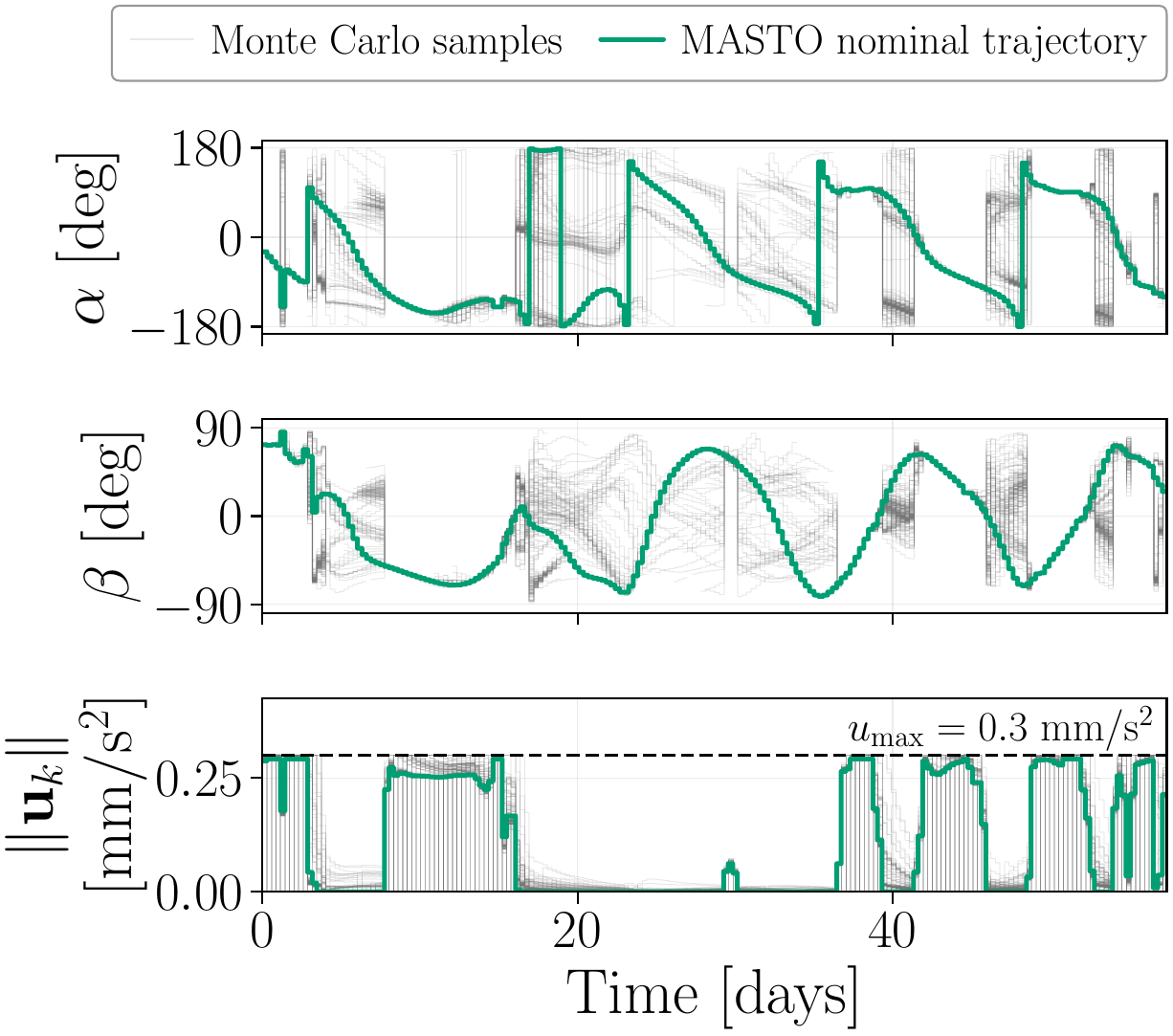}
        \caption{Monte Carlo control angles and magnitude.}
        \label{fig:nh_mc_ctrl}
    \end{subfigure}
    \hfill
    \begin{subfigure}[b]{0.48\textwidth}
        \centering
        \includegraphics[height=0.74\linewidth]{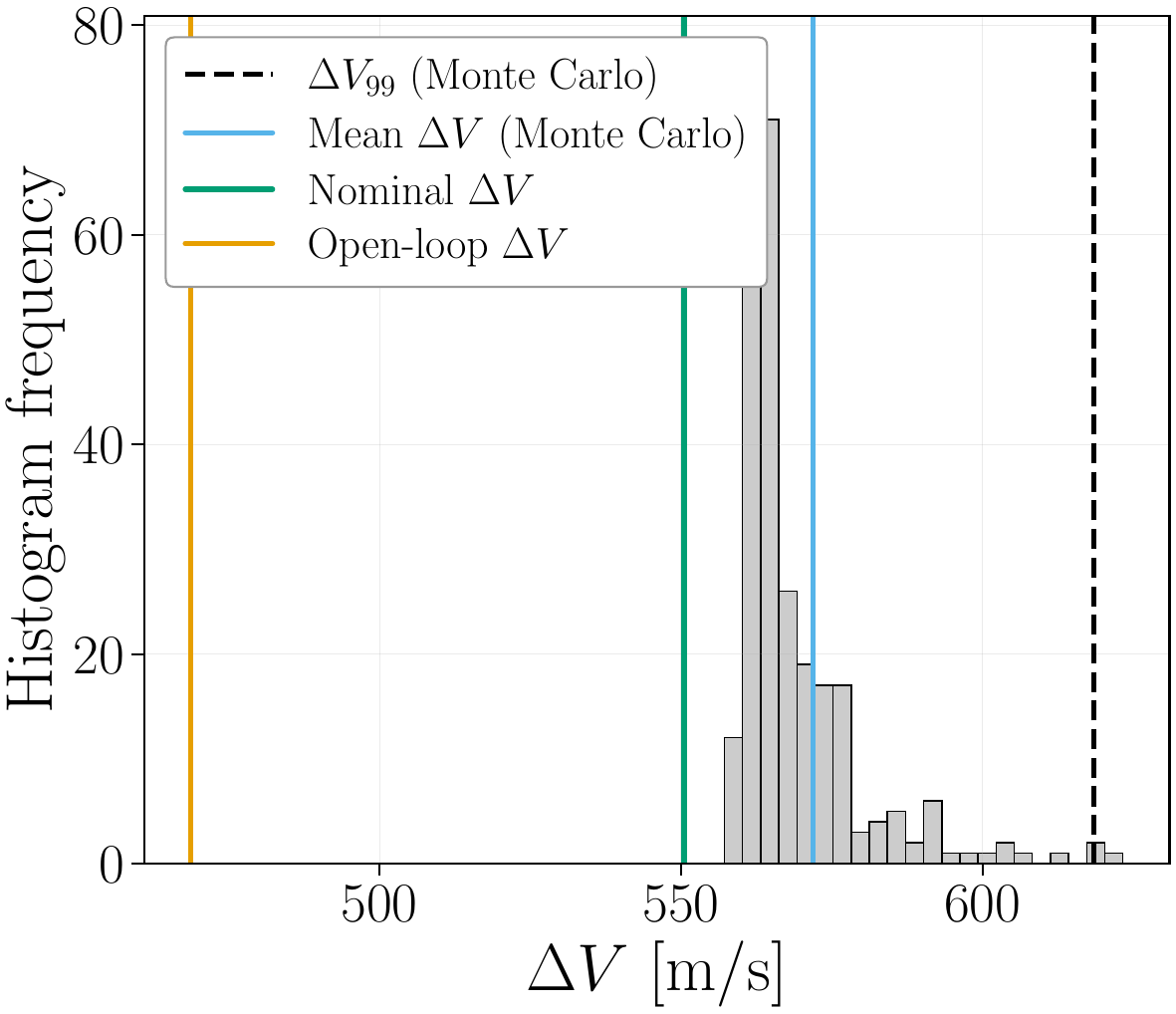}
        \caption{$\Delta v_{99}$ histogram.}
        \label{fig:nh_dv_hist}
    \end{subfigure}
    \caption{Monte Carlo control histories and $\Delta v$ distribution for the L2 NRHO to L1 Halo transfer.}
    \label{fig:nh_ctrl_dv}
\end{figure}

We now present results for $250$ samples of nonlinear Monte Carlo simulations under the same MTE distribution and present the results in \Cref{fig:nh_mc}. \Cref{fig:nh_mc_xy} and \Cref{fig:nh_mc_xz} show the $xy$ and $xz$ projections, respectively, \Cref{fig:nh_mc_3d} the three-dimensional trajectory, and \Cref{fig:nh_mc_disp} the terminal state dispersion. The Monte Carlo samples show good convergence to the desired orbit and the terminal dispersion plots show that the sample and predicted dispersion at final state agree well while mostly falling inside the prescribed bound. In position space, $91.6\%$ of samples fall within $3\sigma_{r_f}$ with one sample diverging (it exceeds $20\sigma_{r_f}$). The delivered terminal $\sigma_r$ is $106.1$ km against the $100$ km specification, and the one failed sample (with terminal position error $\geq20\sigma_{r_f}$ at final time) faces eight MTEs, with the first four MTEs coinciding with the second thrusting arc. This negatively impacted the critical pre-ballistic phase of the trajectory with more than 2 days of lost thrusting time. For visual clarity, we do not include this failed case in the Monte Carlo plots. We note that several other successful samples face more MTEs but recover under the MASTO policy. The corresponding control history plot for the Monte Carlo simulations is shown in \Cref{fig:nh_mc_ctrl}, with the $\Delta v_{99}$ histogram shown in \Cref{fig:nh_dv_hist}. The control angle history shows high variation when the nominal thrusting arc is absent, which implies the executed control action mostly comprises the corrective maneuvers. The MASTO nominal $\Delta v=550.5$ m/s is $17.5\%$ higher than the deterministic open-loop minimum-fuel reference $\Delta v=468.7$ m/s. This is in contrast to all the transfers we have considered before (the first example in this study, and the transfer examples in Reference \cite{shrivastavaMASTOMISSEDTHRUSTAWARE2026}), where the nominal trajectory introduced robustness with the MASTO framework without significant increase in the $\Delta v$ costs. This indicates that this L2 NRHO to L1 Halo transfer is very sensitive to MTEs, and introducing robustness requires significant increase in fuel consumption. The Monte Carlo $\Delta v_{99}$ cost was $618.4$ m/s, which is $31.9\%$ higher than the open-loop reference. This is also the largest increase from the open-loop $\Delta v$ costs to the $\Delta v_{99}$ costs we have seen in all the transfers considered before. 

\begin{figure}[htbp!]
    \centering
    \includegraphics[width=0.6\textwidth]{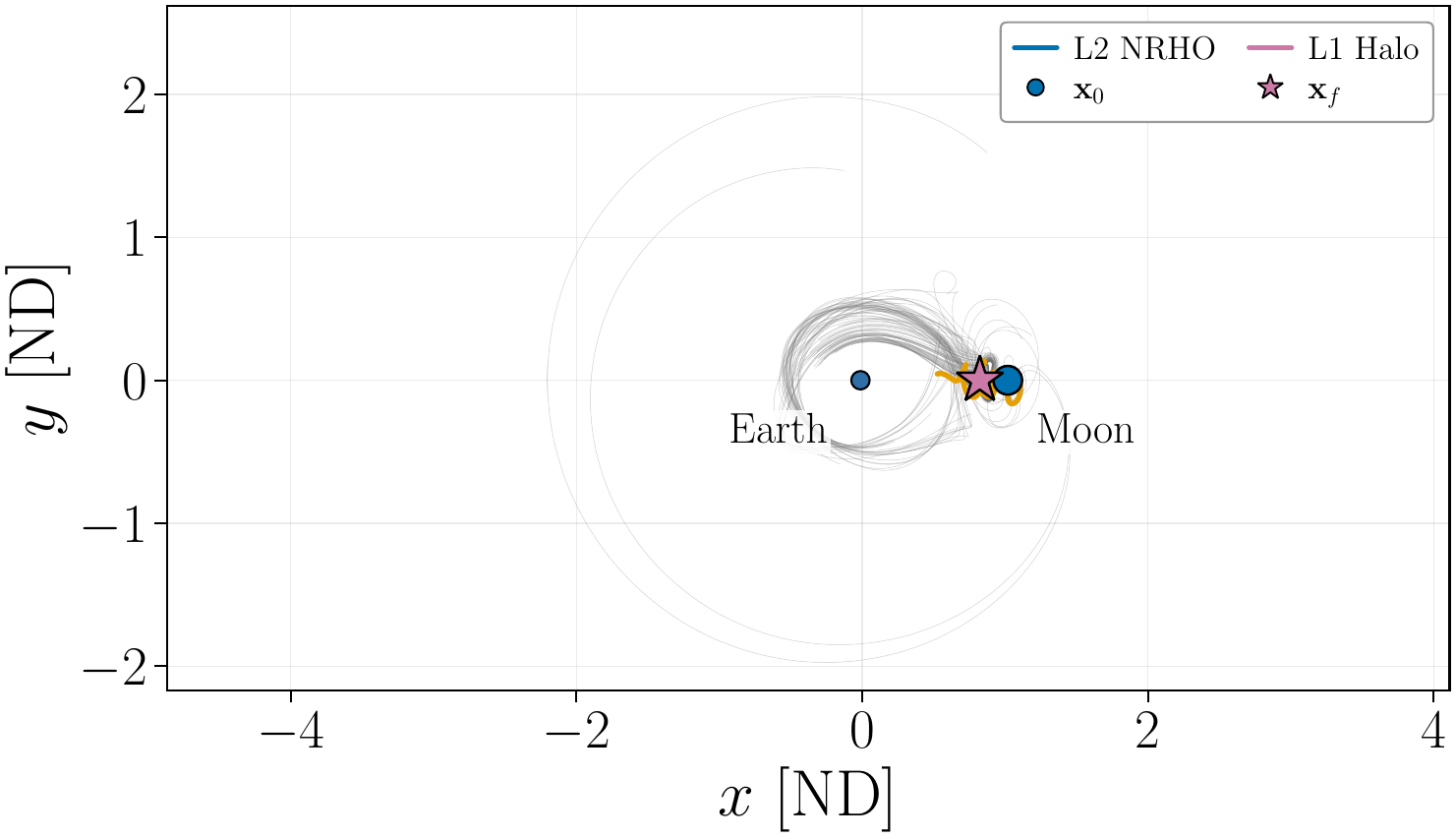}
    \caption{Single-mode Monte Carlo simulations ($xy$-projection) for L2 NRHO to L1 Halo transfer with MTEs.}
    \label{fig:nh_mc_singlemode_xy}
\end{figure}

\begin{figure}[htbp!]
    \centering
    \includegraphics[width=0.6\textwidth]{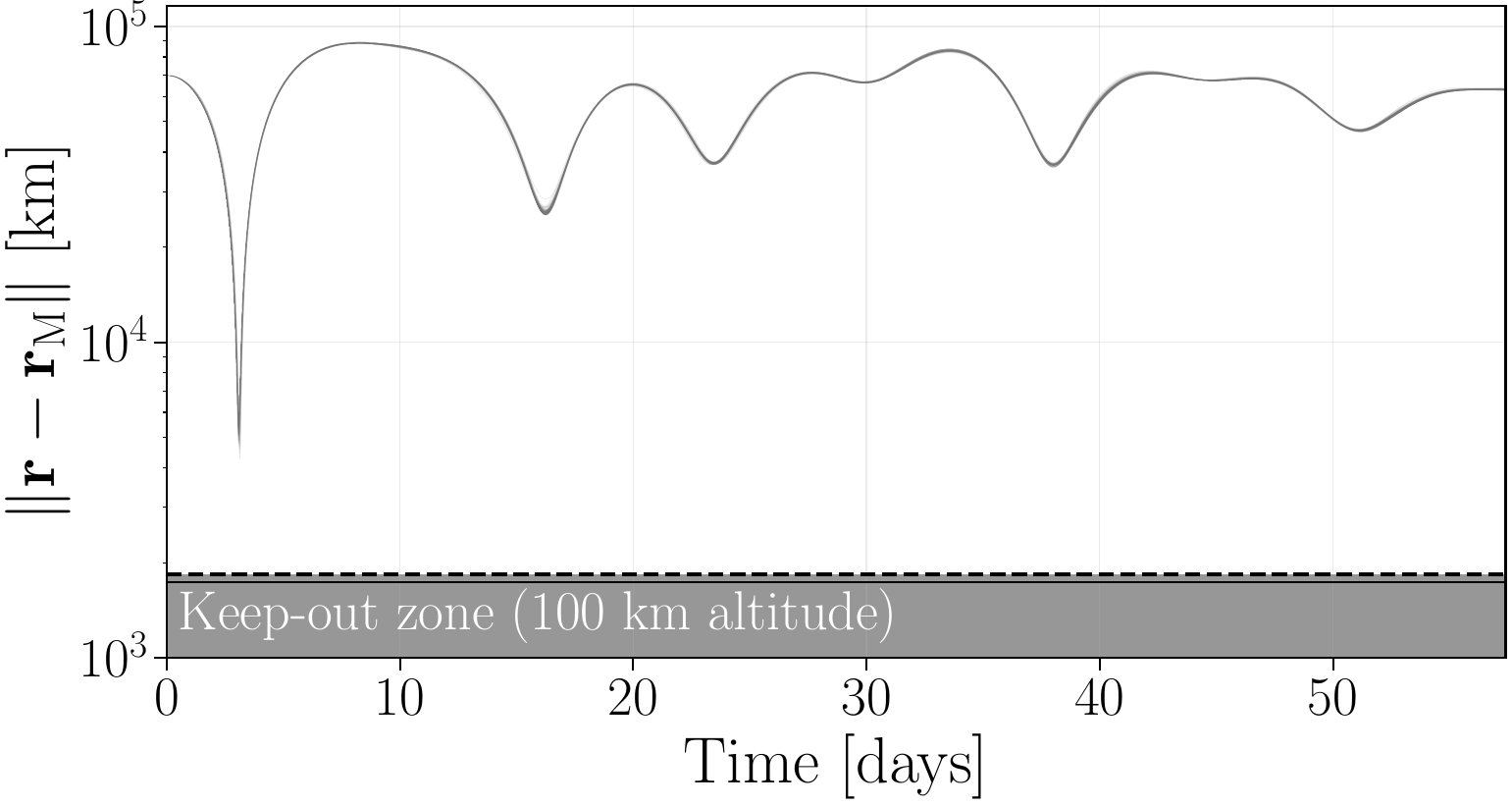}
    \caption{Distance from the Moon for all Monte Carlo samples under MASTO policy in the L2 NRHO to L1 Halo transfer. Shaded band: lunar keep-out zone; dashed line: $R_{\mathrm{safe}}$.}
    \label{fig:moon_distance_nh}
\end{figure}

For comparison, we also perform nonlinear simulations using the non-MTE aware single-mode covariance steering policy, under the same sequence of MTEs generated in the Monte Carlo simulations before. \Cref{fig:nh_mc_singlemode_xy} shows the $xy$ projection of the single-mode trajectories, with $143$ of the total $250$ samples exceeding $20\sigma_{r_f}$. While the feedback action recovers the spacecraft from some MTEs, it lacks the ability to account for an arbitrary number of MTEs which can occur at dynamically sensitive regions, leading to its poor performance. We emphasize that designing such trajectories with existing MTE-aware trajectory design methods, such as the virtual-swarm method~\cite{venigallaMultiObjectiveLowThrustTrajectory2022}, becomes computationally prohibitive, as it would lead to an exponential growth in the combinatorial possibilities to account for during the optimization process. Despite the close lunar flyby, the spacecraft satisfies the lunar surface clearance chance constraint as shown in \Cref{fig:moon_distance_nh}. Our results show that the MASTO trajectory was able to recover from multiple MTEs in this complex cislunar transfer that involved several thrusting arcs and a close lunar flyby, while respecting the control and state chance constraints. 

\subsection{Sensitivity to MTE Distribution}

While we use a given set of MTE probability parameters in this study for the two cislunar transfers considered, we emphasize that they are input parameters to the MASTO framework and can be changed according to mission requirements. The MASTO framework allows us to consider various MTE failure and recovery statistics thanks to its favorable scalability properties. To further analyze the effect of these parameters, we compare MASTO trajectories for five different values of $\text{TTF}_{90}\in\{15,20,30,40,50\} \text{ days}$ and $\text{TTR}_{95}\in\{6,10,12,18,24\}\text{ hours}$ with three different control authority values $u_{\max}\in\{0.2,0.3,0.5\}$ mm/s$^2$ for the L2 Halo to L1 Halo transfer. Therefore, we solve $5\times5\times3=75$ optimization problems, and perform $250\times75=18750$ nonlinear Monte Carlo simulations. This trade study takes approximately fourteen hours to complete. 

\begin{figure}[htbp!]
    \centering
    \includegraphics[width=\linewidth]{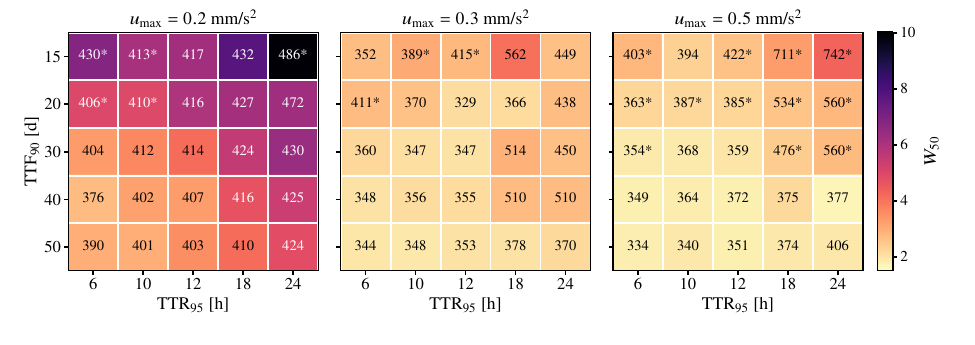}
    \caption{L2 Halo to L1 Halo trade study: cell color is $W_{50}$, annotation is $\Delta v_{99}$ in m/s; an asterisk marks cells where the SCP did not satisfy the optimality convergence criteria within the iteration limit.}
    \label{fig:hh_ts_w}
\end{figure}

\Cref{fig:hh_ts_w} reports the results of our trade study, where the colorbar denotes a capped mean Mahalanobis distance of the terminal state, $W_{50}=\frac{1}{N}\sum_{i=1}^{N}\min\left(d_i,50\right)$ with $d_i=\sqrt{\delta\mathbf{x}_i^\top\Sigma_f^{-1}\delta\mathbf{x}_i}$, where $\delta\mathbf{x}_i$ is the terminal deviation of sample $i$ from the planned terminal mean and $\Sigma_f$ is the planned terminal covariance. The numbers in each cell indicate the total $\Delta v_{99}$ fuel consumption for that case. 

Across the three grids in \Cref{fig:hh_ts_w}, we note that the grid in top right corner in the left grid with lowest control authority $u_{\max}=0.2 \text{ mm/s}^2$ represents the most acute case with lowest $\text{TTF}_{90}$ and highest $\text{TTR}_{95}$ values. This case represents statistically most frequent MTEs with the longest recovery times, while they have least control authority. Unsurprisingly, this is where the dispersion peaks, as shown by our metric $W_{50}$. The dispersion trends are generally higher for the lower control authority cases, while they are generally the lowest for $u_{\max}=0.5\text{ mm/s}^2$. 

For $\Delta v_{99}$ trends, we note that it is not fair to compare across the three different grids. This is because the higher control authority might lead to more built-in conservatism, but the higher $u_{\max}$ will itself lead to higher total fuel consumption. However, for a given $u_{\max}$, we see a clear trend that is generally increasing with longer MTEs and longer recovery times. For instance, in the $u_{\max}=0.5\text{ mm/s}^2$ case, it is most pronounced, where total $\Delta v_{99}$ cost more than doubles from $334\text{ m/s}$ to $742 \text{ m/s}$ as we go from $\{\text{TTF}_{90},\text{TTR}_{95}\}=\{50 \text{ days},6\text{ hours}\}$ case to $\{\text{TTF}_{90},\text{TTR}_{95}\}=\{15 \text{ days},24\text{ hours}\}$. Similar trends are also observed for $u_{\max}=0.2 \text{ mm/s}^2$ and $u_{\max}=0.3 \text{ mm/s}^2$, although to a lesser effect. We also note some outliers in the trade study, such as the $\{\text{TTF}_{90},\text{TTR}_{95}\}=\{15 \text{ days},18\text{ hours}\}$ case for $u_{\max}=0.3 \text{ mm/s}^2$ which had the highest dispersion and $\Delta v_{99}$ for that grid. One possible explanation for such outliers is that our original MTE-resilient trajectory optimization problem is a highly nonconvex problem. Like most optimization frameworks, MASTO can also get stuck in a local minima, leading to these outliers. They can also be affected by number of samples proving insufficient for statistical convergence of the trends, and greater number of samples can lead to more consistent trends with lesser outliers.

We note that while all $75$ cells yield a feasible solution, twenty of them (marked with an asterisk in \Cref{fig:hh_ts_w}) do not satisfy the SCP optimality convergence criteria within the iteration limit. We do not perform any \emph{case-by-case} tuning of SCP parameters in any of these 75 cases, and we use the same SCP parameters as our original baseline shown in \Cref{tab:scvx_params}. This showcases MASTO convergence behavior's robustness to SCP parameters, requiring no additional tuning to produce feasible trajectories for such extensive trade studies. The scalability of the MASTO framework allows us to perform such trade studies across a range of different MTE distribution parameters. This can offer valuable information about how the MTE-resilient trajectory optimization problem behaves with respect to MTE modeling parameters. This is a valuable tool for mission designers, as they trade MTE-resilience with fuel optimality depending on their design requirements. To the authors’ best knowledge, the existing MTE-resilient frameworks is not scalable to account for such rich trade studies. Consequently, the common practice in mission design has resorted to expensive Monte Carlo simulations, with the need for human-in-the-loop heuristics that lack safety guarantees.

\begin{remark}
    In some cases, mission designers may wish to consider time-varying MTE probabilities. For instance, failure events can be more likely in the initial post-launch phase and become less frequent as the mission settles into routine operations. While the numerical simulations in this work assume a time-homogeneous MTE Markov chain, the MASTO framework accommodates time-varying models directly: introducing time-varying failure and recovery rates in \cref{tpm_exp} yields a time-varying TPM $\mathcal{P}_k$, while the rest of the formulation remains unchanged. 
\end{remark}
    
\begin{remark}    
    Some missions require more risk-averse trajectory design, with increased resilience against MTEs during sensitive thrusting arcs such as flybys. Standard approaches for such problems consider 1-2 MTEs at such critical times, and leverage scenario-based optimization or Monte Carlo simulation-based design\cite{venigallaMultiObjectiveLowThrustTrajectory2022,sekineRobustTrajectoryDesign2026,grebowDawnSafingApproach2016,mccartyMissedThrustAnalysis2020}. However, such methods suffer from a combinatorial growth in problem size as we consider different scenarios. The MASTO framework can define such problems directly, by modeling the MTE probability as near unity close to the flyby and near zero elsewhere. MASTO therefore generalizes existing MTE-aware trajectory optimization frameworks: it can model the more risk-averse approach, which considers only a few critical MTEs but demands increased robustness, as well as the stochastic approach adopted in this study, which accounts for an arbitrary number of MTEs throughout the mission.
\end{remark}

\section{Conclusions}

In this work, we propose a stochastic optimal control-based framework for robust spacecraft trajectory optimization in the presence of multiple missed thrust events (MTEs). The Missed-thrust Aware Stochastic Trajectory Optimization (MASTO) framework extends chance-constrained covariance steering theory to nonlinear Markov Jump Systems using a SCP-based algorithm. The MASTO framework allows us to generate MTE-resilient trajectories that are able to recover from multiple stochastic missed thrust events. We demonstrate the framework's resilience against multiple MTEs in complex cislunar transfers involving lunar flybys. Further, our study also outlines controller performance across a range of MTE failure and recovery distributions. Therefore, the MASTO framework presents a scalable approach for practitioners to design robust nominal trajectories with guaranteed feasible recovery strategies (through corrective feedback maneuvers) against an arbitrary number of MTEs, and explore the trade space for MTE-resilient complex transfer design problems.




\section*{Funding Sources}

\noindent{This work was supported by the U.S. Air Force Office of Scientific Research through research grant
FA9550-23-1-0512.}
\section*{Acknowledgments}
 The authors thank Yanis Sidhoum and Naoya Kumagai for providing the minimum-fuel reference trajectories for the L2 Halo to L1 Halo transfer and the L2 NRHO to L1 Halo transfer, respectively.  

We used AI tools for coding assistance, schematic drawing, and grammar correction. However, all content was edited and checked by the authors, and we take full responsibility for the content of the paper.

\bibliography{references2}

@inproceedings{nuccioStatisticalTrajectoryOptimization2026,
	title = {Statistical {Trajectory} {Optimization} {Under} {Uncertainty} for {Icy} {Moon} {Exploration} with {Probability} of {Impact} {Constraints}},
	url = {https://re.public.polimi.it/bitstream/11311/1317067/1/NUCCG01-26.pdf},
	urldate = {2026-11-20},
	booktitle = {30th {International} {Symposium} on {Space} {Flight} {Dynamics} ({ISSFD})},
	author = {Nuccio, G. and Oguri, K. and Topputo, F.},
	year = {2026},
	pages = {1--11},
}

@inproceedings{SafeAutonomousGuidance,
	title = {Safe {Autonomous} {Guidance} {Under} {Uncertainty} via {Sequential} {Covariance} {Steering} for {Small}-{Body} {Proximity} {Operations}},
	language = {en},
	urldate = {2026-11-20},
	booktitle = {30th {International} {Symposium} on {Space} {Flight} {Dynamics} ({ISSFD})},
	author = {Frassinella, Luca and Oguri, Kenshiro and Topputo, Francesco},
	year = {2026},
}

@article{alkhenaniChanceConstrainedAutonomousSpacecraft,
	title = {Chance-{Constrained} {Autonomous} {Spacecraft} {Formation} {Reconfiguration} and {Rendezvous} with {Passive} {Abort} {Safety}},
	volume = {0},
	issn = {0731-5090},
	url = {https://doi.org/10.2514/1.G010121},
	doi = {10.2514/1.G010121},
	number = {0},
	urldate = {2026-11-20},
	journal = {Journal of Guidance, Control, and Dynamics},
	publisher = {American Institute of Aeronautics and Astronautics},
	author = {Alkhenani, Abdullah K. and Oguri, Kenshiro},
	pages = {1--18},
}

@incollection{shibukawaDDPBasedStochasticTrajectory,
	title = {{DDP}-{Based} {Stochastic} {Trajectory} {Optimization} {With} {Safe} {Mode} {Uncertainty} {Modeled} by an {Alternative} {Renewal} {Process}},
	doi = {10.2514/6.2026-1461},
	urldate = {2026-10-21},
	booktitle = {{AIAA} {SCITECH} 2026 {Forum}},
	publisher = {American Institute of Aeronautics and Astronautics},
	author = {Shibukawa, Masato and Ozaki, Naoya},
	year = {2026},
}

@article{onoProbabilisticPlanningContinuous2013,
	title = {Probabilistic {Planning} for {Continuous} {Dynamic} {Systems} under {Bounded} {Risk}},
	volume = {46},
	issn = {1076-9757},
	doi = {10.1613/jair.3893},
	language = {en},
	urldate = {2025-02-27},
	journal = {J. of Artif. Intell. Res.},
	author = {Ono, M. and Williams, B. C. and Blackmore, Lars},
	month = mar,
	year = {2013},
	pages = {511--577},
}

@article{sidhoumPontryaginBellmanDifferential2026,
	title = {Pontryagin–{Bellman} differential dynamic programming for low-thrust trajectory optimization with path constraints},
	volume = {247},
	issn = {00945765},
	doi = {10.1016/j.actaastro.2026.04.067},
	language = {en},
	urldate = {2026-10-20},
	journal = {Acta Astronautica},
	author = {Sidhoum, Yanis and Oguri, Kenshiro},
	month = oct,
	year = {2026},
	pages = {10--24},
}

@misc{qiNonGaussianDistributionSteering2025,
	title = {Non-{Gaussian} {Distribution} {Steering} in {Nonlinear} {Dynamics} with {Conjugate} {Unscented} {Transformation}},
	doi = {10.48550/arXiv.2510.12946},
	urldate = {2026-03-29},
	publisher = {arXiv},
	author = {Qi, Daniel C. and Oguri, Kenshiro and Singla, Puneet and Akella, Maruthi R.},
	month = oct,
	year = {2025},
	note = {arXiv:2510.12946 [eess]},
}

@article{apsMOSEKOptimizerAPI,
	title = {{MOSEK} {Optimizer} {API} for {Python}},
	language = {en},
	author = {ApS, MOSEK},
	year = {2022},
}

@inproceedings{oguriSuccessiveConvexificationFeasibility,
	title = {Successive {Convexification} with {Feasibility} {Guarantee} via {Augmented} {Lagrangian} for {Non}-{Convex} {Optimal} {Control} {Problems}},
	doi = {10.1109/CDC49753.2023.10383462},
	urldate = {2025-03-12},
	booktitle = {2023 62nd {IEEE} {Conf}. on {Decis}. and {Control} ({CDC})},
	author = {Oguri, Kenshiro},
	year = {2023},
}

@book{kullbackInformationTheoryStatistics1997,
	title = {Information theory and statistics},
	urldate = {2025-03-14},
	publisher = {Courier Corp},
	author = {Kullback, Solomon},
	year = {1997},
}

@inproceedings{rapakouliasDiscreteTimeOptimalCovariance,
	title = {Discrete-{Time} {Optimal} {Covariance} {Steering} via {Semidefinite} {Programming}},
	issn = {2576-2370},
	doi = {10.1109/CDC49753.2023.10384118},
	urldate = {2024-11-07},
	booktitle = {2023 62nd {IEEE} {Conf}. on {Decis}. and {Control} ({CDC})},
	author = {Rapakoulias, George and Tsiotras, Panagiotis},
	year = {2023},
	pages = {1802--1807},
}

@article{liuOptimalCovarianceSteering2024,
	title = {Optimal {Covariance} {Steering} for {Discrete}-{Time} {Linear} {Stochastic} {Systems}},
	issn = {1558-2523},
	doi = {10.1109/TAC.2024.3472788},
	urldate = {2025-01-31},
	journal = {IEEE Trans. on Autom. Control},
	author = {Liu, Fengjiao and Rapakoulias, George and Tsiotras, Panagiotis},
	year = {2024},
	pages = {1--16},
}

@article{carpentierLowThrustInterplanetaryTrajectories,
	title = {Low-{Thrust} {Interplanetary} {Trajectories} with {Missed} {Thrust} {Events}: {A} {Numerical} {Approach}},
	volume = {0},
	issn = {0731-5090},
	shorttitle = {Low-{Thrust} {Interplanetary} {Trajectories} with {Missed} {Thrust} {Events}},
	doi = {10.2514/1.G009777},
	number = {0},
	urldate = {2026-07-03},
	journal = {Journal of Guidance, Control, and Dynamics},
	publisher = {American Institute of Aeronautics and Astronautics},
	author = {Carpentier, Pierre and Chancelier, Jean-Philippe and Cohen, Guy and Dargent, Thierry and Epenoy, Richard},
	year = {2026},
	pages = {1--18},
}

@article{zavoliReinforcementLearningRobust2021,
	title = {Reinforcement {Learning} for {Robust} {Trajectory} {Design} of {Interplanetary} {Missions}},
	volume = {44},
	issn = {0731-5090},
	doi = {10.2514/1.G005794},
	number = {8},
	urldate = {2023-10-17},
	journal = {Journal of Guidance, Control, and Dynamics},
	publisher = {American Institute of Aeronautics and Astronautics},
	author = {Zavoli, Alessandro and Federici, Lorenzo},
	year = {2021},
	pages = {1440--1453},
}

@inproceedings{olympioDeterministicMethodSpace,
	title = {Deterministic {Method} for {Space} {Trajectory} {Design} with {Mission} {Margin} {Constraints}},
	language = {en},
	booktitle = {61st {International} {Astronautical} {Congress},},
	author = {Olympio, Joris T},
	year = {2010},
}

@inproceedings{olympioDESIGNINGROBUSTLOWTHRUST,
	title = {Designing {Robust} {Low}-{Thrust} {Interplanetary} {Trajectories} {Subject} {To} {One} {Temporary} {Engine} {Failure}},
	language = {en},
	booktitle = {20th {AAS}/{AIAA} {Space} {Flight} {Meeting}},
	author = {Olympio, Joris T},
	year = {2010},
}

@article{raymanDesignFirstInterplanetary2002,
	title = {Design of the {First} {Interplanetary} {Solar} {Electric} {Propulsion} {Mission}},
	volume = {39},
	issn = {0022-4650},
	doi = {10.2514/2.3848},
	number = {4},
	urldate = {2023-04-07},
	journal = {Journal of Spacecraft and Rockets},
	publisher = {American Institute of Aeronautics and Astronautics},
	author = {Rayman, Marc D. and Williams, Steven N.},
	year = {2002},
	pages = {589--595},
}

@incollection{ohAnalysisSystemMargins,
	title = {Analysis of {System} {Margins} on {Deep} {Space} {Missions} {Using} {Solar} {Electric} {Propulsion}},
	doi = {10.2514/6.2008-5286},
	urldate = {2023-04-07},
	booktitle = {44th {AIAA}/{ASME}/{SAE}/{ASEE} {Joint} {Propulsion} {Conference} \& {Exhibit}},
	publisher = {American Institute of Aeronautics and Astronautics},
	author = {Oh, David and Landau, Damon and Randolph, Tom and Timmerman, Paul and Chase, James and Sims, John and Kowalkowski, Theresa},
	year = {2008},
}

@misc{shrivastavaRobustHaloOrbit2026,
	title = {Robust {Halo} {Orbit} {Stationkeeping} {Under} {Missed} {Thrust} {Events} and {Maneuver} {Execution} {Errors}},
	url = {http://arxiv.org/abs/2610.12906},
	doi = {10.48550/arXiv.2610.12906},
	urldate = {2026-10-21},
	publisher = {arXiv},
	author = {Shrivastava, Shaurya and Oguri, Kenshiro},
	month = oct,
	year = {2026},
	note = {arXiv:2610.12906 [math.OC]},
}

@inproceedings{sekineRobustTrajectoryDesign2026,
	title = {Robust {Trajectory} {Design} for {Earth}-to-lunar {Orbit} {Transfers} via {Sequential} {Convex} {Programming}},
	booktitle = {{AAS} {Guidance}, {Navigation}, and {Control} {Conference}},
	author = {Sekine, Hirotaka and Oguri, Kenshiro and Kawabata, Yosuke and Funase, Ryu},
	month = feb,
	year = {2026},
}

@article{rubinsztejnNeuralNetworkOptimal2020,
	title = {Neural network optimal control in astrodynamics: {Application} to the missed thrust problem},
	volume = {176},
	issn = {0094-5765},
	shorttitle = {Neural network optimal control in astrodynamics},
	doi = {10.1016/j.actaastro.2020.05.027},
	language = {en},
	urldate = {2023-04-07},
	journal = {Acta Astronautica},
	author = {Rubinsztejn, Ari and Sood, Rohan and Laipert, Frank E.},
	month = nov,
	year = {2020},
	pages = {192--203},
}

@incollection{grebowDawnSafingApproach2016,
	series = {{AIAA} {SPACE} {Forum}},
	title = {Dawn {Safing} {Approach} to {Ceres} {Re}-{Design}},
	doi = {10.2514/6.2016-5426},
	urldate = {2025-01-08},
	booktitle = {{AIAA}/{AAS} {Astrodynamics} {Specialist} {Conference}},
	author = {Grebow, Daniel and Whiffen, Gregory J. and Han, Dongsuk and Kennedy, Brian},
	month = sep,
	year = {2016},
}

@inproceedings{shrivastavaMASTOMissedthrustAware2026,
	title = {{MASTO}: {Missed}-thrust {Aware} {Stochastic} {Trajectory} {Optimization}},
	shorttitle = {{MASTO}},
	booktitle = {2026 {AAS}/{AIAA} {Astrodynamics} {Specialist} {Conference}},
	author = {Shrivastava, Shaurya and Oguri, Kenshiro},
	month = jul,
	year = {2026},
}

@article{sinhaInitialGuessGeneration2025a,
	title = {Initial {Guess} {Generation} for {Low}-{Thrust} {Trajectory} {Design} with {Robustness} to {Missed}-{Thrust} {Events}},
	volume = {48},
	issn = {0731-5090},
	doi = {10.2514/1.G009050},
	number = {12},
	urldate = {2026-03-28},
	journal = {Journal of Guidance, Control, and Dynamics},
	publisher = {American Institute of Aeronautics and Astronautics},
	author = {Sinha, Amlan and Beeson, Ryne},
	year = {2025},
	pages = {2831--2848},
}

@article{rubinsztejnDesigningTrajectoriesResilient2021,
	title = {Designing trajectories resilient to missed thrust events using expected thrust fraction},
	volume = {115},
	issn = {1270-9638},
	doi = {10.1016/j.ast.2021.106780},
	language = {en},
	urldate = {2023-04-07},
	journal = {Aerospace Science and Technology},
	author = {Rubinsztejn, Ari and Sandel, Carrie Grace and Sood, Rohan and Laipert, Frank E.},
	month = aug,
	year = {2021},
	pages = {106780},
}

@inproceedings{oguriChanceConstrainedControlSafe2024,
	title = {Chance-{Constrained} {Control} for {Safe} {Spacecraft} {Autonomy}: {Convex} {Programming} {Approach}},
	issn = {2378-5861},
	shorttitle = {Chance-{Constrained} {Control} for {Safe} {Spacecraft} {Autonomy}},
	doi = {10.23919/ACC60939.2024.10645008},
	urldate = {2025-01-08},
	booktitle = {2024 {Amer}. {Control} {Conf}. ({ACC})},
	author = {Oguri, Kenshiro},
	month = jul,
	year = {2024},
	pages = {2318--2324},
}

@article{laipertAutomatedMissedThrustPropellant2015,
	title = {Automated {Missed}-{Thrust} {Propellant} {Margin} {Analysis} for {Low}-{Thrust} {Trajectories}},
	volume = {52},
	issn = {0022-4650},
	doi = {10.2514/1.A33264},
	number = {4},
	urldate = {2024-09-23},
	journal = {Journal of Spacecraft and Rockets},
	publisher = {American Institute of Aeronautics and Astronautics},
	author = {Laipert, Frank E. and Longuski, James M.},
	year = {2015},
	pages = {1135--1143},
}

@article{nemirovskyTrajectoryDesignControl2026,
	title = {Trajectory {Design} and {Control} for {Missed}-{Thrust} {Rendezvous}},
	volume = {49},
	issn = {0731-5090},
	doi = {10.2514/1.G008831},
	number = {4},
	urldate = {2026-08-31},
	journal = {Journal of Guidance, Control, and Dynamics},
	publisher = {American Institute of Aeronautics and Astronautics},
	author = {Nemirovsky, Meir and Gurfil, Pini},
	month = apr,
	year = {2026},
	pages = {990--1004},
}

@article{kumagaiRobustCislunarLowThrust,
	title = {Robust {Cislunar} {Low}-{Thrust} {Trajectory} {Optimization} {Under} {Uncertainties} via {Sequential} {Covariance} {Steering}},
	volume = {48},
	issn = {0731-5090},
	doi = {10.2514/1.G009092},
	number = {12},
	urldate = {2025-09-24},
	journal = {Journal of Guidance, Control, and Dynamics},
	author = {Kumagai, Naoya and Oguri, Kenshiro},
	month = dec,
	year = {2025},
	pages = {1--19},
}

@inproceedings{imkenModelingSpacecraftSafe2018,
	title = {Modeling spacecraft safe mode events},
	doi = {10.1109/AERO.2018.8396383},
	urldate = {2024-04-26},
	booktitle = {2018 {IEEE} {Aerospace} {Conference}},
	author = {Imken, Travis and Randolph, Thomas and DiNicola, Michael and Nicholas, Austin},
	month = mar,
	year = {2018},
	pages = {1--13},
}

@misc{wangCovarianceSteeringDiscreteTime2026,
	title = {Covariance {Steering} of {Discrete}-{Time} {Markov} {Jump} {Linear} {Systems} with {Multiplicative} {Noise}},
	doi = {10.48550/arXiv.2604.19994},
	urldate = {2026-05-11},
	publisher = {arXiv},
	author = {Wang, Fangji and Ganguly, Siddhartha and Tsiotras, Panagiotis},
	month = apr,
	year = {2026},
	note = {arXiv:2604.19994 [math.OC]},
}

@inproceedings{ridderhofChanceconstrainedCovarianceControl2020,
	title = {Chance-constrained covariance control for low-thrust minimum-fuel trajectory optimization},
	urldate = {2025-03-13},
	booktitle = {{AAS}/{AIAA} {Astrodynamics} {Specialist} {Conf}.},
	author = {Ridderhof, Jack and Pilipovsky, Joshua and Tsiotras, Panagiotis},
	year = {2020},
	pages = {9--13},
}

@article{nathersonReachabililtybasedMissedThrust2026,
	title = {Reachabililty-based missed thrust analysis},
	volume = {247},
	issn = {00945765},
	url = {https://linkinghub.elsevier.com/retrieve/pii/S0094576526003413},
	doi = {10.1016/j.actaastro.2026.05.030},
	language = {en},
	urldate = {2026-07-03},
	journal = {Acta Astronautica},
	author = {Natherson, Robyn A. and Scheeres, Daniel J.},
	month = oct,
	year = {2026},
	pages = {565--578},
}

@inproceedings{oguriConvexBlockCholeskyApproach2026,
	title = {Convex {Block}-{Cholesky} {Approach} to {Risk}-{Constrained} {Low}-thrust {Trajectory} {Design} under {Operational} {Uncertainty}},
	url = {https://arxiv.org/abs/2602.18416v1},
	language = {en},
	urldate = {2026-03-29},
	booktitle = {{arXiv}.org},
	author = {Oguri, Kenshiro and Lantoine, Gregory},
	month = feb,
	year = {2026},
}

@article{ozakiStochasticDifferentialDynamic2018,
	title = {Stochastic differential dynamic programming with unscented transform for low-thrust trajectory design},
	volume = {41},
	issn = {0731-5090},
	number = {2},
	journal = {Journal of Guidance, Control, and Dynamics},
	publisher = {American Institute of Aeronautics and Astronautics},
	author = {Ozaki, Naoya and Campagnola, Stefano and Funase, Ryu and Yam, Chit Hong},
	year = {2018},
	pages = {377--387},
}

@misc{sinhaBiLevelOptimalControl2025,
	title = {Bi-{Level} {Optimal} {Control} {Framework} {For} {Missed}-{Thrust}-{Design} {With} {First}-{Order} {Bounds} {On} {Maximum} {Missed}-{Thrust}-{Duration}},
	url = {https://arxiv.org/abs/2512.18984v1},
	language = {en},
	urldate = {2026-03-28},
	journal = {arXiv.org},
	author = {Sinha, Amlan and Beeson, Ryne},
	month = dec,
	year = {2025},
}

@article{levTechnologicalCommercialExpansion2019,
	title = {The technological and commercial expansion of electric propulsion},
	volume = {159},
	issn = {0094-5765},
	url = {https://www.sciencedirect.com/science/article/pii/S0094576518319672},
	doi = {10.1016/j.actaastro.2019.03.058},
	urldate = {2026-03-28},
	journal = {Acta Astronautica},
	author = {Lev, Dan and Myers, Roger M. and Lemmer, Kristina M. and Kolbeck, Jonathan and Koizumi, Hiroyuki and Polzin, Kurt},
	month = jun,
	year = {2019},
	pages = {213--227},
}

@misc{elangoSuccessiveConvexificationPassivelySafe2025,
	title = {Successive {Convexification} for {Passively}-{Safe} {Spacecraft} {Rendezvous} on {Near} {Rectilinear} {Halo} {Orbit}},
	url = {http://arxiv.org/abs/2505.17251},
	doi = {10.48550/arXiv.2505.17251},
	urldate = {2026-01-07},
	publisher = {arXiv},
	author = {Elango, Purnanand and Vinod, Abraham P. and Kitamura, Kenji and Açıkmeşe, Behçet and Cairano, Stefano Di and Weiss, Avishai},
	month = may,
	year = {2025},
	note = {arXiv:2505.17251 [math]},
}

@book{costaDiscreteTimeMarkovJump2005,
	title = {Discrete-{Time} {Markov} {Jump} {Linear} {Systems}},
	isbn = {978-1-85233-761-2},
	language = {en},
	publisher = {Springer Sci. \& Bus. Media},
	author = {Costa, O. L. V. and Fragoso, M. D. and Marques, R. P.},
	month = feb,
	year = {2005},
}

@misc{shrivastavaChanceConstrainedCovarianceSteering2025,
	title = {Chance-{Constrained} {Covariance} {Steering} for {Discrete}-{Time} {Markov} {Jump} {Linear} {Systems}},
	url = {http://arxiv.org/abs/2503.13675},
	doi = {10.48550/arXiv.2503.13675},
	urldate = {2025-03-21},
	publisher = {arXiv},
	author = {Shrivastava, Shaurya and Oguri, Kenshiro},
	month = mar,
	year = {2025},
	note = {arXiv:2503.13675 [math]},
}

@article{diamondCVXPYPythonEmbeddedModeling2016,
	title = {{CVXPY}: {A} {Python}-{Embedded} {Modeling} {Language} for {Convex} {Optimization}},
	volume = {17},
	issn = {1533-7928},
	shorttitle = {{CVXPY}},
	url = {http://jmlr.org/papers/v17/15-408.html},
	number = {83},
	urldate = {2025-03-13},
	journal = {J. of Mach. Learn. Res.},
	author = {Diamond, Steven and Boyd, Stephen},
	year = {2016},
	pages = {1--5},
}

@inproceedings{minkaExpectationPropagationApproximate2001,
	address = {San Francisco, CA, USA},
	series = {{UAI}'01},
	title = {Expectation propagation for approximate {Bayesian} inference},
	isbn = {978-1-55860-800-9},
	urldate = {2025-03-14},
	booktitle = {Proc. of the {Seventeenth} conf. on {Uncertainty} in {Artif}. {Intell}.},
	author = {Minka, Thomas P.},
	month = aug,
	year = {2001},
	pages = {362--369},
}

@incollection{starekSpacecraftAutonomyChallenges2016,
	address = {Berlin, Heidelberg},
	title = {Spacecraft {Autonomy} {Challenges} for {Next}-{Generation} {Space} {Missions}},
	volume = {460},
	isbn = {978-3-662-47693-2 978-3-662-47694-9},
	language = {en},
	urldate = {2025-01-08},
	booktitle = {Advances in {Control} {System} {Technology} for {Aerospace} {Applications}},
	publisher = {Springer Berlin Heidelberg},
	author = {Starek, Joseph A. and Açıkmeşe, Behçet and Nesnas, Issa A. and Pavone, Marco},
	editor = {Feron, Eric},
	year = {2016},
	pages = {1--48},
}

@article{spreenTrajectoryDesignTargeting2021,
	title = {Trajectory {Design} and {Targeting} for {Applications} to the {Exploration} {Program} in {Cislunar} {Space}},
	url = {https://docs.lib.purdue.edu/dissertations/AAI30504506},
	journal = {Theses and Dissertations Available from ProQuest},
	author = {Spreen, Emily M. Zimovan},
	month = jan,
	year = {2021},
	pages = {1--249},
}

@inproceedings{kumagaiSequentialChanceConstrainedCovariance2024,
	title = {Sequential {Chance}-{Constrained} {Covariance} {Steering} {For} {Robust} {Cislunar} {Trajectory} {Design} {Under} {Uncertainties}},
	booktitle = {{AAS}/{AIAA} {Astrodynamics} {Specialist} {Conference}},
	author = {Kumagai, Naoya and Oguri, Kenshiro},
	year = {2024},
}

@article{connorhowellThreedimensionalPeriodicHalo1984,
	title = {Three-dimensional, periodic, ‘halo’ orbits},
	volume = {32},
	issn = {1572-9478},
	url = {https://doi.org/10.1007/BF01358403},
	doi = {10.1007/BF01358403},
	language = {en},
	number = {1},
	urldate = {2024-12-23},
	journal = {Celestial mechanics},
	author = {Connor Howell, Kathleen},
	month = jan,
	year = {1984},
	pages = {53--71},
}

@inproceedings{shrivastavaFaultTolerantFeedbackControl2023,
	title = {Fault-{Tolerant} {Feedback} {Control} for {Spacecraft} {Rendezvous} using {Markov} {Jump} {Model}},
	booktitle = {{AAS}/{AIAA} {Astrodynamics} {Specialist} {Conference}},
	author = {Shrivastava, Shaurya and Oguri, Kenshiro},
	month = aug,
	year = {2023},
}

@article{fuModelReferenceTracking2017,
	title = {Model reference tracking control of continuous-time periodic linear systems with actuator jumping fault and its applications in orbit maneuvering},
	volume = {15},
	issn = {2005-4092},
	url = {https://doi.org/10.1007/s12555-016-0104-9},
	doi = {10.1007/s12555-016-0104-9},
	language = {en},
	number = {5},
	urldate = {2024-09-18},
	journal = {International Journal of Control, Automation and Systems},
	author = {Fu, Yan-Ming and Lu, Yang and Zhang, Mao-Rui},
	month = oct,
	year = {2017},
	pages = {2182--2192},
}

@inproceedings{mccartyMissedThrustAnalysis2020,
	title = {Missed {Thrust} {Analysis} and {Design} {For} {Low} {Thrust} {Cislunar} {Transfers}},
	booktitle = {{AAS}/{AIAA} {Astrodynamics} {Specialist} {Conference}},
	author = {McCarty, Steven and Grebow, Daniel},
	month = aug,
	year = {2020},
}

@misc{goulartClarabelInteriorpointSolver2024,
	title = {Clarabel: {An} interior-point solver for conic programs with quadratic objectives},
	shorttitle = {Clarabel},
	url = {http://arxiv.org/abs/2405.12762},
	doi = {10.48550/arXiv.2405.12762},
	urldate = {2024-09-10},
	publisher = {arXiv},
	author = {Goulart, Paul J. and Chen, Yuwen},
	month = may,
	year = {2024},
	note = {arXiv:2405.12762},
}

@article{venigallaMultiObjectiveLowThrustTrajectory2022,
	title = {Multi-{Objective} {Low}-{Thrust} {Trajectory} {Optimization} with {Robustness} to {Missed} {Thrust} {Events}},
	volume = {45},
	issn = {0731-5090},
	url = {https://arc.aiaa.org/doi/10.2514/1.G006056},
	doi = {10.2514/1.G006056},
	number = {7},
	urldate = {2023-04-07},
	journal = {Journal of Guidance, Control, and Dynamics},
	publisher = {American Institute of Aeronautics and Astronautics},
	author = {Venigalla, Chandrakanth and Englander, Jacob A. and Scheeres, Daniel J.},
	month = jul,
	year = {2022},
	pages = {1255--1268},
}

@article{grecoRobustSpaceTrajectory2022,
	title = {Robust {Space} {Trajectory} {Design} {Using} {Belief} {Optimal} {Control}},
	copyright = {Copyright © 2022 by the authors. Published by the American Institute of Aeronautics and Astronautics, Inc., with permission. All requests for copying and permission to reprint should be submitted to CCC at www.copyright.com; employ the eISSN 1533-3884 to initiate your request. See also AIAA Rights and Permissions www.aiaa.org/randp.},
	url = {https://arc.aiaa.org/doi/10.2514/1.G005704},
	doi = {10.2514/1.G005704},
	language = {en},
	urldate = {2022-09-19},
	journal = {Journal of Guidance, Control, and Dynamics},
	publisher = {American Institute of Aeronautics and Astronautics},
	author = {Greco, Cristian and Campagnola, Stefano and Vasile, Massimiliano},
	month = apr,
	year = {2022},
}

\end{document}